\documentclass[11pt]{amsart}

\usepackage{graphicx}
\usepackage{amssymb,amsmath, amsthm} 
\usepackage{amsfonts}
\usepackage{mathrsfs}
\usepackage{enumerate}
\usepackage{esint}
\usepackage[svgnames]{xcolor} 
\usepackage[colorlinks,citecolor=red,pagebackref,hypertexnames=false,breaklinks]{hyperref}
\usepackage{pgf,tikz}
\usepackage{pdfsync}
\usepackage[normalem]{ulem}

\usepackage[top = 4 cm, bottom=4cm, left=3.5cm, right=3.5cm]{geometry} 

\numberwithin{equation}{section}

\theoremstyle{plain}
\newtheorem{theorem}[equation]{Theorem}
\newtheorem{lemma}[equation]{Lemma}
\newtheorem{corollary}[equation]{Corollary}
\newtheorem{proposition}[equation]{Proposition}

\theoremstyle{definition}
\newtheorem{definition}[equation]{Definition}

\theoremstyle{remark}
\newtheorem{remark}[equation]{Remark}

\newcommand{\ep}{\varepsilon}
\newcommand{\ra}{\rightarrow}
\newcommand{\n}[1]{\mathscr{#1}}
\newcommand{\m}[1]{\mathcal{#1}}
\newcommand{\bb}[1]{\mathbb{#1}}
\newcommand{\mbf}[1]{\mathbf{#1}}
\newcommand{\f}[1]{\mathfrak{#1}}

\def\XXint#1#2#3{{\setbox0=\hbox{$#1{#2#3}{%
\int}$ }
\vcenter{\hbox{$#2#3$ }}\kern-.6\wd0}}

\usepackage{pgf}
\usepackage{color}

\calclayout
\allowdisplaybreaks

\newcommand{\supp}{\operatorname{supp}}
\newcommand{\dist}{\operatorname{dist}}

\newcommand{\rn}{{\mathbb{R}^n}}

\newcommand{\RR}{{\mathbb{R}}}

\newcommand{\D}{D_{\mathbf a}}

\newcommand{\bp}{\noindent {\it Proof}.\,\,}

\DeclareMathOperator{\diam}{diam}

\usepackage{color}
\definecolor{mygreen}{rgb}{0,0.5,0}
\definecolor{purple}{rgb}{0.5,0,0.5}

\begin{document}

\author[S. Mayboroda]{Svitlana Mayboroda}
\email[]{svitlana@math.umn.edu}
\address{School of Mathematics, University of Minnesota, Minneapolis, MN, USA}

\author[B. Poggi]{Bruno Poggi}
\email[]{poggi008@umn.edu}
\address{School of Mathematics, University of Minnesota, Minneapolis, MN, USA}

\thanks {
	The first author is supported in part by the NSF INSPIRE Award DMS 1344235, NSF CAREER Award DMS 1220089, Simons Fellowship, and the Simons Foundation grant 563916, SM. Both authors would like to thank the Mathematical Sciences Research Institute (NSF grant DMS 1440140) for support and hospitality.}

\title[Exponential decay estimates]{Exponential decay estimates for fundamental solutions of Schr\"odinger-type operators}
\date{February 3, 2019.}

\begin{abstract} In the present paper we establish sharp exponential decay estimates for operator and integral kernels of   the (not necessarily self-adjoint) operators $L=-(\nabla-i\mathbf{a})^TA(\nabla-i\mathbf{a})+V$. The latter class includes, in particular,  the magnetic Schr\"odinger operator $-\left(\nabla-i\mathbf{a}\right)^2+V$ and the generalized electric Schr\"odinger operator $-{\rm div }A\nabla+V$. Our exponential decay bounds rest on a generalization of the Fefferman-Phong uncertainty principle to the present context and are governed by the Agmon distance associated to the corresponding maximal function. In the presence of a scale-invariant Harnack inequality, for instance, for the generalized electric Schr\"odinger operator with real coefficients, we establish both lower and upper estimates for fundamental solutions, thus demonstrating sharpness of our results. The only previously known estimates of this type pertain to the classical Schr\"odinger operator $-\Delta +V$ \cite{S1}  .
\end{abstract}

\maketitle

\tableofcontents

\section{Introduction}
The exponential decay of solutions to the Schr\"odinger operator in the presence of a positive potential is an important property underpinning foundation of quantum physics. However, establishing a precise rate of decay for complicated potentials is a challenging open problem to this date. For instance, the Landis conjecture from around 1965 claimed that bounded solutions of $-\Delta u+Vu =0$, with a bounded potential $V$, should decay no faster than $e^{-C|x|^{1+}}$. In the context of complex potentials, this conjecture was disproved by Meshkov \cite{M}, and later on the sharp lower bound on decay rate,  $e^{-C|x|^{4/3}\log |x|}$, was proved by notoriously difficult Carleman-type arguments in \cite{BK}. The latter, however, does not distinguish between real and complex potentials, and it is not known to this date if Landis conjecture is true for real-valued $u$ and $V$. Notably, none of these results are carefully adapted to the local properties of the potential which drives the decay. 

In the present paper we tackle the problem of optimal decay rate from a different angle. We use a generalized version of the Fefferman-Phong uncertainty principle to show that the fundamental solution of a generalized Schr\"odinger operator with a non-negative potential enjoys the same bounds from above and below. For instance,  if $A$ is an elliptic matrix with real, bounded coefficients, and $V\in RH_{\frac n2}$, then 
\begin{equation}\label{lowerbounded2-intro}
\frac{c_1e^{-\ep_1d(x,y,V)}}{|x-y|^{n-2}}\leq \Gamma_E(x,y)\leq \frac{c_2e^{-\ep_2d(x,y,V)}}{|x-y|^{n-2}}, 
\end{equation}
where $\Gamma_E$ is an integral kernel of the generalized electric Schr\"odinger operator, that is, the fundamental solution to $L_E=-{\rm div }A\nabla+V$, $x, y\in \RR^n$, interpreted in a suitable weak sense, and $d$ is a certain distance function depending on $V$. In fact, we establish the upper estimates for a considerably more general class of operators, which can be formally written as $L=-(\nabla-i\mathbf{a})^TA(\nabla-i\mathbf{a})+V$ including, in particular, magnetic Schr\"odinger. Let us discuss this in some details and properly define all the notions used in this statement. 

The first results expressing upper estimates on the solutions in terms of a certain distance associated to the potential $V$ go back to Agmon \cite{Agmon}. He has introduced a distance function which now bears his name and which we will discuss below, and showed that the solution decays exponentially in the region where $V\geq 0$. Agmon's estimates, however, are clearly non-sharp for most non-trivial potentials, for a simple reason that solutions carry some amount of regularity and low values of $V$ in small regions should not drastically affect their decay properties. This vague statement is very hard to quantify (and there are other interfering mechanisms at play which we will not discuss in this paper). In many situations, however, the behavior of solutions, notably of the eigenfunctions, is rather precisely governed by the uncertainty principle. The latter has a few manifestations. In particular, the most recent one in \cite{ADFJM} yielded astonishingly accurate estimates on eigenfunctions even for the prototype of the Anderson model based on disordered potentials. Here, however, we will use a much earlier result due to Fefferman and Phong which has been later extended, e.g., to the context of the magnetic Schr\"odinger operator - one of our main objects of study in the present paper. This extension unfortunately currently seems beyond the reach of the methods in \cite{ADFJM}.

Let $A=(A_{ij})_{i,j=1}^n$ be an $n\times n$ matrix with complex bounded measurable coefficients satisfying the uniform ellipticity conditions
\begin{equation}
\lambda|\xi|^{2}\leq\Re e\,\langle A(x)\xi,\xi\rangle
\equiv \Re e\sum_{i,j=1}^{n}A_{ij}(x)\xi_{j}\bar{\xi_{i}}\,\,\,\text{and}\quad
\Vert A\Vert_{L^{\infty}(\mathbb{R}^{n})}\leq\Lambda,
\label{eq1.2-intro}
\end{equation}
for some $\lambda>0$, $\Lambda<\infty$ and for all $\xi\in\mathbb{C}^{n}$, $x\in\mathbb{R}^{n}$. Let $\mathbf{a}=(a_1,\ldots,a_n)$ be a vector of real-valued $L^2_{loc}(\RR^n)$ functions and assume that $V\in L^1_{loc}(\RR^n)$ is scalar, real valued, and positive almost everywhere on $\RR^n$.  Following the tradition and physical significance, we will refer to ${\bf a}$ as the magnetic potential and to $V$ as the electric potential. We consider the operator formally given by 
\begin{equation}\label{generalized}
L=-(\nabla-i\mathbf{a})^TA(\nabla-i\mathbf{a})+V, 
\end{equation}
where the superscript $T$ stands for the transpose of the corresponding vector. Let us further denote
\[
\D=\nabla-i\mathbf{a},
\]
and the magnetic field by  $\mathbf{B}$, so that
\begin{equation}\label{magfield}
\mbf B=\text{ curl }\mbf{a}=(b_{jk})_{1\leq j,k\leq n}
\end{equation}
where
\begin{equation}\label{curlterms}
b_{jk}(x)=\frac{\partial a_j}{\partial x_k}-\frac{\partial a_k}{\partial x_j}.
\end{equation}
Due to the gauge invariance property, one expects ${\bf B}$ rather than ${\bf a}$ to be the primary relevant parameter. As mentioned in the abstract, the most important particular cases that will be highlighted throughout the paper are the magnetic Schr\"odinger operator $-\left(\nabla-i\mathbf{a}\right)^2+V$ and the generalized electric Schr\"odinger operator $-{\rm div }A\nabla+V$. We remark that our operators are not necessarily self-adjoint, in particular, the matrix $A$ is not required to be symmetric (or even real-valued in the first part of the paper).

Our goal is to treat as general a situation as possible taking no  regularity assumptions on $A$ or on $V$. The Fefferman-Phong Uncertainty Principle requires a mild control on oscillations of $V$ and ${\bf B}$, manifested, for instance, in terms of a membership to suitable weight spaces.   
We say that $w\in L^p_{loc}(\bb R^n)$, with $w>0$ a.e., belongs to the Reverse H\"older class
$RH_p=RH_p(\bb R^n)$ if there exists a constant $C$ so that for any ball $B\subset\bb R^n$,
\begin{equation}\label{reverseholder-intro}
\left(\fint_B w^p\right)^{1/p}\leq C\fint_B w,
\end{equation}
and the reader will witness below that we typically assume that $V+|\mbf B|\in RH_{n/2}$. Obviously, such potentials are not necessarily smooth and not necessarily bounded. We will review below some key historical points but let us briefly mention here that a lot of attention has been devoted to certain polynomial potentials - in particular, because they serve as a toy model in related problems in semiclassical analysis. Any non-negative polynomial belongs to $RH_{n/2}$ class, with the constant depending on the degree and the dimension. In fact, for any polynomial $P$ and $\alpha>0$ we have $|P|^\alpha\in RH_p$ for any $p>1$, with the constant depending only on $\alpha, n,$ and the degree of $P$, and $|x|^{\alpha}\in RH_{n/2}$ for any $\alpha>-2$.

For a function $w\in RH_p, p\geq\frac n2$, define the maximal function $m(x,w)$ by
\begin{equation}\label{maximal-intro}
\frac1{m(x, w)}:=\sup_{r>0}\Bigl\{r:\, \frac1{r^{n-2}}\int_{B(x,r)} w\leq1\Bigr\},
\end{equation}
and the distance function
\begin{equation}\label{distfn}
d(x,y, w)=\inf_{\gamma}\int\limits_0^1m(\gamma(t), w)|\gamma'(t)|\,dt,
\end{equation}
where $\gamma:[0,1]\ra\bb R^n$ is absolutely continuous and $\gamma(0)=x, \gamma(1)=y$. Finally, for any  $U\subset\bb R^n$, we let $d(x,U, w):=\inf\limits_{y\in U}d(x,y, w).$ The function $m$ comes from the Uncertainty Principle, which is generalized in the present paper to serve the operators \eqref{generalized}. Its explicit formula (\ref{maximal-intro}) was introduced in \cite{S4}. Since it is one of the main points underpinning many of our results, let us say a few more words. The function $m$ measures the sum of the contributions of the kinetic energy $ \Re e A \D f \overline{\D f}$ and potential energy $V|f|^2$, reaching optimum when $f$ is a bump. In this vein, we find the definition \eqref{maximal-intro} more telling than any of the particular representations, but let us mention nonetheless that in the aforementioned case $V=|P|^\alpha$, $\alpha>0$, we have $m(\cdot, V)\approx\sum_{|\beta|\leq k} |\partial^\beta P|^{\frac{\alpha}{\alpha|\beta|+2}},$ where $k$ is the degree of $P$ \cite{S3}.  In the context of polynomial-like potentials, by methods crucially relying on smoothness, the Uncertainty Principle has been proved in \cite{Sm, HM, HN, MN}. However, the real breakthrough in this direction came when Fefferman and Phong treated the non-smooth potentials \cite{F}. This approach has been further formalized in  \cite{S4} to address $V\in RH_{n/2}$ and in \cite{S3}, \cite{benali}, to treat magnetic potentials.


With this notation in mind, we list our main results. 

For any operator $L$ given by \eqref{generalized} and for any $f\in L^2(\RR^n)$ with compact support, there exist constants $\tilde d,\ep, C>0$ such that
\begin{multline}
\int_{\big\{x\in\bb  R^n|d(x,\supp f,V+|\mathbf{B}|)\geq\tilde d\big\}}m\left(\cdot,V+|\mathbf{B}|\right)^2|L^{-1}f|^2e^{2\ep d(\cdot,\supp f,V+|\mathbf{B}|)}\\\leq C\int_{\bb R^n}|f|^2\frac1{m(x,V+|\mathbf{B}|)^2},\label{l2resulthom-intro}
\end{multline}
provided that $A$ is an elliptic matrix with complex bounded measurable coefficients, and either ${\bf a}=0$ and $V\in RH_{n/2}$, or, more generally, ${\mathbf a}\in L^2_{loc}(\bb R^n)$, $V>0$ a.e. on $\bb R^n$, and
\begin{equation}\label{magassumptions-intro}
\left\{\begin{matrix}V+|\mbf B|\in RH_{n/2},\\0\leq V\leq c\,m(\cdot,V+|\mbf B|)^2,\\|\nabla \mbf B|\leq c'\,m(\cdot,V+|\mbf B|)^3.\end{matrix}\right.
\end{equation}
An analogous estimate holds  for the resolvent operator $(I+t^2L)^{-1}$, $t>0$: 
\begin{eqnarray}\nonumber
\int\limits_{\big\{x\in\bb  R^n|d(x,\supp f,V+|\mathbf{B}|+\frac {1}{t^2})\geq\tilde d\big\}}m\Big(\cdot,V+|\mathbf{B}|+\frac {1}{t^2}\Big)^2\left|(I+t^2L)^{-1}f\right|^2e^{2\ep d(\cdot,\supp f,V+|\mathbf{B}|+\frac {1}{t^2})}\\[4mm]\leq C\int\limits_{\bb R^n}|f|^2m\Big(\cdot,V+|\mathbf{B}|+\frac {1}{t^2}\Big)^2.\label{l2result-intro}
\end{eqnarray}

In other words, $L^{-1} f$ decays as $e^{-\ep d(\cdot,\supp f,V+|\mathbf{B}|)}$ away from the support of $f$ and the resolvent decays as $e^{-\ep d(\cdot,\supp f,V+|\mathbf{B}|+ \frac {1}{t^2})}$. The strongest previously known result for the resolvent (almost) in this generality is due to Germinet and Klein \cite{GK}. Their work is restricted to self-adjoint operators, but otherwise, modulo some technical differences, they treat considerably more general elliptic systems than we do, including the Maxwell equation, and they go much farther towards the Combes-Thomas estimates. However, the exponential decay that they postulate is a much weaker estimate with $\frac {1}{t^2}$ in place of our $V+|\mathbf{B}|+ \frac {1}{t^2}$. They do not treat the operator $L^{-1}$. Actually, an estimate with $\frac {1}{t^2}$ in place of our $V+|\mathbf{B}|+ \frac {1}{t^2}$ has also appeared in many sources before but under stronger assumptions on the operator, and we do not attempt to review the corresponding literature. 

%

Due to a possible lack of local boundedness of solutions to \eqref{generalized}, the $L^2$ estimates in \eqref{l2resulthom-intro} are of the nature of the best possible. However, for operators whose solutions satisfy Moser and/or Harnack inequality stronger pointwise bounds can be obtained. For instance, if $\mathbf{a}\in L^2_{loc}(\bb R^n)$, and assumptions (\ref{magassumptions-intro}) are satisfied, then 
\begin{equation}\label{upperbound1-intro}
|\Gamma_{M}(x,y)|\leq\frac{Ce^{-\ep d(x,y,V+|\mathbf{B}|)}}{|x-y|^{n-2}} \quad \mbox{ for a.e. }x, y \in \rn,
\end{equation}
where $\Gamma_M$ is an integral kernel of the magnetic Schr\"odinger operator $L_M=-\left(\nabla-i\mathbf{a}\right)^2+V$, that is, the solution to $L_M \Gamma (x, y)=\delta_y(x)$, $x, y\in \RR^n$, interpreted in a suitable weak sense. More generally, this bound is valid for any operator \eqref{generalized} with an elliptic matrix $A$ of complex bounded measurable coefficients and \eqref{magassumptions-intro}, assuming, in addition, local boundedness of solutions and a classical estimate on the fundamental solution by $|x-y|^{2-n}$ -- see Theorem \ref{Upper Bound}.

Finally, if fundamental solutions are bounded from above and below by a multiple of $|x-y|^{2-n}$, e.g., if 
$A$ is a real, bounded, elliptic matrix, $V\in RH_{\frac n2}$, and ${\bf a}=0$, then we establish both upper and lower estimates \eqref{lowerbounded2-intro}  -- see Corollaries \ref{upperboundcor} and \ref{lowerboundcor}. This covers the case of the generalized electric Schr\"odinger operator. 

Before passing to the proof of these results, let us make a few more points about the existing literature. Needless to say, \eqref{lowerbounded2-intro} underlines  sharpness of the emerging estimates. The only context in which \eqref{lowerbounded2-intro} have been proved before is that of the classical Schr\"odinger operator $-\Delta +V$ \cite{S1}, and we, of course, build on the ideas from \cite{S1}. 
As we mentioned in the beginning, to the best of our knowledge, no sharp results on the exponential decay of the kernels to the magnetic Schr\"odinger operator or generalized Schr\"odinger operator existed in the literature. In fact, even the existence of the fundamental solution to the magnetic Schr\"odinger operator for ${\bf a}\in L^2_{loc}(\RR^n)$ and $V\in L^1_{loc}(\RR^n)$, subject to the usual bound by a multiple of $|x-y|^{2-n}$ (Theorem~\ref{fundsolapprox}), seemed to be out of reach, as previous treatises normally imposed somewhat ad hoc conditions of smoothness for the magnetic field ${\bf a}\in C^2$ or at least $\mathbf{a}\in L^4_{loc}(\bb R^n), \text{div }\mathbf{a}\in L^2_{loc}(\bb R^n)$, and $V\in L^{\infty}_{loc}(\bb R^n)$ \cite{KS, benali}. As we will see in Section \ref{theory}, both situations are considerably simpler than ours but not completely natural, for ${\bf a}\in L^2_{loc}(\RR^n)$ and $V\in L^1_{loc}(\RR^n)$ are the minimal restrictions allowing one to make sense of the bilinear form associated to $L_M$ in the weak sense. 

Furthermore, certain polynomial upper estimates on the fundamental solutions in terms of $m$ have been established in a variety of contexts, in particular, in \cite{S4} for $-\Delta+V$, and in \cite{S5} for the magnetic Schr\"odinger operator, under assumptions similar to ours. Polynomial decay is sufficient for establishing key properties of the associated Riesz transforms and similar operators - the main goal of the majority of these papers, but is obviously not sharp. An attempt to get  exponential decay has been made at \cite{K}. In this paper the author treated the heat kernel estimates for $-{\rm div} A\nabla+V$ ($A$ real and symmetric) and $-\left(\nabla-i\mathbf{a}\right)^2+V$ and integrating them obtained for these two operators 
$$|\Gamma(x,y)|\leq\frac{Ce^{-\ep (1+m(x,V+|\mathbf{B}|)|x-y|)^{\frac{2}{2k_0+3}}}}{|x-y|^{n-2}} \quad \mbox{ for a.e. }x, y \in \rn,$$
for some $k_0>0$, which is, once again, not sharp, as can be seen from \eqref{lowerbounded2-intro} and \eqref{upperbound1-intro}. Finally, without an attempt of a comprehensive review of the theory, we mention that resolvent estimates are routinely used in many aspects of semiclassical analysis, which roughly speaking, concentrates on the behavior of solutions to $-\hbar^2\Delta +V$ and analogous operators as $\hbar\to 0$, but these bounds are typically independent of local features of $V$ and ${\bf B}$. The major achievement of the present paper is similar estimates from above and below which is only possible by a careful account of the impact of the electric and magnetic potentials. 


The outline of the paper is as follows. In Section~\ref{theory}, we present a theory of a generalized magnetic Schr\"odinger operator $L$ from \eqref{generalized}, we define the resolvent, the heat semigroup, the inverse of $L$, and other notions. In Section~\ref{mproperties}, we provide auxiliary estimates on the maximal function $m$ and distance $d$; most of the material in this section is well-known. In Section~\ref{l2sec}, we establish exponential decay in $L^2$ for the resolvent of the operator $L$ and for $L^{-1}$, including \eqref{l2resulthom-intro} and \eqref{l2result-intro}. In Section \ref{fundsolnsec}, we provide a construction of the fundamental solution to the magnetic Schr\"odinger operator with ${\bf a}\in L^2_{loc}(\RR^n)$ and $V\in L^1_{loc}(\RR^n)$, $V>0$ a.e., together with the basic bound by $C|x-y|^{2-n}$. In Section~\ref{upperboundsec}, we establish exponential upper pointwise bounds on the fundamental solution, including, in particular, \eqref{upperbound1-intro} and the upper bound in  \eqref{lowerbounded2-intro}.  In Section~\ref{lowerboundsec}, we give exponential lower bounds for the fundamental solution, including, in particular, the lower bound in  \eqref{lowerbounded2-intro}. 

\section{The theory of the generalized magnetic Schr\"odinger operator}\label{theory}

\subsection{Preliminaries}
{  We always assume $n\geq3$. Let $L$ be the operator formally given by (\ref{generalized}), where $\mathbf a\in L^2_{loc}(\bb R^n)$ is a real-valued vector function, $A$ is an elliptic matrix with complex, bounded measurable coefficients, $V\in L^1_{loc}(\bb R^n)$ is scalar, complex. We will write $\Re e\, V=V^+-V^-$ where $V^{\pm}\geq0$. The negative part, $V^-$, must satisfy
\begin{equation}\label{gkcond}
\int_{\bb R^n}V^-|u|^2\leq c_1\int_{\bb R^n}\Re e\, A\D u\overline{\D u}\,+\,c_2\Vert u\Vert_{L^2(\bb R^n)}^2,\qquad\text{for each } u\in C_c^{\infty}(\bb R^n),
\end{equation}
where $c_1\in[0,1)$ and $c_2\in[0,\infty)$.

Let $H:=L^2(\bb R^n)$. Corresponding to the operator $L$ given in (\ref{generalized}), we consider for $u,v\in C_c^{\infty}(\bb R^n)$ the form
\begin{equation}\label{lform}
\f l(u,v):=\int_{\bb R^n}A\D u\overline{\D v}+Vu\overline v.
\end{equation}
and then define the domain of $\f l$ as the completion of $C_c^{\infty}(\bb R^n)$ with respect to the norm
\begin{equation}\label{lnorm}
\Vert u\Vert_{\f l}:=\sqrt{\Re e\,\f l(u,u)+(1+c_2)\Vert u\Vert^2_H},
\end{equation}
which will henceforth be known as $D(\f l)$. This can be done because by adding $0\leq c_1\int_{\bb R^n}[\Re e\,V+V^-]|u|^2$ to (\ref{gkcond}), we see that
\begin{equation}\label{negdombyform}
\int_{\bb R^n}V^-|u|^2\leq\frac{c_1}{1-c_1}\Re e\,\f l(u,u)\,+\,\frac{c_2}{1-c_1}\Vert u\Vert_H^2,\qquad\text{for each } u\in D(\f l),
\end{equation} 
and so $\Re e\,\f l(u,u)+c_2\Vert u\Vert^2_H\geq0$. It immediately follows that
\[
\Vert u\Vert_{\f l}=0~\implies~u=0\text{ a.e. on }\bb R^n,
\]
and $(D(\f l),\Vert\cdot\Vert_{\f l})$ is a normed space. From (\ref{negdombyform}) we can see that
\begin{equation}\label{negdombynorm}
\int_{\bb R^n}V^-|u|^2\leq\frac{1}{1-c_1}\Vert u\Vert_{\f l}^2,\qquad\text{for each } u\in D(\f l).
\end{equation}
Moreover, using (\ref{gkcond}) and the fact that $\Re e\,V=V^+-V^-$, it is easy to conclude that
\begin{equation}\label{posdombynorm}
\int_{\bb R^n}V^+|u|^2\leq\Vert u\Vert_{\f l}^2,\qquad\text{for each } u\in D(\f l).
\end{equation}
We will also need to consider the following condition on the imaginary part of $V$: 
\begin{equation}\label{imVreq}
\int_{\bb R^n}|\text{Im }V|\,|u||v|\leq c_3\sqrt{\Re e\,\f l(u,u)+c_4\Vert u\Vert_H^2}~\sqrt{\Re e\,\f l(v,v)+c_4\Vert v\Vert_H^2},~\text{for each } u,v\in D(\f l),
\end{equation}
where $c_3>0$ and $c_4$ is either $0$ or $1$. Of course, the condition (\ref{imVreq}) with $c_4=0$ implies the one with $c_4=1$, but we'll see that for the non-homogeneous setting, the case $c_4=1$ is enough for us. As a start, we recall the \emph{diamagnetic inequality}, which can be formulated as

\begin{lemma}\label{diamagnetictheorem} Suppose that $\mathbf{a}$ is a measurable function on $\bb R^n$ and let $u$ be a measurable function on $\bb R^n$ such that $\D u$ is a measurable function on $\bb R^n$. Then $\nabla|u|$ is a measurable function on $\bb R^n$, and
\begin{equation}\label{diamagnetic}
\Big|\nabla|u|(x)\Big|\leq\Big|\D u(x)\Big|,
\end{equation}
for almost every $x\in\bb R^n$.
\end{lemma}


Let $2^*:=\frac{2n}{n-2}$. By the Sobolev Embedding, we have that
\[
\Vert u\Vert_{L^{2^*}(\bb R^n)}\leq C\Vert\nabla u\Vert_{L^2(\bb R^n)}.
\]
Observe that the map $\Vert\nabla\cdot\Vert_{L^2(\bb R^n)}$ is a norm on $C_c^{\infty}(\bb R^n)$. Define $Y^{1,2}$ as the completion of $C_c^{\infty}(\bb R^n)$ under this norm. The diamagnetic inequality implies
\begin{corollary}\label{coolembedding} Suppose that $\mathbf{a}\in L^2_{loc}(\bb R^n)$ and $u\in C_c^{\infty}(\bb R^n)$. Then
\[
\Vert u\Vert_{L^{2^*}(\bb R^n)}\leq C\Vert\D u\Vert_{L^2(\bb R^n)}.
\]
In particular, $D(\f l)\hookrightarrow L^{2^*}(\bb R^n)$.
\end{corollary}

The diamagnetic inequality is especially useful for the aforementioned embedding. At the moment, we turn to the theory of the form $\f l$: the proof of the following proposition is standard:

\begin{proposition}\label{formlemma} Let $\mathbf{a}\in L^2_{loc}(\bb R^n)$, let $A$ be an elliptic matrix with complex, bounded measurable coefficients, and let $V\in L^1_{loc}(\bb R^n)$ satisfying (\ref{gkcond}) and (\ref{imVreq}) with $c_4$ either $0$ or $1$. Then the form $\f l$ is densely defined, bounded below, continuous, and closed. If $c_2\equiv0$ in (\ref{gkcond}), then $\f l$ is accretive.
\end{proposition}

Unless stated otherwise, take all the assumptions of the previous proposition, with $c_2\equiv0$ in (\ref{gkcond}). Using Definition 1.21 of \cite{Ou} (with $A\mapsto L$, $\f a\mapsto\f l$ in the notation of \cite{Ou}), we can define an unbounded operator $L:D(L)\ra H$ where $D(L)$ is given as
\[
D(L)=\Big\{u\in D(\f l) \text{ s.t. }\exists v\in H: \f l(u,\phi)=(v,\phi)_H~\forall\phi\in D(\f l)\Big\},\qquad Lu:=v.
\]
The operator $L$ is called the \emph{operator associated with the form $\f l$}. Then Proposition 1.22 in \cite{Ou} applies and we conclude that $L$ is densely defined, for every $\ep>0$ the operator $L+\ep$ is invertible from $D(L)$ into $H$, and its inverse $(L+\ep)^{-1}$ is a bounded operator on $H$. In addition,
\[
\Vert\ep(L+\ep)^{-1}f\Vert_H\leq\Vert f\Vert_H,\qquad\text{for each }\ep>0, f\in H.
\]
We will denote by $\f l^*$ the adjoint form of $\f l$, and by $L^*$ the operator associated to $\f l^*$ (see Proposition 1.24 in \cite{Ou}). We also note (Lemma 1.25 in \cite{Ou}) that $D(L)$ is a core of $\f l$; that is, $D(L)$ is dense in $D(\f l)$ under the norm $\Vert\cdot\Vert_{\f l}$. Moreover, since by the aforementioned results, the resolvent set $\rho(-L)$ is not empty, then by Proposition 1.35 in \cite{Ou}, we see that $-L$ is a closed operator.

We now see that $L$ is an accretive operator (see Definition 1.46 in \cite{Ou}) since $\f l$ is an accretive form. Since $(L+\ep)$ is invertible, then in particular it has dense range. So, by Theorem 1.49 in \cite{Ou}, it follows that $L$ is m-accretive, and that $-L$ is the generator of a strongly continuous contraction semigroup on $H$.


\subsection{The homogeneous operator}

So far we have seen that the expression $(L+\ep)^{-1}f$ makes sense for $\ep>0$, but our previous construction cannot work for the homogeneous case: the operator $L$ as defined above is not necessarily invertible as a map from $D(L)$ to $H=L^2(\bb R^n)$. It is imperative therefore to construct a homogeneous theory. The following argument is inspired by that of Section $3$ of \cite{AB}. 

Let us use the notation $\m V_{\mathbf{a},V}:=D(\f l)$, where we will omit the subscript if the magnetic and electric potentials are clear from context. Observe that by the diamagnetic inequality we have that $(\Re e\int A\D u\overline{\D u})^{\frac12}$ is a norm on $C_c^{\infty}(\bb R^n)$.  If $c_2\equiv0$, we define the space $\dot{\m V}$ as the completion of $C_c^{\infty}(\bb R^n)$ under the norm
\[
\Vert u\Vert_{\dot{\m V}}:=\sqrt{\Re e\,\Big(\int_{\bb R^n}A\D u\overline{\D u}+V|u|^2\Big)}.
\]
Indeed, if $u\in\dot{\m V}$ with $\Vert u\Vert_{\dot{\m V}}=0$, then by (\ref{gkcond}), the diamagnetic inequality and the Sobolev embedding we obtain
\[
\Vert u\Vert_{L^{\frac{2n}{n-2}}(\bb R^n)}\lesssim\Vert\D u\Vert_{L^2(\bb R^n)}\lesssim\Vert u\Vert_{\dot{\m V}}=0,
\]
whence we must have $u=0$ a.e. on $\bb R^n$. Thus we have $\dot{\m V}\hookrightarrow L^{2^*}(\bb R^n)$. For instance, if $\mathbf{a}\equiv0, V\equiv0$, then $\dot{\m V}=Y^{1,2}$. The form $\dot{\f l}$ is given by the same formula as $\f l$ in (\ref{lform}) for $u,v\in\dot{\m V}$, and $\dot{\f l}$ is a coercive, bounded form on $\dot{\m V}$. Now also suppose that $c_4\equiv0$. The estimates
\[
\Vert\D u\Vert_{L^2(\bb R^n)}^2\leq C(\lambda,c_1)\Vert u\Vert_{\dot{\m V}}^2,
\]
\[
\int_{\bb R^n}|V||u|^2\leq C(c_1,c_3)\Vert u\Vert_{\dot{\m V}}^2,
\]
\[
\Vert u\Vert_{\dot{\m V}}\leq\Lambda\Vert\D u\Vert_{L^2(\bb R^n)}+\Vert|V|^{\frac12}u\Vert_{L^2(\bb R^n)}
\]
hold for $u\in\dot{\m V}$. Actually, if we further assume that $V$ is real-valued, then the map
\[
\langle u, v\rangle_{\dot{\f l}_I}:=\int_{\bb R^n}\D u\overline{\D v}+Vu\overline v,\qquad\text{for each } u,v\in\dot{\m V},
\]
is an inner product on $\dot{\m V}$, and the induced norm is equivalent to $\Vert\cdot\Vert_{\f l}$. Hence $\dot{\m V}$ can be seen as a Hilbert space when $V$ is real-valued.

Define the operator $\dot L:\dot{\m V}\ra\dot{\m V}'$ in the following way: for $u\in\dot{\m V}$, $\dot Lu$ is the functional $f\in\dot{\m V}'$ given by
\[
f(v)=\dot{\f l}(u,v),\qquad\text{for each } v\in\dot{\m V}.
\]
Clearly, $\dot L$ is a linear operator, which is bounded on $\dot{\m V}$ (this is proven similarly to the continuity of $\f l$). By the Lax-Milgram Theorem, it is also invertible, so that $\dot L^{-1}:\dot{\m V}'\ra\dot{\m V}$ exists and is unique. This means that for all $f\in\dot{\m V}'$, there exists a unique $u\in\dot{\m V}$ such that
\[
\int_{\bb R^n}A\D u\overline{\D v}+Vu\overline{v}=(f,v),\qquad\text{for each } v\in\dot{\m V},
\]
where $(f,v)$ is the duality pairing of $\dot{\m V}'$ with $\dot{\m V}$.

The following proposition says that the space of compactly supported $L^2(\bb R^n)$ functions can be seen as a subspace of $\dot{\m V}'$. The proof is omitted.

\begin{proposition} Assume that $\mathbf{a}\in L^2_{loc}(\bb R^n)$, $A$ is an elliptic matrix with complex, bounded, measurable coefficients, $V\in L^{1}_{loc}(\bb R^n)$ satisfies (\ref{imVreq}) and (\ref{gkcond}) with $c_2\equiv c_4\equiv0$. Suppose $f\in L^2(\bb R^n)$ is compactly supported. Then $f\in\dot{\m V}'$, and
\begin{equation}\label{embeddingdual}
\Vert f\Vert_{\dot{\m V}'}\leq C|\supp f|^{\frac1n}\Vert f\Vert_{L^2(\bb R^n)},
\end{equation}
\begin{equation}\label{embeddingdual1}
\Vert\dot{L}^{-1}f\Vert_{\dot{\m V}}\leq C|\supp f|^{\frac1n}\Vert f\Vert_{L^2(\bb R^n)},
\end{equation}
where $C$ is a constant depending only on $\lambda, c_1,$ and $n$.
\end{proposition}

We now prove
\begin{lemma}\label{opstrongconv} Assume that $\mathbf{a}\in L^2_{loc}(\bb R^n)$, $A$ is an elliptic matrix with complex, bounded, measurable coefficients, $V\in L^{1}_{loc}(\bb R^n)$ satisfies (\ref{imVreq}) and (\ref{gkcond}) with $c_2\equiv c_4\equiv0$. Suppose $f\in\dot{\m V}'\cap L^2(\bb R^n)$. For $\ep>0$, let $u_{\ep}=(L+\ep)^{-1}f\in D(\f l)$. Then $\{u_{\ep}\}$ is a bounded sequence in $\dot{\m V}$ which converges strongly in $\dot{\m V}$ to $\dot L^{-1}f$. In particular, $\{u_{\ep}\}$ converges to $\dot L^{-1}f$ strongly in the topology of $L^2_{loc}(\bb R^n)$, and a subsequence converges pointwise a.e. on $\bb R^n$. 
\end{lemma}

\noindent\emph{Proof.} By definition of the sequence $\{u_{\ep}\}$, we have
\begin{equation}\label{weaky}
\int_{\bb R^n}A\D u_{\ep}\overline{\D v}+(V+\ep)u_{\ep}\overline v=\int_{\bb R^n}f\,\overline v,\qquad\text{for each } v\in\m V=D(\f l),
\end{equation}
and in particular, since $u_{\ep}\in\m V$, we can write
\[
\int_{\bb R^n}A\D u_{\ep}\overline{\D u_{\ep}}+(V+\ep)|u_{\ep}|^2=\int_{\bb R^n}f\,\overline u_{\ep},
\]
and from this we obtain (since the right-hand side can be re-written as the duality pairing $(f,\overline{u_{\ep}})$)
\[
\Vert u_{\ep}\Vert_{\dot{\m V}}^2=\Re e\,\int_{\bb R^n}A\D u_{\ep}\overline{\D u_{\ep}}+V|u_{\ep}|^2\leq\Vert f\Vert_{\dot{\m V}'}\Vert u_{\ep}\Vert_{\dot{\m V}},
\]
yielding the boundedness of the sequence $\{u_{\ep}\}$ in $\dot{\m V}$, with
\begin{equation}\label{seqbounded}
\Vert u_{\ep}\Vert_{\dot{\m V}}\leq\Vert f\Vert_{\dot{\m V}'}.
\end{equation}
Hence, the sequence has a weak limit, say, $u\in\dot{\m V}$. Since by the diamagnetic inequality (\ref{diamagnetic}) we have
\begin{align*}
\Bigl|~\int_{\bb R^n}\ep u_{\ep}\overline{\phi}\Bigr|&\leq C_{\phi}\ep\Vert u_{\ep}\Vert_{L^{2^*}(\bb R^n)}\leq C_{\phi}\ep\Vert\nabla|u_{\ep}|\Vert_{L^{2}(\bb R^n)}\\[4mm]&\leq C_{\phi}\ep\Vert\D u_{\ep}\Vert_{L^{2}(\bb R^n)}\leq C_{\phi}\ep\Vert u_{\ep}\Vert_{\dot{\m V}}\leq C_{\phi,f}\ep,\qquad\text{for each }\phi\in C_c^{\infty}(\bb R^n),
\end{align*}
then by taking limit as $\ep\ra0$ on (\ref{weaky}), we get that
\[
\int_{\bb R^n}A\D u\overline{\D\phi}+Vu\overline\phi=\int_{\bb R^n}f\,\overline \phi,\qquad\text{for each } \phi\in C_c^{\infty}(\bb R^n),
\]
and hence for all $\phi\in\dot{\m V}$, as $u\in\dot{\m V}$. In other words, $u=\dot L^{-1}f\in\dot{\m V}$. By the uniqueness of $\dot L^{-1}f$, it follows the whole sequence $\{u_{\ep}\}$ converges weakly to $\dot L^{-1}f$. Now,
\begin{align}
\Re e\,(f,u)=\Vert u\Vert_{\dot{\m V}}^2&\leq\liminf\limits_{\ep\ra0}\Vert u_{\ep}\Vert_{\dot{\m V}}^2\leq\limsup\limits_{\ep\ra0}\Vert u_{\ep}\Vert_{\dot{\m V}}^2\nonumber\\[4mm]&\leq\limsup\limits_{\ep\ra0}\Big[\Vert u_{\ep}\Vert_{\dot{\m V}}^2+\int_{\bb R^n}\ep|u_{\ep}|^2\Big]=\limsup\limits_{\ep\ra0}\Re e\, (f,u_{\ep})=\Re e\,(f,u)\label{strongconvarg},
\end{align}
which implies $\Vert u_{\ep}\Vert_{\dot{\m V}}\ra\Vert u\Vert_{\dot{\m V}}$. This, together with the weak convergence, gives the strong convergence in $\dot{\m V}$.\hfill{$\square$}\\

\subsection{Local solutions to the magnetic Schr\"odinger operator and their properties}
It will also be of interest to define local solutions for the operator $L$ given in (\ref{generalized}); under certain conditions, such local solutions will enjoy a Caccioppoli-type estimate and a Moser estimate (but we don't prove the latter result until Section \ref{fundsolnsec}). First, if $\Omega\subset\bb R^n$, we define
\[
\m V_{\mathbf{a},V}(\Omega)=\Big\{u\text{ measurable s.t. }\D u\in L^2(\Omega)\text{ and }|V|^{\frac12}u\in L^2(\Omega)\Big\},
\]
and also
\[
\m V_{\mathbf{a},V,loc}(\Omega)=\Big\{u\text{ measurable s.t. }\D u\in L^2_{loc}(\Omega)\text{ and }|V|^{\frac12}u\in L^2_{loc}(\Omega)\Big\},
\]
with $\m V_{\mathbf{a},V,0}(\Omega)$, the completion of $C_c^{\infty}(\Omega)$ under the topology associated to $\m V_{\mathbf{a},V}(\Omega)$. We omit the subscripts when possibility of confusion is slim. We remark that, in particular, the space of functions with compact support which lie in $\m V(\Omega)$ is a subset of $\m V_{0}(\Omega)$, and if $(\ref{gkcond}), (\ref{imVreq})$ are satisfied by $V$ with $c_2\equiv c_4\equiv0$, then the elements of $\dot{\m V}$ lie in $\m V(\Omega)$ for any $\Omega\subset\bb R^n$. Furthermore, elements of $\m V_{loc}(\Omega)$ are locally square integrable. Indeed, since $\D u\in L^2_{loc}(\Omega)$, then by the diamagnetic inequality (\ref{diamagnetic}), $\nabla|u|\in L^2_{loc}(\Omega)$, and from this, we conclude $u\in L^2_{loc}(\Omega)$.  Now, if $f\in(\m V_{0}(\Omega))',$ we say that a measurable function $u$ solves $Lu=f$ on $\Omega$ \emph{in the weak sense} if $u\in\m V_{loc}(\Omega)$ and
\begin{equation}\label{weak}
\int_{\Omega}A\D u\overline{\D\phi}+Vu\overline{\phi}=(f,\phi)
\end{equation}
for every $\phi\in C_c^{\infty}(\Omega)$. By a standard limiting argument it is clear that if $u$ solves $Lu=f$ on $\Omega$ in the weak sense, then (\ref{weak}) is satisfied for $\phi\in\m V_{0}(\Omega)$. It's also clear that $u=\dot L^{-1}f$ also solves $Lu=f$ on $\Omega$ in the weak sense. We present a generalized Caccioppoli inequality, whose proof is standard:

\begin{theorem}\label{Cacciopoli} Assume that $\mathbf{a}\in L^2_{loc}(\bb R^n)$, $A$ is an elliptic matrix with complex, bounded, measurable coefficients, $V\in L^1_{loc}(\bb R^n)$ satisfies (\ref{gkcond}) and (\ref{imVreq}) (with $c_2\in[0,\infty)$ and $c_4$ either $0$ or $1$). Suppose $f\in\Big(\m V(B_R)\Big)'\cap L^2_{loc}(B_R)$, and that $Lu=f$ on $B(x_0,R)\subset\bb R^n$ in the weak sense. Then
\begin{equation}\label{cacc}
\int_{B(x_0,r)}|\D u|^2\leq C\Bigl\{\Bigl[\frac1{(R-r)^2}+c_2\Bigr]\int_{B(x_0,R)}|u|^2+\int_{B(x_0,R)}|f||u|\Bigr\}.
\end{equation}
for every $r,\,0<r<R$, where $C$ is a constant depending only on $c_1,\lambda,\Lambda$.
\end{theorem}

{\section{The Fefferman-Phong-Shen maximal function and related properties}\label{mproperties}
Let $\n B$ be the collection of all balls in $\bb R^n$, and define
\[
\Vert w\Vert_{RH_p}:=\sup\limits_{B\in\n B}\frac{\left(\fint_Bw^p\right)^{1/p}}{\fint_Bw}.
\]

The following results are well-known in the theory of the Reverse-H\"older classes (see, for instance, \cite{St}):

\begin{proposition} If $w\in RH_p$ for some $p\geq1$, then $w\in RH_r$ for each $r\in[1,p]$.
\end{proposition}


\begin{proposition}\label{openend}\cite{G} If $w\in RH_p$, for some $p\geq1$, then there exists $\ep>0$ depending only on $\Vert w\Vert_{RH_p}, p,$ and $n$, such that $w\in RH_{r}$ for every $r\in[p,p+\ep)$. Moreover, $\Vert w\Vert_{RH_{r}}$ depends only on $\Vert w\Vert_{RH_p},p,$ and $n$.
\end{proposition}

\begin{proposition}\label{RHisdoubling} If $w\in RH_p$ for some $p>1$, then there exists a constant $C_0$ depending only on $\Vert w\Vert_{RH_p}, p,$ and $n$ so that for any ball $B\subset\bb R^n$,
\begin{equation}\label{doubling}
\int_{2B}w\leq C_0\int_{B}w,
\end{equation}
where $2B$ denotes the ball with same center as $B$ and twice the radius of $B$.
\end{proposition}

\emph{Proof.} Fix a point $x\in\bb R^n$ and $r>0$. Then, for any $R>r$, we observe that
\begin{multline}
\int_{B(x,R)\backslash B(x,r)}w\leq\Big(\int_{B(x,R)\backslash B(x,r)}w^p\Big)^{1/p}|B(x,R)\backslash B(x,r)|^{1-\frac1p}\\[4mm]\leq\Big(\frac{R^n}{|B(0,1)|R^n}\int_{B(x,R)}w^p\Big)^{1/p}|B(0,1)|(R^n-r^n)^{1-\frac1p}\\[4mm]\leq\Big(\frac{\Vert w\Vert_{RH_p}}{|B(0,1)|R^n}\int_{B(x,R)}w^p\Big)|B(0,1)|R^{n/p}(R^n-r^n)^{1-\frac1p}\\[4mm]= \Vert w\Vert_{RH_p} \left[1-\frac{r^n}{R^n}\right]^{1-\frac1p}\int_{B(x,R)}w^p, \label{bounds}
\end{multline}
where in the first inequality we used H\"older's Inequality, and in the third one we used the fact that $w\in RH_p$. Let $R_w>0$ satisfy
\[
\Vert w\Vert_{RH_p}\left[1-\frac{r^n}{R_w^n}\right]^{1-\frac1p}=\frac12,
\]
that is,
\begin{equation}\label{rw}
\frac{R_w}{r}=\left[1-\left(\frac1{2\Vert w\Vert_{RH_p}}\right)^{\frac p{p-1}}\right]^{-1/n}=:\alpha,
\end{equation}
and we note that $\alpha$ does not depend on $x$ or $r$, and $\alpha>1$. Then, for each $R\in(r,\alpha r]$, due to (\ref{bounds}) we observe that
\begin{gather}
\int_{B(x,R)}w=\int_{B(x,r)}w+\int_{B(x,R)\backslash B(x,r)}w\nonumber\\[4mm]\leq\int_{B(x,r)}w +\frac12\int_{B(x,R)}w\nonumber\\[5mm]\implies \int_{B(x,R)}w\leq2\int_{B(x,r)}w.\label{bounds2}
\end{gather}
Estimate (\ref{bounds2}) holds independently of $x$ and $r$, provided $R\in(r,\alpha r]$. Since $\alpha>1$, there exists $M\in\bb N$ such that $2\leq\alpha^M$. For example, we can take
\[
M:=\left\lfloor\frac{\ln 2}{\ln\alpha}\right\rfloor+1,
\]
where $\lfloor\cdot\rfloor$ is the floor function. Thus, upon applying (\ref{bounds2}) $M$ times, we have
\[
\int_{B(x,2r)}w\leq\int_{B(x,\alpha^Mr)}w\leq2^M\int_{B(x,r)}w\leq2^{1+\frac{\ln2}{\ln\alpha}}\int_{B(x,r)}w,
\]
which gives (\ref{doubling}) with
\[
C_0:=2^{1+\frac{\ln2}{\ln\alpha}}.
\]\hfill{$\square$}

For a function $w\in RH_p, p\geq\frac n2$, recall we define the maximal function $m(x,w)$ by (\ref{maximal-intro}), and the distance function $d(x,y,w)$ by (\ref{distfn}). For each $r>0$, let
\[
\Phi(r):=\frac1{r^{n-2}}\int_{B(x,r)} w.
\]
Whenever $w$ is understood from context, we will simply write $d(x,y)=d(x,y,w)$. It is straightforward to prove that $d$ satisfies the triangle inequality.
%

%

The following results were proven in \cite{S4}; here we expose the results while keeping a careful account of the constants.

\begin{proposition}\label{m=1} \cite{S4} Suppose $ w\in L^1_{loc}(\bb R^n)$. Fix $x\in\bb R^n$, and let $\hat r=\frac1{m(x, w)}$. Suppose that $\hat r\in(0,+\infty)$. Then
\[
\Phi(\hat r)=\frac1{\hat r^{n-2}}\int_{B(x,\hat r)} w=1.
\]	
\end{proposition}

\bp Suppose otherwise. We note that $\Phi(r)$ is continuous on $(0,+\infty)$ (because it is the product of two continuous functions). So if $\Phi(\hat r)<1$, then there exists $\delta>0$ such that $\Phi(r)<1$ for all $r\in(\hat r-\delta,\hat r+\delta)$. If $\hat r<+\infty$, then this contradicts the definition of $\hat r$ as the supremum of $\n R:=\{r>0~|~\Phi(r)\leq1\}$. So $\Phi(\hat r)\geq1$. But if $\Phi(\hat r)>1$, then again there exists $\delta>0$ such that $\Phi(r)>1$ for all $r\in(\hat r-\delta,\hat r+\delta)$, so there can be no sequence $\{r_n\}\subset\n R$ which converges to $\hat r$. It follows that $\hat r$ cannot be a supremum of $\n R$. Thus $\Phi(\hat r)=1$.\hfill{$\square$}

\begin{lemma}\label{Vprop} \cite{S4} Let $ w\in RH_p,~ p\geq\frac n2$, and $0<r<R$. Then
\[
\Phi(r)\leq\Vert w\Vert_{RH_p}\Bigl(\frac Rr\Bigr)^{\frac np-2}\Phi(R).
\]
That is,
\begin{gather}\label{Vprop1}
\frac1{r^{n-2}}\int_{B(x,r)} w\leq \Vert w\Vert_{RH_p}\Bigl(\frac Rr\Bigr)^{\frac np-2}\frac1{R^{n-2}}\int_{B(x,R)} w.
\end{gather}
\end{lemma}\ \\
\bp We observe that
\begin{multline*}
\frac1{|B(x,r)|}\int_{B(x,r)} w\leq\Bigl(\frac1{|B(x,r)|}\int_{B(x,r)} w^p\Bigr)^{\frac1p}\leq\Bigl(\frac Rr\Bigr)^{\frac np}\Bigl(\frac1{|B(x,R)|}\int_{B(x,R)} w^p\Bigr)^{\frac1p}\\[3mm]\leq \Vert w\Vert_{RH_p}\Bigl(\frac Rr\Bigr)^{\frac np}\Bigl(\frac1{|B(x,R)|}\int_{B(x,R)} w\Bigr),
\end{multline*}
where the first inequality follows from H\"older's Inequality, and the last one is due to $ w\in RH_p$. Multiplying both sides by $R^2r^2$ now gives (\ref{Vprop1}).\hfill{$\square$}

\begin{remark} From Lemma \ref{Vprop}, it is easy to conclude that, if $ w\in RH_p,\,p\geq\frac n2$ and $w>0$ a.e. on $\bb R^n$, then $m(\cdot, w)$ only takes values in $(0,+\infty)$.
\end{remark}

Observe that by Proposition \ref{openend}, if $w\in RH_{\frac n2}$, we may actually assume that $w\in RH_p$ for some $p>\frac n2$.

\begin{proposition}\label{similar} \cite{S4} Let $ w\in RH_p, p>\frac n2$. Fix $x\in\bb R^n$ and let $\hat r=\frac1{m(x, w)}$. Suppose that $r>0$ satisfies $\Phi(r)\sim1$. That is, there exists a constant $\hat C\geq1$ such that $\frac1{\hat C}\leq\Phi(r)\leq\hat C$. Then $r\sim\hat r$, with constant $\tilde C:=\Big(\Vert w\Vert_{RH_p}\hat C\Big)^{\frac1{2-\frac np}}$. Furthermore, the converse is true: if $r\sim\hat r$, then $\Phi(r)\sim1$.\\
\end{proposition}

\bp First suppose that $r<\hat r$. Then we can apply Lemma \ref{Vprop} to see that
\begin{gather*}
\frac1{\hat C}\leq\Phi(r)\leq \Vert w\Vert_{RH_p}\Bigl(\frac{\hat r}r\Bigr)^{\frac np-2}\Phi(\hat r)=\Vert w\Vert_{RH_p}\Bigl(\frac{\hat r}r\Bigr)^{\frac np-2}
\end{gather*}
since $\Phi(\hat r)=1$. Now, $p>\frac n2$ implies $\frac np-2<0$. As such, we have
\[
\hat r\leq\Big(\Vert w\Vert_{RH_p}\hat C\Big)^{\frac1{2-\frac np}}r.
\]
So in this case we are able to write
\begin{equation}\label{equivalent}
\frac1{(\Vert w\Vert_{RH_p}\hat C)^{\frac1{2-\frac np}}}\leq\frac{\hat r}r\leq\Big(\Vert w\Vert_{RH_p}\hat C\Big)^{\frac1{2-\frac np}}.
\end{equation}
If instead we have $r>\hat r$, similarly we get
\[
\frac1{(\Vert w\Vert_{RH_p}\hat C)^{\frac1{2-\frac np}}}<\frac r{\hat r}\leq\Big(\Vert w\Vert_{RH_p}\hat C\Big)^{\frac1{2-\frac np}}
\]
which gives (\ref{equivalent}) in this case too. It follows that (\ref{equivalent}) is true regardless of whether $r<\hat r$. This proves the result.\hfill{$\square$}

\begin{lemma}\label{propertiesm} Let $ w\in RH_p,~ p>\frac n2$. Then, for any constant $C>0$ and any $x,y\in\bb R^n$ we have
\begin{align*}
\Big(\Vert w\Vert_{RH_p}C_0^{1+\log_2(C+1)}\Big)^{-\frac1{2-\frac np}}\leq\frac{m(x, w)}{m(y, w)}\leq \Big(\Vert w\Vert_{RH_p}&C_0^{1+\log_2(C+1)}\Big)^{\frac1{2-\frac np}}\\[4mm]&\text{if}\quad|x-y|\leq\frac{C}{m(x, w)},
\end{align*}
where $C_0$ is the constant from (\ref{doubling}), which depends on $\Vert w\Vert_{RH_p}, p,$ and $n$ only. Therefore,
\begin{equation}\label{slow}
m(x, w)\sim m(y, w)\quad\text{if}\quad|x-y|\leq\frac{C}{m(x, w)}.
\end{equation}
\end{lemma}

\bp Let $r:=\frac1{m(x, w)}$. Suppose that $|x-y|\leq Cr$. Then $B(y,r)\subset B(x,(C+1)r)$. Therefore,
\[
\frac1{r^{n-2}}\int_{B(y,r)} w\leq\frac1{r^{n-2}}\int_{B(x,(C+1)r)} w\leq C_0^{1+\log_2(C+1)}\frac1{r^{n-2}}\int_{B(x,r)} w= C_0^{1+\log_2(C+1)}.
\]
Similarly, we have
\[
1=\frac1{r^{n-2}}\int_{B(x,r)} w\leq C_0^{1+\log_2(C+1)} \frac1{r^{n-2}}\int_{B(y,r)} w,
\]
so that
\[
\frac1{C_0^{1+\log_2(C+1)}}\leq\frac1{r^{n-2}}\int_{B(y,r)} w\leq C_0^{1+\log_2(C+1)},
\]
which by Proposition \ref{similar} implies that
\[
\frac1{\Big(\Vert w\Vert_{RH_p}C_0^{1+\log_2(C+1)}\Big)^{\frac1{2-\frac np}}}\leq rm(y, w)\leq\Big(\Vert w\Vert_{RH_p}C_0^{1+\log_2(C+1)}\Big)^{\frac1{2-\frac np}}.
\]
This means $r\sim\frac1{m(y, w)}$, or equivalently, $m(x, w)\sim m(y, w)$.\hfill{$\square$}\\

\begin{remark}\label{mbounded} From Lemma \ref{propertiesm} it follows that if $\Omega\subset\bb R^n$ is a bounded set, then there exists a constant $C$, depending on $\Omega$, such that
\begin{equation}\label{bounded}
\frac1C\leq m(x,w)\leq C.
\end{equation}
\end{remark}

The following results are proved in \cite{S4}.

\begin{lemma}\label{polygrowth} \cite{S4} Let $ w\in RH_p, p\geq\frac n2$. Then there exist constants $C, c, k_0>0$, depending only on $p,\Vert w\Vert_p, n$, such that for any $x,y\in\bb R^n$,
\begin{equation}\label{polygrowth1}
m(x, w)\leq Cm(y, w)[1+|x-y|m(y, w)]^{k_0},
\end{equation}
and
\begin{equation}\label{polygrowth2}
m(x, w)\geq\frac{cm(y, w)}{\big[1+|x-y|m(y, w)\big]^{\frac{k_0}{k_0+1}}}.
\end{equation}
\end{lemma}

\begin{proposition}\label{intbound}\cite{S4} Let $n\geq3$, $~ w\in RH_p,~ p>\frac n2$. Then there is a constant $C$, which is a numerical multiple of $C_0$ from (\ref{doubling}), such that for every $x\in\bb R^n$ and every $R>0$,
\begin{equation}\label{Kato}
\int_{B(x,R)}\frac{ w(z)\,dz}{|z-x|^{n-2}}\leq C\frac1{R^{n-2}}\int_{B(x,R)} w(z)\,dz.
\end{equation}
\end{proposition}

\bp If $p=\frac n2$, then by Proposition \ref{openend}, $w\in RH_{p+\ep}$ for some $\ep>0$. Therefore, we can assume $p>\frac n2$ without loss of generality. Let $q$ be the H\"older conjugate of $p$. Since $p>\frac n2$, then $q<\frac n{n-2}$. By H\"older's Inequality we observe
\begin{gather*}
\int_{B(x,R)}\frac{ w(z)\,dz}{|z-x|^{n-2}}\leq\left(\int_{B(x,R)} w^p\right)^{\frac1p}\left(\int_{B(x,R)}\frac{dz}{|z-x|^{q(n-2)}}\right)^{\frac1q}\\[3mm]\leq \left(\int_{B(x,R)}\frac{dz}{|z-x|^{q(n-2)}}\right)^{\frac1q}C_0|B(x,R)|^{\frac1p-1}\int_{B(x,R)} w.
\end{gather*}
But the integral on the right-hand side can be estimated in spherical coordinates:
\begin{gather*}
\int_{B(x,R)}\frac{dz}{|z-x|^{q(n-2)}}\leq C(n)\int_0^R\frac1{r^{q(n-2)}}\,r^{n-1}\,dr=C(n)\,\frac1{n-q(n-2)}\,r^{n-q(n-2)}\Big|_0^R\\[3mm]\leq C(n)\,\frac1{n-q(n-2)}\,R^{n-q(n-2)}.
\end{gather*}
It follows that
\[
\int_{B(x,R)}\frac{ w(z)\,dz}{|z-x|^{n-2}}\leq \left(C(n)\frac1{n-q(n-2)}\right)^{\frac1q}R^{\frac nq}R^{-(n-2)}C_0(C(n))^{-\frac1q}R^{-\frac nq}\int_{B(x,R)} w,
\]
which gives (\ref{Kato}), with constant $C=\Big(n-q(n-2)\Big)^{-\frac1q}C_0$.\hfill{$\square$}

%
%
%
%
%
%
%
%

We also note the following useful observation:

\begin{proposition}\label{propdbounded} Let $ w\in RH_p,~ p>\frac n2$. Let $x,y\in\bb R^n$. Then for any constant $C>0$ there exists a constant $\m C$ depending on $C,\Vert w\Vert_{RH_p}, p,$ and $n$ only such that
\begin{equation}\label{dbounded}
d(x,y, w)\leq\m C\qquad\text{if }\qquad|x-y|\leq\frac{C}{m(x, w)}.
\end{equation}
\end{proposition}
\bp Let $\gamma$ be the straight line that connects $x$ to $y$. Then
\[
|x-\gamma(t)|\leq|x-y|\leq\frac{C}{m(x, w)},\qquad\text{for each } t\in[0,1],
\]
per the hypothesis. By (\ref{slow}), we see that $m(\gamma(t), w)\sim m(x, w)$ for each $t\in[0,1]$, so that there exists a constant $C_1>0$ satisfying $m(\gamma(t), w)\leq C_1m(x, w)$ for all $t\in[0,1]$. By definition of $d$, we have
\[
d(x,y, w)\leq\int_0^1m(\gamma(t), w)|\gamma'(t)|\,dt\leq C_1\int_0^1m(x, w)|\gamma'(t)|\,dt=C_1m(x, w)|x-y|\leq C_1C
\]
Finally, take $\m C:=C_1C$.\hfill{$\square$}

The following proposition will be useful in proving the lower bound exponential decay estimate of Section \ref{lowerboundsec}.

\begin{proposition}\label{Aprop}\cite{S1} Let $ w\in RH_p,~p>\frac n2$. Then there exists a constant $\m C>0$, depending on $\Vert w\Vert_{RH_p}$ and $n$ only, which satisfies
\begin{equation}\label{Aexists}
x\notin B\Bigl(y,\frac2{m(y, w)}\Bigr)\qquad\text{whenever }~|x-y|\geq\frac{\m C}{m(x, w)}.
\end{equation}
\end{proposition}
\bp Suppose otherwise, so that for each $j\in\bb N$, there are points $x_j,y_j\in\bb R^n$ which satisfy
\[
|x_j-y_j|\geq\frac j{m(x_j, w)}\qquad\text{but }~x_j\in B\Bigl(y_j,\frac2{m(y_j, w)}\Bigr)
\]
whence we have
\begin{equation}\label{chain}
\frac j{m(x_j, w)}\leq|x_j-y_j|\leq\frac2{m(y_j, w)},\qquad\text{for each } j\in\bb N.
\end{equation}

Now, the fact that $x_j\in B\Bigl(y_j,\frac2{m(y_j, w)}\Bigr)$ for each $j\in\bb N$ implies by Lemma \ref{propertiesm} that there is a constant $C$, independent of $j$, such that
\[
m(x_j, w)\leq Cm(y_j, w),\qquad\text{for each } j\in\bb N.
\]
Using this result on (\ref{chain}) gives
\[
m(y_j, w)j\leq2m(x_j, w)\leq2Cm(y_j, w),\quad\text{for each } j\in\bb N
\]
\[
\implies j\leq2C,\qquad\text{for each } j\in\bb N
\]
which is a contradiction. This establishes the result.\hfill{$\square$}

%
%
%

Next, we present an estimate often known as a \emph{Fefferman-Phong inequality} for the magnetic Schr\"odinger operator. We cite the statement from \cite{benali}; the proof of the statement is an easy generalization of the proof first written in \cite{S3}:


\begin{theorem}\label{fptheorem} Suppose that $\mathbf{a}\in L^2_{loc}(\bb R^n)^n$, and moreover assume (\ref{magassumptions-intro}).
Then, for all $u\in C^1_c(\bb R^n)$,
\begin{equation}\label{fpmag}
\int_{\bb R^n}m^2(x,V+|\mbf B|)|u|^2\,dx\leq C\int_{\bb R^n}(|\D u|^2+V|u|^2)\,dx,
\end{equation}
where $C$ depends on the constants $c,c'$ from (\ref{magassumptions-intro}) and on $\Vert V+|\mathbf{B}|\Vert_{RH_{\frac n2}}$.
\end{theorem}

At this juncture we observe that the Fefferman Phong Inequality of Theorem \ref{fptheorem} is preserved in the class of functions $\dot{\m V}_{\mathbf{a},V}$:

\begin{corollary}\label{fpcorollary} Suppose that $\mathbf{a}\in L^2_{loc}(\bb R^n)$, $A$ is an elliptic matrix with complex, bounded, measurable coefficients, and assume (\ref{magassumptions-intro}).	Then, for all $u\in\dot{\m V}_{\mathbf{a},V}$,
\begin{equation}\label{fpmag1}
\int_{\bb R^n}m^2(x,V+|\mbf B|)|u|^2\,dx\leq C\int_{\bb R^n}(\Re e\,A\D u\overline{\D u}+V|u|^2)\,dx,
\end{equation}
where $C$ depends on the constants $c,c'$ from (\ref{magassumptions-intro}) and on $\lambda,\Vert V+|\mathbf{B}|\Vert_{RH_{\frac n2}}$. Note that $V$ is a real (rather than complex) valued function, which is, in addition, non-negative.
\end{corollary}

}

{ \section{$L^2$ Exponential decay}\label{l2sec}

For notational simplicity, for $t>0$ we will write
\[
L_t:=L+\frac1{t^2},\qquad V_t:=V+\frac1{t^2},\qquad R_t:=(1+t^2L)^{-1}=\frac1{t^2}L_t^{-1},
\]
and identify the case $t=\infty$ with $L_{\infty}:=L, V_{\infty}=V$. Per the discussion in Section \ref{theory}, the family of operators $\{R_t\}$ is uniformly bounded from $L^2(\bb R^n)$ into $L^2(\bb R^n)$. We assume that $V$ is real-valued, $V\geq0$, and emphasize that for $U$ an open bounded set in $\bb R^n$,  in the definition of the spaces $\m V_{\mathbf{a},V,loc}(U)$ we always take the weight to be $V$, even when studying the operators $L_t$. This can be done, because if $\left(V+\frac1{t^2}\right)^{\frac12}u\in L^2_{loc}(U)$, then it necessarily follows that $V^{\frac12}u\in L^2_{loc}(U)$. On the other hand, if $V^{\frac12}u\in L^2_{loc}(U)$ and $\D u\in L^2_{loc}(U)$, then by the diamagnetic inequality it follows that $u\in L^2_{loc}(U)$, which implies $\frac1{t^2}u\in L^2_{loc}(U)$, and so $\left(V+\frac1{t^2}\right)^{\frac12}u\in L^2_{loc}(U)$. For a bounded set $U\subset\bb R^n$ and $c>1$, we define
\[
cU:=\Big\{x\in\bb R^n~|~\exists y\in U\text{ s.t. }|x-y|<\frac{c-1}2\diam U\Big\},
\]
that is, $cU$ is a $\frac{c-1}2\text{diam }U-$ neighborhood of $U$. Since for each $t>0$, $V_t+|\mathbf{B}|\geq V+|\mathbf{B}|$ pointwise a.e. on $\bb R^n$, then by the definition of the Fefferman-Phong-Shen maximal function, it is easy to see that
\begin{equation}\label{tmdomination}
m(x,V+|\mathbf{B}|)\leq m(x,V_t+|\mathbf{B}|),\qquad\text{for each } x\in\bb R^n,
\end{equation}
and hence also
\begin{equation}\label{tddomination}
d(x,y,V+|\mathbf{B}|)\leq d(x,y,V_t+|\mathbf{B}|),\qquad\text{for each } x,y\in\bb R^n.
\end{equation}
Furthermore, if for some $p\in[1,\infty)$ we have $V+|\mathbf{B}|\in RH_p$ and $V+|\mathbf{B}|\neq0$ a.e. on $\bb R^n$, then for any $t>0$, $V_t+|\mathbf{B}|\in RH_p$ with its $RH_{p}-$norm controlled by that of the former function. 


\begin{proposition}\label{proposition1} Assume that $\mathbf{a}\in L^2_{loc}(\bb R^n)$, $A$ is an elliptic matrix with complex, bounded, measurable coefficients, and $V\in L^1_{loc}(\bb R^n)$ is real-valued with $V\geq0$ a.e. on $\bb R^n$. If $\mathbf{a}\equiv0$, assume $V\in RH_{\frac n2}$, otherwise take assumptions (\ref{magassumptions-intro}). Suppose $ U\subset\bb R^n$ is a bounded open set. Let $t\in(0,\infty]$ and let $u\in\m V_{loc}(\bb R^n\backslash U)$ be a solution to $L_tu=0$ in the weak sense on $\bb R^n\backslash U$. Suppose $\phi\in C_c^{\infty}(\bb R^n)$ is such that $\phi\equiv0$ on $2 U$. Let $g=g_t\in C^{\infty}(\bb R^n)$ be a non-negative function satisfying $|\nabla g(x)|\leq C_2m(x,V_t+|\mathbf{B}|)$ for every $x\in\bb R^n$ ($C_2$ independent of $t$). Then
\begin{equation}\label{estimate}
\int_{\bb R^n}m(x,V_t+|\mathbf{B}|)^2|u\phi|^2e^{2\ep g}\,dx\leq C\int_{\bb R^n}|u|^2|\nabla\phi|^2e^{2\ep g}\,dx,
\end{equation}
for any $\ep\in(0,\ep_0)$, where $\ep_0$ and $C$ depend only on $\lambda,\Lambda, C_2, n, \Vert V+|\mathbf{B}|\Vert_{RH_{\frac n2}},$ and the constants from (\ref{magassumptions-intro}), but not on $t$.
\end{proposition}

\bp Let $f=\phi e^{\ep g}$ and $\psi=\psi_{\ep}=uf$. Let us prove that $\psi\in\m V_0(\bb R^n\backslash U)\cap\dot{\m V}$. Since $u\in\m V_{loc}(\bb R^n\backslash U)$, then $\D u\in L^2_{loc}(\bb R^n\backslash U)$ and $V^{\frac12}u\in L^2_{loc}(\bb R^n\backslash U)$. Since $f\in C_c^{\infty}(\bb R^n)$, we have that $V^{\frac12}uf\in L^2_{loc}(\bb R^n\backslash U)$. Owing to Corollary \ref{coolembedding}, we see that $u\in L^2_{loc}(\bb R^n\backslash U)$. Now, for any $\eta\in C_c^{\infty}(\bb R^n\backslash U)$,
\[
\int_{\bb R^n}(fu_{x_i}+uf_{x_i})\eta\,dx=-\int_{\bb R^n}u(f\eta)_{x_i}\,dx+\int_{\bb R^n}uf_{x_i}\eta\,dx=-\int_{\bb R^n}uf\eta_{x_i}\,dx,
\]
which proves that $uf$ is weakly differentiable with gradient $f\nabla u+u\nabla f$. Since
\[
\D\psi=\D(uf)=f\D u+u\nabla f,
\]
it follows that $\D\psi\in L^2_{loc}(\bb R^n\backslash U)$. Hence $\psi\in\m V_{loc}(\bb R^n\backslash U)$. The fact that $\psi$ has compact support within $\bb R^n\backslash U$ now implies that $\psi\in\m V_0(\bb R^n\backslash U)$. Since $f$ is compactly supported in $\bb R^n\backslash U$, we have that $V|\psi|^2=V|u|^2f^2\in L^1(\bb R^n)$. Similarly, $\D\psi\in L^2(\bb R^n)$. Hence, we actually have $\psi\in\dot{\m V}$.

By a virtually identical argument, it is easy to see that $uf^2\in\m V_0(\bb R^n\backslash U)$. We note that
\begin{multline}
\int_{\bb R^n}|\D\psi|^2+V_t|\psi|^2\leq\frac1{\lambda}\int_{\bb R^n}\Re e\,A\D\psi\overline{\D\psi}+V_t|\psi|^2\\[4mm]=\frac1{\lambda}\int_{\bb R^n}\Re e\,A\Big[f^2\D u\overline{\D u}+f\overline{u}\D u\nabla f+fu\nabla f\overline{\D u}+|u|^2\nabla f\nabla f\Big]+V_tu\overline{uf^2}\\[4mm]=\frac1{\lambda}\int_{\bb R^n}\Re e\,A\Big[\D u\overline{\D(uf^2)}+fu\nabla f\overline{\D u}-f\overline{u}\D u\nabla f+|u|^2\nabla f\nabla f\Big]+V_tu\overline{uf^2}.\label{absorb1}
\end{multline}
Since $u$ solves $L_tu=0$ on $\bb R^n\backslash U$ in the weak sense and $uf^2\in\m V_0(\bb R^n\backslash U)$, (\ref{absorb1}) reduces to
\begin{gather}
\int_{\bb R^n}|\D\psi|^2+V_t|\psi|^2\leq C\int_{\bb R^n}|u|^2|\nabla f|^2+C\int_{\bb R^n}\Re e\,A\Big[fu\nabla f\overline{\D u}-f\overline{u}\D u\nabla f\Big]=I+II.\label{absorb2}
\end{gather}
By the boundedness of $A$ and the Cauchy inequality with $\delta>0$, we observe that
\begin{equation}\label{absorb3}
|II|\leq\int_{\bb R^n}\delta f^2|\D u|^2+\int_{\bb R^n}\frac1{\delta}|u|^2|\nabla f|^2.
\end{equation}
Moreover, a straightforward computation yields
\begin{gather*}
|\D(uf)|^2=f^2|\D u|^2+2f\nabla f\Re e\,(\bar u\D u)+|u|^2|\nabla f|^2,
\end{gather*}
so that applying the Cauchy inequality we see that
\begin{gather}\label{absorb4}
\frac12f^2|\D u|^2\leq|\D(uf)|^2+|u|^2|\nabla f|^2.
\end{gather}
Putting together (\ref{absorb4}) with (\ref{absorb3}), (\ref{absorb2}), and the fact that $\psi=uf$, we achieve
\begin{gather}\label{absorb5}
\int_{\bb R^n}|\D\psi|^2+V_t|\psi|^2\leq C\int_{\bb R^n}|u|^2|\nabla f|^2,
\end{gather}
where $C$ is a constant which depends only on $\lambda,\Lambda, n$, and not on $t$. In passing, we note that by the diamagnetic inequality, (\ref{absorb5}) yields
\begin{equation}\label{interesting}
\int_{\bb R^n}|\nabla|\psi||^2+V_t|\psi|^2\leq C\int_{\bb R^n}|u|^2|\nabla f|^2.\nonumber
\end{equation} 
Using the Fefferman-Phong inequality (\ref{fpmag1}) and (\ref{absorb5}), it follows that we can write
\begin{align*}
\int_{\bb R^n}m(\cdot,V_t+|{\bf B}|)^2|\psi|^2&\leq C\int_{\bb R^n}|u|^2|\nabla f|^2\\[3mm]&\leq C\int_{\bb R^n}|u|^2|\nabla\phi|^2e^{2\ep g}+C\ep^2\int_{\bb R^n}|u\phi|^2e^{2\ep g}|\nabla g|^2\\[3mm]&\leq C\int_{\bb R^n}|u|^2|\nabla\phi|^2e^{2\ep g}+ CC_2^2\ep^2\int_{\bb R^n}m(\cdot,V_t+|\mathbf{B}|)^2|u\phi|^2e^{2\ep g},
\end{align*}
with $C$ independent of $t$. This implies that for $\ep$ small enough, we can absorb the right-most term into the left-hand side:
\[
\int m(\cdot,V_t+|\mathbf{B}|)^2|u\phi|^2e^{2\ep g}\leq C\int|u|^2|\nabla\phi|^2e^{2\ep g},
\]
and this proves the proposition.\hfill{$\square$}\\

Heuristically, the idea to prove (\ref{l2result-intro}) is that for fixed $y\in\bb R^n$, $g=d(\cdot,y,V_t+|\mathbf{B}|)$. However, the Agmon distance $d$ is not necessarily a smooth function. The next proposition shows that we can procure a continuous function which is close to $d$ in a uniform way, and which can be approximated by a sequence of bounded continuous functions. Both results are given in \cite{S3}.

\begin{proposition}\label{dapprox} Let $w\in RH_p,~ p\geq\frac n2$. Then for each $y\in\bb R^n$, there exists a non-negative function $\varphi(\cdot,y)=\varphi(\cdot,y,w)\in C^{\infty}(\bb R^n)$ such that for every $x\in\bb R^n$ and $y\in\bb R^n$,
\begin{equation}\label{closetod}
|\varphi(x,y)-d(x,y,w)|\leq C,
\end{equation}
and
\begin{equation}\label{gradientphi}
|\nabla_x\varphi(x,y)|\leq Cm(x,w),
\end{equation}
where the constants in (\ref{closetod}) and (\ref{gradientphi}) depend only on $\Vert w\Vert_{RH_p},p,$ and $n$.
\end{proposition}

\begin{proposition}\label{dapprox1} Let $w\in RH_p,~ p\geq\frac n2$. For each $y\in\bb R^n$, there exists a sequence of non-negative bounded $C^{\infty}$ functions $\{\varphi_{j}(\cdot,y)\}=\{\varphi_{j}(\cdot,y,w)\}$ such that, for every $x\in\bb R^n$,
\[
\varphi_{j}(x,y)\leq\varphi(x,y)\quad\text{and}\quad\varphi_{j}(x,y)\ra\varphi(x,y)\quad\text{as }j\ra\infty,
\]
and
\begin{equation}\label{star}
|\nabla_x\varphi_{j}(x,y)|\leq Cm(x, w),
\end{equation}
where $C$ depends on $\Vert w\Vert_{RH_p}, p,$ and $n$ only.
\end{proposition}


For an open bounded set $ U\subset\bb R^n$, let
\[
\n B_{U,t}:=\Big\{ B\Bigl(x,\frac 1{m(x, V_t+|\mathbf{B}|)}\Bigr)~\Big|~x\in U\Big\}.
\]
By the Besicovitch Covering Theorem (see \cite{DiB} for a proof), there exists a countable subcollection $\n B'_{U,t}$ of $\n B_{U,t}$ which covers $ U$, and for which there is uniformly finite overlap of the balls,
\[
\sum\limits_{B\in\n B'_{U,t}}\chi_B\leq c_n,
\]
with $c_n$ depending only on the dimension $n$ and not on $U$ nor $t$. Since $\overline{ U}$ is compact, there exists a finite subcollection $\n B''_{U,t}=\{B_k\}_{k=1}^K$ of $\n B'_{U,t}$ which covers $ U$. Let $\n F_{U,t}$ be the family of finite subcollections of $\n B_{U,t}$ which cover $ U$ with finite overlap at most $c_n$. Clearly, $\n B''_{U,t}\in\n F_{U,t}$, so this family is not empty. Let us first show
\begin{proposition}\label{dset} Assume that $\mathbf{a}\in L^2_{loc}(\bb R^n)$ and $V\in L^1_{loc}(\bb R^n)$ with $V\geq0$ a.e. on $\bb R^n$. If $\mathbf{a}\equiv0$, assume $V\in RH_{\frac n2}$, otherwise take assumptions (\ref{magassumptions-intro}). Then there exists $\tilde d$ depending on $\Vert V+|\mathbf{B}|\Vert_{RH_{\frac n2}}$ and $n$ only, such that
\[
\Big\{x\in\bb R^n~\big|~d(x,U, V_t+|\mathbf{B}|)\geq\tilde d\Big\}\subseteq\bb R^n\backslash\bigcup\limits_{B\in\n B''_{U,t}}4B,\qquad\text{for any }~\n B''_{U,t}\in\n F_{U,t}.
\]
\end{proposition}

\begin{theorem}\label{l2decay} Assume that $\mathbf{a}\in L^2_{loc}(\bb R^n)$, $A$ is an elliptic matrix with complex, bounded, measurable coefficients, and $V\in L^1_{loc}(\bb R^n)$ with $V\geq0$ a.e. on $\bb R^n$. If $\mathbf{a}\equiv0$, assume $V\in RH_{\frac n2}$, otherwise take assumptions (\ref{magassumptions-intro}). Suppose $f\in L^2(\bb R^n)$ is compactly supported, and let $t\in(0,\infty)$. Then there exists $\tilde d>0$, depending on $\Vert V+|\mathbf{B}|\Vert_{RH_{\frac n2}}$ and $n$ only, and there exists $\ep>0$ such that (\ref{l2result-intro}) holds,
where $\ep$ and $C$ depend on $\lambda,\Lambda, n,\Vert V+|\mathbf{B}|\Vert_{RH_{\frac n2}},$ and the constants from (\ref{magassumptions-intro}), but they are independent of $\supp f,t$. Moreover, there exists $\ep>0$ such that
\begin{multline}
\int_{\big\{x\in\bb  R^n|d(x,\supp f,V_t+|\mathbf{B}|)\geq\tilde d\big\}}m\left(\cdot,V_t+|\mathbf{B}|\right)^2\left|L_t^{-1}f\right|^2e^{2\ep d(\cdot,\supp f,V_t+|\mathbf{B}|)}\\\leq C\int_{\bb R^n}|f|^2\frac1{m(\cdot,V_t+|{\mathbf B}|)^2},\label{l2resultunif}
\end{multline}
where $\ep$ and $C$ depend on $\lambda,\Lambda, n,\Vert V+|\mathbf{B}|\Vert_{RH_{\frac n2}},$ and the constants from (\ref{magassumptions-intro}), but they are independent of $\supp f,t$.
\end{theorem}

\noindent\emph{Proof.} If $f=0$ on $\bb R^n$, then there is nothing to show, so suppose that $|\supp f|>0$. Fix $t>0$. Let $U$ be any open ball such that $\supp f\subset U$, and write $u:=R_tf$. By construction, $u$ is a weak solution to $L_tu=0$ on $\bb R^n\backslash U$ in the weak sense, since
\[
u=R_tf=\frac1{t^2}L_t^{-1}f~\implies~L_tu=\frac f{t^2}.
\]
Let $M>0$ such that $4U\subset B_M=B(0,M)$. Take $\phi\in C_c^{\infty}(\bb R^n)$ with $\phi\equiv0$ on $2 U$. Furthermore, suppose $\phi\equiv1$ on $B_M\backslash4 U$, $\phi\equiv0$ on $\bb R^n\backslash2B_M$, and
\begin{gather*}
|\nabla\phi|\leq\frac 2{\diam U}\quad\text{on }4 U\backslash2 U,\\
|\nabla\phi|\leq\frac2M\quad\text{on }2B_M\backslash B_M.\\
\end{gather*}
Fix $y\in\supp f$, $j\in\bb N$, and let $g=\varphi_j(x,y)=\varphi_j(x,y;V_t+|\mathbf{B}|)$ as in Proposition \ref{dapprox1}. For each $j\in\bb N$, $g=\varphi_j(x,y)$ is an admissible function in Proposition \ref{proposition1}. Then by \ref{estimate}, we have
\begin{multline}\label{greatineq}
\int_{B_M\backslash4 U}m(\cdot,V_t+|\mathbf{B}|)^2|u|^2e^{2\ep g}\leq\int_{\bb R^n}m(\cdot,V_t+|\mathbf{B}|)^2|u\phi|^2e^{2\ep g}\leq C\int_{\bb R^n}|u|^2|\nabla\phi|^2e^{2\ep g}\\[4mm]\leq C\Bigl\{\int_{4 U\backslash2 U}|u|^2\frac1{|\diam U|^2}e^{2\ep g}+\int_{2B_M\backslash B_M}|u|^2\frac1{M^2}e^{2\ep g}\Bigr\},
\end{multline}
where $C$ is independent of $j, M, y, U,t$. Since $g=\varphi_j(x,y)$ is a bounded function on $\bb R^n$ and by definition $u\in L^2(\bb R^n)$ because the resolvent maps $L^2(\bb R^n)$ into $L^2(\bb R^n)$, it follows that
\[
\Bigl|~\int_{2B_M\backslash B_M}|u|^2\frac1{M^2}e^{2\ep g}\Bigr|\leq C_j\frac1{M^2}\int_{\bb R^n}|u|^2\longrightarrow0\quad\text{as}\quad M\ra\infty.
\]
Consequenty, by Fatou's lemma,
\begin{equation}\label{great}
\int_{\bb R^n\backslash 4 U}m(\cdot,V_t+|\mathbf{B}|)^2|u|^2e^{2\ep g}\leq C\int_{4 U\backslash2 U}|u|^2\frac1{|\diam U|^2}e^{2\ep g},
\end{equation}
with $C$ independent of $j,y, U,t$.

Now let $\n B''\in\n F_{\supp f,t}$, and we can write $\n B''=\{B_k\}_{k=1}^K$. Let $\Omega=\bigcup\limits_{k=1}^K4B_k$, $f_k=f\chi_{B_k}$, and $h=\max\Bigl\{1,\sum\limits_{k=1}^K\chi_{B_k}\Bigr\}$. Then it's clear that $4\supp f\subset \Omega$, $\supp f_k\subset\overline{B_k}$, and
\[
1\leq h(x)\leq c_n,\qquad\text{for each } x\in\bb R^n,
\]
\[
\frac1{h}\sum\limits_{k=1}^{K}f_k=\frac1{\max\Bigl\{1,\sum\limits_{k=1}^K\chi_{B_k}\Bigr\}}\sum\limits_{k=1}^Kf\chi_{B_k}=f.
\]
Let $u_k:=R_t\left(\frac1hf_k\right)$. Since
\[
R_t^{-1}\Bigl(\sum\limits_{k=1}^Ku_k\Bigr)=\sum\limits_{k=1}^KR_t^{-1}u_k=\sum\limits_{k=1}^K\frac1hf_k=f,
\]
it follows that $u=R_tf=\sum\limits_{k=1}^Ku_k$ by the uniqueness of such a solution in $\m V$. Let $y_k$ be the center of $B_k$, so that $B_k=B\Big(y_k,\frac1{m(y_k,V_t+|\mathbf{B}|)}\Big)$. By (\ref{great}), for each $j\in\bb N$ and $k=1,\ldots,K$ we have
\begin{multline}
\int_{\bb R^n\backslash 4B_k}m(\cdot,V_t+|\mathbf{B}|)^2|u_k|^2e^{2\ep\varphi_j(\cdot,y_k)}\leq C\int_{4B_k\backslash2B_k}|u_k|^2\frac1{|\diam B_k|^2}e^{2\ep\varphi_j(\cdot,y_k)}\label{great0}.
\end{multline}
We note that for $x\in 4B_k\backslash2B_k$, by definition we have $|x-y_k|\leq\frac4{m(y_k,V_t+|\mathbf{B}|)}$, which implies that $d(x,y_k)\leq C$, and so by (\ref{closetod}) and Proposition \ref{dapprox1}, it follows that $\varphi_j(\cdot,y_k)\leq C$ on $4B_k\backslash2B_k$, with $C$ independent of $j,y_k$. Hence we have
\begin{multline}
\int_{\bb R^n\backslash 4B_k}m(\cdot,V_t+|\mathbf{B}|)^2|u_k|^2e^{2\ep\varphi_j(\cdot,y_k)}\leq Cm\big(y_k,V_t+|\mathbf{B}|\big)^2\int_{\bb R^n}|u_k|^2\\[4mm]\leq Cm\big(y_k,V_t+|\mathbf{B}|\big)^2\int_{\bb R^n}\frac1{h^2}|f_k|^2\leq C\int_{B_k}|f_k|^2m\big(\cdot,V_t+|\mathbf{B}|\big)^2\label{great1},
\end{multline}
where in the second inequality we used the uniform boundedness of the resolvents $R_t$ from $L^2(\bb R^n)$ into $L^2(\bb R^n)$, and in the third inequality we used the fact that $f_k$ is supported on $B_k$ and Lemma \ref{propertiesm}. Here $C$ is independent of $j,y_k$ and $t$. Letting $j\ra\infty$ on (\ref{great1}), using Fatou's Lemma  and (\ref{closetod}), we conclude that for each $k=1,\ldots,K$,
\begin{equation}\label{great2}
\int_{\bb R^n\backslash 4B_k}m(\cdot,V_t+|\mathbf{B}|)^2|u_k|^2e^{2\ep d(\cdot,y_k, V_t+|{\mathbf B}|)}\leq C\int_{B_k}|f_k|^2m\big(\cdot,V_t+|\mathbf{B}|\big)^2.
\end{equation}
We remark that, by using $L_t^{-1}$ instead of $R_t$ whenever it has appeared up to this point, we can prove all the results so far up to (\ref{great0}). In this case, note that
\begin{multline}\label{greathomcalc}
\Bigl(\int_{4B_k\backslash2B_k}\big|L_t^{-1}\frac1hf_k\big|^2\frac1{|\diam B_k|^2}\Bigr)^{\frac12}\leq C\big\Vert L_t^{-1}\frac1hf_k\big\Vert_{L^{2^*}(\bb R^n)}\leq C\big\Vert\D L_t^{-1}\frac1hf_k\big\Vert_{L^2(\bb R^n)}\\[3mm]\leq C\big\Vert L_t^{-1}\frac1hf_k\big\Vert_{\dot{\m V}}\leq C|\supp f_k|^{\frac1n}\Vert f_k\Vert_{L^2(\bb R^n)}\leq C\Bigl(\int_{B_k}|f_k|^2\frac1{m(\cdot,V_t+|{\mathbf B}|)^2}\Bigr)^{\frac12},
\end{multline}
where we first used H\"older's inequality, then the Sobolev inequality and the diamagnetic inequality, then the definition of the norm $\Vert\cdot\Vert_{\dot{\m V}}$, then (\ref{seqbounded}) and (\ref{embeddingdual}), and finally Lemma \ref{propertiesm} and the fact that $\supp f_k\subset B_k$. Using (\ref{greathomcalc}) and the analogous (\ref{great0}), we can thus prove
\begin{equation}\label{greathom}
\int_{\bb R^n\backslash 4B_k}m(\cdot,V_t+|\mathbf{B}|)^2\Big|L_t^{-1}\Bigl(\frac1hf_k\Bigr)\Big|^2e^{2\ep\varphi_j(\cdot,y_k)}\leq C\int_{B_k}|f_k|^2\frac1{m(\cdot,V_t+|{\mathbf B}|)^2},
\end{equation}
where $C$ here is independent of $t, k$ and $f_k$.

Let us prove that for each $x\in\bb R^n\backslash\Omega$,
\begin{equation}\label{eq.decayest}
\sum_{k=1}^Ke^{-2\ep d(x,y_k)}\leq C_{\ep},
\end{equation}
where $C_{\ep}$ is independent of $x$, $t$, $K$, and the support of $f$. Fix $x\in\bb R^n\backslash\Omega$ and let
\[
A_{\ell}=\Big\{y\in\bb R^n\,:\,|x-y|\in\Big[\frac{\ell}{m(x,V_t+|\mathbf{B}|)},\frac{\ell+1}{m(x,V_t+|\mathbf{B}|)}\Big)\Big\},\qquad\ell\in\bb Z, \ell\geq0.
\]
Now, since $x\notin 4B_k$ for any $k$, note that
\begin{equation}\label{eq.fix1}
\sum_{k=1}^Ke^{-2\ep d(x,y_k)}=\sum_{\ell=4}^{\infty}\sum_{y_k\in A_{\ell}}e^{-2\ep d(x,y_k)}\leq\sum_{\ell=4}^{\infty}\Big(\max_{y_k\in A_{\ell}}e^{-2\ep d(x,y_k)}\Big)\sum_{y_k\in A_{\ell}}1.
\end{equation}
Since for all $y_k\in A_{\ell}$ we have that $|x-y_k|\approx\frac{\ell}{m(x,V_t+|\mathbf{B}|)}$, it follows by using Lemma \ref{polygrowth} that $d(x,y_k)\gtrsim\ell^{\delta}$ for all $y_k\in A_{\ell}$ and for some $\delta>0$. Hence
\[
\max_{y_k\in A_{\ell}}e^{-2\ep d(x,y_k)}\leq e^{-c_{\ep}\ell^{\delta}},\qquad\text{for each }\ell\geq4.
\]
For each $\ell\geq4$ and $y\in A_{\ell}$, write $B_y=B(y,\frac1{m(y,V_t+|\mathbf{B}|)})$. Then
\begin{multline*}
	\sum_{y_k\in A_{\ell}}1=\sum_{y_k\in A_{\ell}}\frac{|B_k|}{|B_k|}\leq\frac1{\min_{y\in A_{\ell}}|B_y|}\sum_{y_k\in A_{\ell}}|B_k|\\[3mm]\lesssim\frac1{\min_{y\in A_{\ell}}|B_y|}\Big|B\Big(x,\frac{\ell}{m(x,V_t+|\mathbf{B}|)}\Big)\Big|\lesssim\max_{y\in A_{\ell}}[m(y,V_t+|\mathbf{B}|)^n]\frac{\ell^n}{m(x,V_t+|\mathbf{B}|)^n}\\[3mm]\lesssim(m(x,V_t+|\mathbf{B}|)^n\ell^{\frac32n})\frac{\ell^n}{m(x,V_t+|\mathbf{B}|)^n}\lesssim\ell^{3n}.
\end{multline*}
Therefore, from (\ref{eq.fix1}) we conclude that
\[
\sum_{k=1}^Ke^{-2\ep d(x,y_k)}\lesssim\sum_{\ell=4}^{\infty}\ell^{3n}e^{-c_{\ep}\ell^{\delta}}\lesssim C_{\ep},
\]
as claimed.
	
	 Next, since $4B_k\subset\Omega$, and $d(\cdot,\supp f,V_t+|{\mathbf B}|)\leq d(\cdot,y_k,V_t+|{\mathbf B}|)$, we note that for $\ep'=\frac12\ep$ and $\ep$ small enough,
\begin{multline}\label{great3}
\int_{\bb R^n\backslash\Omega}m(\cdot,V_t+|\mathbf{B}|)^2|u|^2e^{2\ep' d(\cdot,\supp f)}\\[4mm]=\int_{\bb R^n\backslash\Omega}\Big|\sum_{k=1}^Km(\cdot,V_t+|\mathbf{B}|)u_ke^{\ep'[d(\cdot,y_k)+d(\cdot,\supp f)]}e^{-\ep'd(\cdot,y_k)]}\Big|^2\\[4mm]\leq\int_{\bb R^n\backslash\Omega}\Big[\sum\limits_{k=1}^Km(\cdot,V_t+|\mathbf{B}|)^2|u_k|^2e^{2\ep d(\cdot,y_k)}\Big]\sum_{j=1}^Ke^{-2\ep' d(\cdot,y_j)}\\[4mm]\leq C\sum\limits_{k=1}^K~\int_{\bb R^n\backslash\Omega}m(\cdot,V_t+|\mathbf{B}|)^2|u_k|^2e^{2\ep d(\cdot,y_k)} \\[4mm]\leq C\sum\limits_{k=1}^K~\int_{\bb R^n\backslash4B_k}m(\cdot,V_t+|\mathbf{B}|)^2|u_k|^2e^{2\ep d(\cdot,y_k)}\\[3mm]\leq C\sum\limits_{k=1}^K~\int_{B_k}|f_k|^2m(\cdot,V_t+|\mathbf{B}|)^2\\[4mm]\leq C\int_{\bb R^n}|f|^2m(\cdot,V_t+|\mathbf{B}|)^2 ,
\end{multline}
where first we used the Cauchy-Schwartz inequality and (\ref{eq.decayest}), then (\ref{great2}). We note that $C$ does not depend on $\supp f$, nor on the subcollection $\n B''\in\n F$ used. Upon using Proposition \ref{dset}, (\ref{l2result-intro}) follows immediately.  Similarly, by using (\ref{greathom}) in place of (\ref{great2}) in the argument leading up to (\ref{great3}), it is clear that we achieve (\ref{l2resultunif}).\hfill{$\square$}\\


We remark that (\ref{l2result-intro}) implies a Gaffney-type estimate:
\begin{multline}\label{gaffney}
\Big\Vert(I+t^2L)^{-1}f\Big\Vert_{L^2(E)}\\[3mm]\leq Ct^2\exp\Big(-\ep\frac{\dist(E,F)}{t}\Big)\Big\Vert m\Big(\cdot,V+|\mathbf{B}|+\frac {1}{t^2}\Big)f\Big\Vert_{L^2(F)}
\end{multline}
where $f\in L^2(\bb R^n)$ is compactly supported with $\supp f\subset F$, and $E\subset\bb R^n$ satisfies that $d(E,F,V+|{\mathbf B}|)$ is large enough (here and elsewhere, $\dist(\cdot,\cdot)$ denotes the Euclidean distance). We prove (\ref{gaffney}) in the following corollary, which also includes a proof of (\ref{l2resulthom-intro}).

\begin{corollary}\label{l2decayhom} Assume that $\mathbf{a}\in L^2_{loc}(\bb R^n)$, $A$ is an elliptic matrix with complex, bounded, measurable coefficients, and $V\in L^1_{loc}(\bb R^n)$ with $V\geq0$ a.e. on $\bb R^n$. If $\mathbf{a}\equiv0$, assume $V\in RH_{\frac n2}$, otherwise take assumptions (\ref{magassumptions-intro}). Suppose $f\in L^2(\bb R^n)$ is compactly supported. Then there exists $\tilde d>0$, depending on $\Vert V+|\mathbf{B}|\Vert_{RH_{\frac n2}}$ and $n$ only, and there exists $\ep>0$ such that (\ref{l2resulthom-intro}) holds,
where $\ep, C$ depend on $\lambda,\Lambda, n,\ep,\Vert V+|\mathbf{B}|\Vert_{RH_{\frac n2}}$, and the constants from (\ref{magassumptions-intro}), and are independent of $\supp f$. Moreover, suppose $E,F$ are open sets in $\bb R^n$ where $\supp f\subset F$ and $d(E\,,F\,,V+|{\mathbf B}|)>\tilde d$. Then for each $t>0$, (\ref{gaffney}) holds with constants $\ep, C$ depending only on $\lambda,\Lambda, n,\ep,\Vert V+|\mathbf{B}|\Vert_{RH_{\frac n2}}$, and the constants from (\ref{magassumptions-intro}). 
\end{corollary}

\noindent\emph{Proof.} Using (\ref{tmdomination}) and (\ref{tddomination}) on  (\ref{l2resultunif}), we observe that
\begin{multline}
\int_{\big\{x\in\bb  R^n|\,d(x,\supp f,V+|\mathbf{B}|)\geq\tilde d\big\}}m\left(\cdot,V+|\mathbf{B}|\right)^2\left|L_t^{-1}f\right|^2e^{2\ep d(\cdot,\supp f,V+|\mathbf{B}|)}\\\leq C\int_{\bb R^n}|f|^2\frac1{m(x,V_t+|{\mathbf B}|)^2}.\label{corproof1}
\end{multline}
The right-hand side of (\ref{corproof1}) converges to the right-hand side of (\ref{l2resulthom-intro}) as $t\ra\infty$ by the Lebesgue Monotone Convergence Theorem, since it can be proven that for each $x\in\bb R^n$,
\[
m(x,V_t+|{\mathbf B}|)\searrow m(x,V+|{\mathbf B}|),\qquad\text{as }t\ra\infty.
\]
Therefore, we can now achieve (\ref{l2resulthom-intro}) by using Lemma \ref{opstrongconv} and Fatou's Lemma on the left-hand side of (\ref{corproof1}).

To see that (\ref{gaffney}) is true, first note that for all $x\in\bb R^n$, if we let $r_x=\frac1{m(x,V_t+|{\mathbf B}|)}$, then
\[
1=\frac1{r_x^{n-2}}\int_{B(x,r_x)}V+\frac1{t^2}\geq\frac{r_x^2}{t^2}|B(0,1)|,
\]
and so
\[
\frac1{m(x,V_t+|{\mathbf B}|)}\leq\frac1{|B(0,1)|^{\frac12}}\,t.
\]
Using the above fact, (\ref{l2result-intro}), (\ref{tddomination}), and the hypothesis on the sets $E,F$, we can write
\begin{equation*}
\int_E|(I+t^2L)^{-1}f|^2e^{2\ep\frac{\dist(x,F)}{t}}\,dx\leq Ct^2\int_F|f|^2m(x,V_t+|{\mathbf B}|)^2\,dx,
\end{equation*}  
and from the above inequality, it is clear that (\ref{gaffney}) follows.\hfill{$\square$}

%
%
%
%
%
%
%
%
%
%
}

\section{The fundamental solution of the magnetic Schr\"odinger operator and its properties}\label{fundsolnsec}

In this section we aim to pass to the pointwise estimates.

\subsection{The semigroup theory, Kato-Simon inequality, and the heat kernel}

We use all the notation and definitions from Section \ref{theory}. Let $\mathbf{a}\in L^2_{loc}(\bb R^n)$, let $A$ be an elliptic matrix with complex, bounded measurable coefficients, and let $V\in L^1_{loc}(\bb R^n)$ satisfy (\ref{gkcond}) and (\ref{imVreq}) with $c_2\geq0$ and $c_4$ either $0$ or $1$. First, recall the Resolvent formula (see the remarks following Theorem 1.43 in \cite{Ou}):
\begin{equation}\label{resolvent}
(L+\ep)^{-1}f=\int_0^{\infty}e^{-\ep t}e^{-tL}f\,dt,\qquad\text{for each }\ep>0,~ f\in H=L^2(\bb R^n)
\end{equation}
and the following identity 
\begin{equation}\label{semigrouplimit}
e^{-tL}f=\lim\limits_{m\ra\infty}\Bigl(1+\frac tmL\Bigr)^{-m}f,\qquad\text{for each } f\in H
\end{equation}
where the limits are in the topology of $H$ (that is, in the $L^2(\bb R^n)-$sense).
Define the form
\[
\f b(u,v)=\int_{\bb R^n}\nabla u\cdot\overline{\nabla v},\qquad\text{for each } u,v\in D(\f b)=W^{1,2}(\bb R^n).
\]
In a very similar fashion as in Section \ref{theory}, it is easy to see that $\f b$ is densely defined, continuous, closed, and accretive. Accordingly, we can define the unbounded operator $-\Delta:D(-\Delta)\ra H$, where $D(-\Delta)$ is given as
\[
D(-\Delta)=\Big\{u\in D(\f b) \text{ s.t. }\exists v\in H: \f b(u,\phi)=(v,\phi)_H~\forall\phi\in D(\f b)\Big\},\qquad -\Delta u:=v.
\]
Hence $D(-\Delta)$ is dense in $W^{1,2}(\bb R^n)$, $(-\Delta+\ep)^{-1}$ is invertible on $H$ for every $\ep>0$, and there is a strongly continuous contraction semigroup $e^{t\Delta}$ associated to $-\Delta$. Accordingly, identities analogous to (\ref{resolvent}) and (\ref{semigrouplimit}) hold. Furthermore, the operator $-\Delta$ is known to enjoy several more properties: its heat semigroup is given by integration against a non-negative heat kernel $p_{-\Delta}$; that is, $p_{-\Delta}(\cdot,\cdot~;t)$ is a measurable function on $\bb R^n\times\bb R^n$ such that
\[
e^{t\Delta}f(x)=\int_{\bb R^n}p_{-\Delta}(x,y;t)f(y)\,dy,\qquad\text{ for a.e. } x\in\bb R^n~\text{and for every }t>0, f\in L^2(\bb R^n),
\]
and
\[
|p_{-\Delta}(x,y;t)|\leq Ct^{-n/2}\qquad\text{for all }t>0.
\]
The operator $-\Delta$ also has a homogeneous realization, which by a slight abuse of notation we denote as $-\Delta$. The \emph{fundamental solution} of $-\Delta$, $\Gamma_{-\Delta}$, is well-known to exist and to satisfy
\[
\frac{c_n}{|x-y|^{n-2}}=\Gamma_{-\Delta}(x,y)=\int_0^{\infty}p_{-\Delta}(x,y;t)\,dt,
\]
where $c_n$ is a dimensional constant. For each compactly supported $f\in L^{\infty}(\bb R^n)$, we can write
\[
(-\Delta)^{-1}f(x)=\int_{\bb R^n}\Gamma_{-\Delta}(x,y)f(y)\,dy,\qquad\text{a.e. }x\in\bb R^n.
\]
From the non-negativity of the heat kernel $p_{-\Delta}$ and the resolvent formula, we deduce that if $0\leq f\leq g$ with $f,g\in L^2(\bb R^n)$, then for all $\ep>0$ one has
\begin{equation}\label{nonnegdom}
(-\Delta+\ep)^{-1}f\leq(-\Delta+\ep)^{-1}g
\end{equation}
in the a.e. sense on $\bb R^n$. 

In the case that $A\equiv I$, ${\mathbf a}\in L^2_{loc}(\bb R^n)$ is real-valued, and $V\in L^1_{loc}(\bb R^n)$  with $V\geq0$ a.e. on $\bb R^n$, we have the Kato-Simon inequality, the following formulation of which can be found in \cite{LS}, Lemma 6:
\begin{equation}\label{katosimon}
|(L+\ep)^{-1}f|\leq(-\Delta+\ep)^{-1}|f|,\qquad\text{for each } f\in H=L^2(\bb R^n).
\end{equation}
A function $p_L(x,y;t):\bb R^n\times\bb R^n\times\bb R_+\ra\bb R$ is called a \emph{heat kernel} of the semigroup $e^{-tL}$ if for each $t>0$, $p_L(\cdot,\cdot~;t)$ is a measurable function on $\bb R^n\times\bb R^n$, and
\[
e^{-tL}f(x)=\int_{\bb R^n}p_{L}(x,y;t)f(y)\,dy,\qquad\text{ for a.e. } x\in\bb R^n~\text{and for every }t>0, f\in L^2(\bb R^n).
\]
It is clear that if a heat kernel to the semigroup $e^{-tL}$ exists, then it must be unique. Let us prove
\begin{proposition} Suppose that $\mathbf{a}\in L^2_{loc}(\bb R^n)$, $A\equiv I$, and $V\in L^1_{loc}(\bb R^n)$ with $V$ real-valued and $V\geq0$ a.e. on $\bb R^n$. Then
\begin{equation}\label{semigroupdom}
|e^{-tL}f|\leq e^{t\Delta}|f|,\qquad\text{for each } t>0, f\in H.
\end{equation}
Furthermore, $e^{-tL}$ can be seen as a bounded map from $L^1(\bb R^n)$ into $L^{\infty}(\bb R^n)$. The heat kernel $p_L$ of the semigroup $e^{-tL}$ exists, and moreover it satisfies
\begin{equation}\label{heatkerneldom}
|p_L(x,y;t)|\leq p_{-\Delta}(x,y;t)
\end{equation}
for all $t>0$, a.e. $x\in\bb R^n$, a.e. $y\in\bb R^n$.
\end{proposition} 

\noindent\emph{Proof.} It is immediate that (\ref{katosimon}) can be rewritten as
\[
\Bigl|\Bigl(\frac1{\ep}L+1\Bigr)^{-1}f\Bigr|\leq\Bigl(\frac1{\ep}(-\Delta)+1\Bigr)^{-1}|f|,
\]
for $\ep>0$ and $f\in H$. Using the observation (\ref{nonnegdom}) and the Kato-Simon inequality, we note that
\[
\Bigl|\Bigl(\frac1{\ep}L+1\Bigr)^{-m}f\Bigr|\leq\Bigl(\frac1{\ep}(-\Delta)+1\Bigr)^{-1}\Bigl|\Bigl(\frac1{\ep}L+1\Bigr)^{-m+1}f\Bigr|\leq\cdots\leq\Bigl(\frac1{\ep}(-\Delta)+1\Bigr)^{-m}|f|
\]
for any $m\in\bb N$, $\ep>0$ and $f\in H$. It thus follows that
\[
\Bigl|\Bigl(\frac tmL+1\Bigr)^{-m}f\Bigr|\leq\Bigl(\frac tm(-\Delta)+1\Bigr)^{-m}|f|
\]
for all $t>0$, all $m\in\bb N$ and all $f\in H$. We take limit as $m\ra\infty$ using (\ref{semigrouplimit}), and (\ref{semigroupdom}) follows.

Now, we'll show $e^{-tL}$ can be seen as a bounded operator mapping $L^1(\bb R^n)$ into $L^{\infty}(\bb R^n)$. To see this, note that for $f\in C_c^{\infty}(\bb R^n)$, we have by (\ref{semigroupdom}),
\[
|e^{-tL}f(x)|\leq e^{t\Delta}|f|(x)\leq\Vert e^{t\Delta}|f|\Vert_{L^{\infty}(\bb R^n)}\leq Ct^{-\frac n2}\Vert f\Vert_{L^1(\bb R^n)},\qquad\text{a.e. }x\in\bb R^n
\]
whence
\[
\Vert e^{-tL} f\Vert_{L^{\infty}(\bb R^n)}\leq Ct^{-\frac n2}\Vert f\Vert_{L^1(\bb R^n)},\qquad\text{for each } f\in C_c^{\infty}(\bb R^n),
\]
and therefore, using the density of $C_c^{\infty}(\bb R^n)$ in $L^1(\bb R^n)$, it follows that we can see $e^{-tL}$ as a map of $L^1(\bb R^n)$ into $L^{\infty}(\bb R^n)$, and a similar argument gives that $e^{-tL}:L^2(\bb R^n)\ra L^{\infty}(\bb R^n)$ is also a bounded map, with
\[
\Vert e^{-tL}\Vert_{L^2\ra L^{\infty}}\leq Ct^{-\frac n4},\qquad\text{for each } t>0.
\]
By Dunford's Theorem (see \cite{dp}), this implies the existence of a measurable function $p_L(\cdot,\cdot,t)$ on $\bb R^n\times\bb R^n$, which satisfies
\begin{equation}\label{heatkernelrep}
e^{-tL}f(x)=\int_{\bb R^n}p_L(x,y;t)f(y)\,dy,\qquad\text{for each } t>0, \text{a.e. } x\in\bb R^n, f\in L^2(\bb R^n).
\end{equation}
At this point we show the domination of $p_L$ by $p_{-\Delta}$. Fix $(x,t)\in\bb R^n\times(0,\infty)$ such that $p_L(x,\cdot\,;t)$ is measurable in $\bb R^n$. Suppose there exists $B\subset\bb R^n$ such that
\begin{equation}\nonumber
|p_L(x,y;t)|>p_{-\Delta}(x,y;t),\qquad y\in B.
\end{equation}
Using (\ref{semigroupdom}) and (\ref{heatkernelrep}) with $f=\chi_B\frac{\overline{p_L}}{|p_L|}$, it is readily observed that for a.e. $x\in\bb R^n$, all $t>0$,
\[
\int_B|p_L(x,y;t)|\,dy\leq\int_Bp_{-\Delta}(x,y;t)\,dy,
\]
and so it follows that $|B|=0$. More precisely, for all $t>0$, for almost every $x\in\bb R^n$, for almost every $y\in\bb R^n$, the inequality (\ref{heatkerneldom}) holds.\hfill{$\square$}\\

\subsection{The fundamental solution of the magnetic Schr\"odinger operator, part I}

In the following definition, $\n L$ is either $L+\ep$ for $\ep>0$ or $\dot L$, with $\n H,\n R$ the domain and range respectively of these operators. Specifically, when $\n L=L+\ep$, we write $\n H=D(L)$, $\n R=L^2(\bb R^n)$. When $\n L=\dot L$, we write $\n H=\dot{\m V}$, $\n R=\dot{\m V}'$.

\begin{definition}\label{fundsoln} We say that a measurable function $\Gamma(x,y)$ defined on $\bb R^n\times\bb R^n$ is the fundamental solution to the invertible operator $\n L:\n H\ra\n R$ if the following conditions are satisfied:
\begin{enumerate}
\item For each $f\in C_c^{\infty}(\bb R^n)$, the function
\[
u_f=\int_{\bb R^n}\Gamma(\cdot,y)f(y)\,dy
\]
lies in $\n H$ and satisfies $u_f=\n L^{-1}f$.
\item For a.e. $y\in\bb R^n$, $\Gamma(\cdot,y)$ solves $\n L\Gamma=0$ in the weak sense locally on $\bb R^n\backslash\{y\}$.
\end{enumerate}

\end{definition}
We carefully note that we avoid for now a stronger traditional statement that $\n L\Gamma=\delta$ in the sense of distributions. We will discuss this further below.

At this point, for $\ep\geq0$ note due to (\ref{heatkerneldom}) that
\begin{equation}\label{justifyabs}
\int_0^{\infty}e^{-\ep t}|p_L(x,y;t)|\,dt\leq\int_0^{\infty}|p_{-\Delta}(x,y;t)|\,dt=\int_0^{\infty}p_{-\Delta}(x,y;t)\,dt=\Gamma_{-\Delta}(x,y)<\infty
\end{equation}
for a.e. $x,y\in\bb R^n, x\neq y$. For each $\ep\geq0$, we define the measurable function $\Gamma_{L+\ep}(x,y)$ given by
\begin{equation}\label{fundsoldef}
\Gamma_{L+\ep}(x,y):=\int_0^{\infty}e^{-\ep t}p_L(x,y;t)\,dt,
\end{equation}
and we will eventually see that this function is the fundamental solution to the operator $L+\ep$ (we use the notation $\Gamma_L$ for $\ep=0$; $\Gamma_L$ will be seen to be the fundamental solution to the operator $\dot{L}$). Due to (\ref{justifyabs}), a function given by (\ref{fundsoldef}) is well-defined and finite a.e.. The following lemma captures the expected convergence result:

\begin{lemma}\label{pwconv} Suppose that $\mathbf{a}\in L^2_{loc}(\bb R^n)$, $A\equiv I$, and $V\in L^1_{loc}(\bb R^n)$ with $V$ real-valued and $V\geq0$ a.e. on $\bb R^n$. Then for every compactly supported $f\in L^{\infty}(\bb R^n)$,
\begin{equation}\label{convop}
(L+\ep)^{-1}f(x)\longrightarrow\int_{\bb R^n}\Gamma_L(x,y)f(y)\,dy,\quad\text{for a.e. }x\in\bb R^n,
\end{equation}
and
\begin{equation}\label{fundsolconv}
\Gamma_{L+\ep}(x,y)\longrightarrow\Gamma_L(x,y)\quad\text{for a.e. }x,y\in\bb R^n, x\neq y,
\end{equation}
as $\ep\searrow0$.
\end{lemma}

\noindent\emph{Proof.} Fix $f\in L^{\infty}(\bb R^n)$ with compact support. For each $\ep\geq0$, we have
\begin{equation}\label{fubiniarg}
\int_0^{\infty}\int_{\bb R^n}e^{-\ep t}|p_L||f|\,dy\,dt\leq\int_0^{\infty}\int_{\bb R^n}e^{-\ep t}p_{-\Delta}|f|\,dy\,dt\leq\int_{\bb R^n}\Gamma_{-\Delta}|f|\,dy=(-\Delta)^{-1}|f|,
\end{equation}
where in the second inequality we used Tonelli's Theorem. Since $(-\Delta)^{-1}|f|$ is a measurable finite a.e. function, it follows that for almost every $x\in\bb R^n$ and every $\ep>0$, Fubini's Theorem can be applied to (\ref{resolvent}) when (\ref{heatkernelrep}) is used. Hence, using (\ref{fundsoldef}), we arrive at the identity
\begin{equation}\label{fundsolisakernel}
(L+\ep)^{-1}f(x)=\int_{\bb R^n}\Gamma_{L+\ep}(x,y)f(y)\,dy,
\end{equation}
valid a.e. in $\bb R^n$, for $\ep>0$ and $f\in L^{\infty}(\bb R^n)$ with compact support. It is clear that
\[
e^{-\ep t}p_{L}(x,y;t)f(y)\longrightarrow p_{L}(x,y;t)f(y)
\]
pointwise in $(y,t)$ for almost every $x\in\bb R^n$, as $\ep\searrow0$. Moreover, note that for a.e. $y\in\bb R^n$,
\[
|e^{-\ep t}p_{L}(x,y;t)f(y)|\leq|p_{L}(x,y;t)f(y)|\leq p_{-\Delta}(x,y;t)|f(y)|
\]
and, since for almost every $x\in\bb R^n$,
\[
\int_{\bb R^n}\int_0^{\infty}p_{-\Delta}(x,y;t)|f(y)|\,dt\,dy<+\infty,
\]
then for almost every $x\in\bb R^n$, we can apply the Lebesgue Dominated Convergence Theorem to get
\begin{equation}\nonumber
\int_{\bb R^n}\int_0^{\infty}e^{-\lambda_nt}p_{-H}(x,y;t)f(y)\,dt\,dy\longrightarrow\int_{\bb R^n}\int_0^{\infty}p_{-H}(x,y;t)f(y)\,dt\,dy,
\end{equation}
which is (\ref{convop}). A very similar convergence argument delivers (\ref{fundsolconv}).\hfill{$\square$}\\

Combining the results of this section with Lemma \ref{opstrongconv}, we have

\begin{theorem}\label{firstone}Suppose $\mathbf{a}\in L^2_{loc}(\bb R^n)$, $A\equiv I$, $V\in L^1_{loc}(\bb R^n)$ with $V\geq0$ a.e. on $\bb R^n$. Then the function $\Gamma_L(x,y)$ given in (\ref{fundsoldef}) satisfies the following properties:
\begin{enumerate}[a)]
\item For any $f\in L^{\infty}(\bb R^n)$ with compact support, the function defined by
\[
u_f=\int_{\bb R^n}\Gamma_L(\cdot,y)f(y)\,dy
\]
lies in $\dot{\m V}$ and is the unique element in $\dot{\m V}$ satisfying $Lu_f=f$ in the sense of distributions on $\bb R^n$. Moreover the sequence $\{(L+\ep)^{-1}f\}$ converges pointwise almost everywhere on $\bb R^n$, and strongly in $\dot{\m V}$, to $u_f$ as $\ep\searrow0$.
\item There exists a constant $C$ depending only on $n,\lambda,\Lambda$, such that for a.e. $x,y\in\bb R^n$, 
\begin{equation}\label{pointwiseupper}
|\Gamma_L(x,y)|\leq\frac{C}{|x-y|^{n-2}}.
\end{equation}
\item For a.e. $y\in\bb R^n$, $\Gamma_L(\cdot,y)\in L^1_{loc}(\bb R^n)$, and $\Gamma_L(\cdot,y)\in L^{\infty}_{loc}(\bb R^n\backslash\{y\})$.
\end{enumerate}
\end{theorem}

\emph{Proof.} Statements b) and c) hold by the definition of $\Gamma_L$ and (\ref{heatkerneldom}). Statement a) holds by (\ref{convop}) and Lemma \ref{opstrongconv}.\hfill{$\square$}

Theorem \ref{firstone} doesn't yet give the existence of a fundamental solution, but it does give the existence of an integral kernel of the operator $\dot L^{-1}$. Missing from the theorem is another important aspect of the fundamental solution, which is that $L\Gamma_L(\cdot,y)=\delta_y$ in the sense of distributions on $\bb R^n$. Though this latter fact may not be accessible to us in the full generality, we will need some variation of $L\Gamma_L(\cdot,y)=0$ in the weak sense on $\bb R^n\backslash\{y\}$ to satisfy the conditions of Definition \ref{fundsoln}. For this purpose, it is necessary to prove $\D\Gamma_L(\cdot,y)\in L^2_{loc}(\bb R^n\backslash\{y\})$.

At this point, we can prove an important property of local weak solutions to the operator $L$; namely, the \emph{local uniform boundedness} of weak solutions to the operator $L$, also known as a \emph{Moser} estimate. We capture this result in
\begin{theorem}\label{moser}Assume that $\mathbf{a}\in L^2_{loc}(\bb R^n)$, $A\equiv I$, and $V\in L^1_{loc}(\bb R^n)$ with $V\geq0$ a.e. on $\bb R^n$. Let $B\equiv B(x_0,R)\subset\bb R^n$, $f\in L^{\infty}(B)$, and suppose $u$ solves $Lu=f$ in the weak sense on $2B:=B(x_0,2R)$. Then $u$ is locally essentially bounded, and
	\[
	\Vert u\Vert_{L^{\infty}(\frac14B)}\leq C\Bigl\{\Bigl(~\fint_{2B}|u|^2\Bigr)^\frac12+R^2\Bigl(~\fint_{2B}|f|^2\Bigr)^{\frac12}\Bigr\}.
	\]
\end{theorem}

\noindent\emph{Proof.} This result is proven in Lemma 1.3 of \cite{S5} for magnetic Schr\"odinger operators with potentials $V+\ep$ where $\ep>0$ is a constant. The proof for the magnetic Schr\"odinger operators satisfying the hypothesis of the theorem as stated here follows as soon as one establishes a Kato-Simon inequality for such operators. More precisely, we want to prove that
\begin{equation}\label{homogeneouskatosimon}
|\dot L^{-1}f|(x)\leq\int_{\bb R^n}\Gamma_{-\Delta}(x,y)|f(y)|\,dy,\qquad\text{for a.e. }x\in\bb R^n,
\end{equation}
is true for each $f\in L^{2}(\bb R^n)$ with compact support (note that (\ref{homogeneouskatosimon}) should make sense even when $\int_{\bb R^n}\Gamma_{-\Delta}(x,y)f(y)\,dy=+\infty$). We note that if $f\in L^{\infty}(\bb R^n)$ with compact support, then we can prove (\ref{homogeneouskatosimon}) as follows: Recall that
\[
(-\Delta+\ep)^{-1}|f|(x)\leq\int_{\bb R^n}\Gamma_{-\Delta}(x,y)|f(y)|\,dy,\qquad\text{for each }\ep>0, \text{ a.e. } x\in\bb R^n.
\]
Using the Kato-Simon inequality (\ref{katosimon}) for operators $L+\ep$ and the previous estimate, we can write that
\[
\Big|(L+\ep)^{-1}f(x)\Big|\leq\int_{\bb R^n}\Gamma_{-\Delta}(x,y)|f(y)|\,dy,\qquad\text{for each }\ep>0, \text{ a.e. }x\in\bb R^n.
\]
Consequently, using a) in Theorem \ref{firstone} and taking $\liminf$ as $\ep\searrow0$ in the above inequality yields 
\[
|\dot L^{-1}f(x)|\leq\int_{\bb R^n}\Gamma_{-\Delta}(x,y)|f(y)|\,dy,\qquad\text{ a.e. }x\in\bb R^n,
\]
which establishes (\ref{homogeneouskatosimon}) in this case. Now suppose $f\in L^2(\bb R^n)$ has compact support, and for each $k\in\bb N$ let
\[
f_k(x):=\left\{\begin{matrix}f(x),&\text{if }|f(x)|\leq k\\k,&\text{if }|f(x)|>k.\end{matrix}\right.
\]
We note that for each $k\in\bb N$, $f\in L^{\infty}(\bb R^n)$, $\supp f_k=\supp f$, $f_k\ra f$ in $L^2(\bb R^n)$ as $k\ra\infty$, and
\[
|f_k(x)|\nearrow|f(x)|,\qquad\text{for a.e. }x\in\bb R^n\text{ as }k\ra\infty.
\]
Observe that for each $k\in\bb N$ and a.e. $x\in\bb R^n$, 
\begin{equation}\label{topass}
|\dot L^{-1}f_k(x)|\leq\int_{\bb R^n}\Gamma_{-\Delta}(x,y)|f_k(y)|\,dy\leq\int_{\bb R^n}\Gamma_{-\Delta}(x,y)|f(y)|\,dy.
\end{equation}
Since $f_k\ra f$ in $L^2(\bb R^n)$ and $\supp f_k=\supp f$, it follows that $f_k\ra f$ in the topology of $\dot{\m V}'$. Hence $\dot L^{-1}f_k\ra\dot L^{-1}f$ in $\dot{\m V}$ as $k\ra\infty$, since $\dot L^{-1}$ is a continuous operator. In particular, a subsequence of $\{\dot L^{-1}f_k\}$ converges to $\dot L^{-1}f$ pointwise a.e. in $\bb R^n$. Hence, passing to infinity along this subsequence in (\ref{topass}) implies the desired estimate (\ref{homogeneouskatosimon}).

With (\ref{homogeneouskatosimon}) at hand, the proof in \cite{S5} can be retraced to prove the theorem in the desired generality.\hfill{$\square$}\\

%

%
%

\subsection{The fundamental solution of the magnetic Schr\"odinger operator with smooth coefficients}

Suppose that $\mathbf{a}\in C_c^{\infty}(\bb R^n)$, $A$ is an elliptic matrix with complex, bounded measurable coefficients, and $V\in C^{\infty}(\bb R^n)\cap L^{\infty}(\bb R^n)$ satisfies (\ref{imVreq}), (\ref{gkcond}) with $c_2\equiv c_4\equiv0$. In this case, the elements of $D(\f l)$ coincide with those of $W^{1,2}(\bb R^n)$. To see this, suppose $u\in W^{1,2}(\bb R^n)$. Then by definition $u\in L^2(\bb R^n)$, thus $|\Re e\,V|^{\frac12}u\in L^2(\bb R^n)$, and the expression $\D u=\nabla u-i\mathbf{a}u$ clearly lies in $L^2(\bb R^n)$, showing $W^{1,2}(\bb R^n)\subset D(\f l)$. On the other hand, if $u\in D(\f l)$, then $u\in L^2(\bb R^n)$ by definition, and it can easily be shown that $\D u\in L^2(\bb R^n)$; hence $\nabla u\in L^2(\bb R^n)$, so that $D(\f l)\subset W^{1,2}(\bb R^n)$. In a similar way, the elements of $\dot{\m V}$ can be shown to coincide with $Y^{1,2}$, which is the space of elements of $\dot W^{1,2}(\bb R^n)$ which lie also in $L^{2^*}(\bb R^n)$.

Next, further assume that $A\equiv I$ and that $V\geq0$. Then we recover all the previous results; so, for instance we have the results of Lemma \ref{opstrongconv} and Theorem \ref{firstone}. Moreover, we can prove

\begin{theorem}\label{testfundsol} Assume $\mathbf{a}\in C_c^{\infty}(\bb R^n)$, $A\equiv I$, and $V\in C^{\infty}(\bb R^n)\cap L^{\infty}(\bb R^n)$ is a non-negative function. Then the function $\Gamma_L(x,y)$ from Theorem \ref{firstone} is the fundamental solution to the operator $\dot L$, and moreover enjoys the following properties:
\begin{enumerate}[a)]
\item For fixed $y\in\bb R^n$, $\Gamma_L(\cdot,y)\in W^{1,2}_{loc}(\bb R^n\backslash\{y\})$, and $\D\cdot\D\Gamma_L~\in L^{\infty}_{loc}(\bb R^n\backslash\{y\})$.
\item For fixed $y\in\bb  R^n$, the equation $L\Gamma_L=\delta_y$ holds in the sense of distributions. That is, for each $\phi\in C_c^{\infty}(\bb R^n)$,
\[
\int_{\bb R^n}\left[ \Gamma_L(x,y)\overline{\D\cdot\D\phi(x)} +V(x)\Gamma_L(x,y)\overline{\phi(x)}\right]\,dx=\overline{\phi(y)}.
\]
\item For almost every $y\in\bb R^n$, the equation
\begin{equation}\label{pointwise}
\D\cdot\D\Gamma_L(x,y)+V(x)\Gamma_L(x,y)=0
\end{equation}
holds (in particular, in the sense of distributions on $\bb R^n\backslash\{y\}$).
\item The identity
\begin{equation}\label{adjointcompare}
\Gamma_L(x,y)\equiv\overline{\Gamma_{L^*}(y,x)}
\end{equation}
holds in the a.e. sense on $\bb R^n\times\bb R^n$.
\end{enumerate}
\end{theorem}

\emph{Proof.} Let $f\in C_c^{\infty}(\bb R^n)$ and $\phi\in C_c^{\infty}(\bb R^n)$. By Theorem \ref{firstone}, it follows that
\[
u_f=\int_{\bb R^n}\Gamma_L(\cdot,y)f(y)\,dy
\]
lies in $Y^{1,2}$ and solves $Lu_f=f$ in the weak sense. That is, for $\phi\in C_c^{\infty}(\bb R^n)$ we have
\begin{align*}
\int_{\bb R^n}\Bigl[\D\Bigl(~\int_{\bb R^n}\Gamma_L(x,y)f(y)\,dy\Bigr)\overline{\D\phi(x)}+V(x)\Bigl(~\int_{\bb R^n}\Gamma_L(x,y)f(y)\,dy\Bigr)\overline{\phi(x)}\Bigr]\,dx\\[4mm]=\int_{\bb R^n}f(x)\overline{\phi(x)}\,dx,
\end{align*}
which we can write as
\begin{align*}
\int_{\bb R^n}\Bigl[\int_{~\bb R^n}\Gamma_L(x,y)f(y)\,dy\overline{\D\cdot \D\phi(x)}+\int_{\bb R^n}V(x)\Gamma_L(x,y)f(y)\overline{\phi(x)}\,dy\Bigr]\,dx\\[4mm]=\int_{\bb R^n}f(x)\overline{\phi(x)}\,dx.
\end{align*}
Using Fubini's Theorem, which is valid due to our assumptions on $\mathbf{a}$ and $V$, the above equality becomes
\begin{align*}
\int_{\bb R^n}f(y)\int_{\bb R^n}\left[ \Gamma_L(x,y)\overline{\D\cdot \D\phi(x)}+V(x)\Gamma_L(x,y)\overline{\phi(x)}\right]\,dx\,dy\\[4mm]=\int_{\bb R^n}f(x)\overline{\phi(x)}\,dx.
\end{align*}
Since $f\in C_c^{\infty}(\bb R^n)$ was arbitrary, the equality
\begin{equation}\label{weakerfund}
\int_{\bb R^n}\left[ \Gamma_L(x,y)\overline{\D\cdot \D\phi(x)} +V(x)\Gamma_L(x,y)\overline{\phi(x)}\right]\,dx=\overline{\phi(y)}
\end{equation}
holds for almost every $y\in\bb R^n$. This proves property b). By letting $h(y)=\dot L^*\phi(y)$ on $\bb R^n$, we see that $h\in C_c^{\infty}(\bb R^n)$. The invertibility of $\dot L^*$ implies that $\phi=(\dot L^*)^{-1}h$, and hence we can rewrite (\ref{weakerfund}) to see that
\[
\overline{(\dot L^*)^{-1}h(y)}=\int_{\bb R^n}\Gamma_L(x,y)h(x)\,dx.
\]
In particular, we have shown that for a.e. $y\in\bb R^n$,
\begin{equation}\label{uniquenessproof}
\int_{\bb R^n}\Big[\Gamma_L(x,y)-\overline{\Gamma_{L^*}(y,x)}\Big]h(x)\,dx=0,\qquad\text{for each } h\in\dot{L}^*(C_c^{\infty}(\bb R^n)).
\end{equation}
Since $C_c^{\infty}(\bb R^n)$ is dense in $\dot{\m V}$, then the bounded operator $\dot L^*:\dot{\m V}\ra\dot{\m V}'$ maps $C_c^{\infty}(\bb R^n)$ to a dense set in $\dot{\m V}'$.  Consequently, (\ref{adjointcompare}) follows.

Now fix $y\in\bb R^n$, $\Omega\subset\bb R^n$ a bounded open set, $r>0$, call $\n O:=\Omega\backslash B(y,r)$, and let $\phi\in C_c^{\infty}(\n O)$. Then, from (\ref{weakerfund}) it follows that for almost every such $y$,
\begin{equation}\label{trick}
\int_{\bb R^n} \Gamma_L(x,y)\overline{\D\cdot \D\phi(x)}\,dx=-\int_{\bb R^n}V(x)\Gamma_L(x,y)\overline{\phi(x)}\,dx,
\end{equation}
whence we have
\begin{equation}\label{funcestimate}
\Bigl|~\int_{\n O} \Gamma_L(x,y)\overline{\D\cdot\D\phi(x)}\,dx\Bigr|\leq C\Vert\phi\Vert_{L^1(\n O)}
\end{equation}
where $C$ depends on $V$ and $\n O$ but does \emph{not} depend on $\phi$. So consider the distribution
\[
\D\cdot\D\Gamma_L,
\]
defined on $C_c^{\infty}(\n O)$. Extend its definition to $L^1(\n O)$ using (\ref{trick})-(\ref{funcestimate}). By (\ref{funcestimate}) we observe that this distribution is a bounded linear functional on $L^1(\n O)$. Therefore, it actually is a measurable function living in $L^{\infty}(\n O)$ and, per (\ref{trick}), (\ref{pointwise}) follows.

To see that $\Gamma_L$ satisfies the conditions of Definition \ref{fundsoln}, it remains to show that for each $y\in\bb R^n$, $\D\Gamma_L\in L^2_{loc}(\bb R^n\backslash\{y\})^n$. Let $\Omega\subset\bb R^n$ be a bounded open set with smooth boundary, and define $\n O=\Omega\backslash B(y,r)$ for some $r>0$. Let $\mathbf{h}\in L^2(\n O)^n$. We can use a Hodge Decomposition as in Lemma 4 of Chapter 4 of \cite{AT}, so that there exists $g=g_{\mathbf{h}}\in\dot W^{1,2}(\n O)$ such that $\Vert\nabla g_{\mathbf{h}}\Vert_{L^2(\n O)}\leq\Vert \mathbf{h}\Vert_{L^2(\n O)}$. Moreover, $\mathbf{h}-\nabla g$ is a divergence-free vector in the weak sense; that is,
\[
\int_{\n O}(\mathbf{h}-\nabla g)\cdot\overline{\nabla\phi}=0,
\]
for each $\phi\in\dot W^{1,2}(\n O)$. We denote by $\n D'(\n O)$ the space of distributions on $\n O$. Let $\{f_k\}\subset C^{\infty}(\n O)$ satisfy
\[
f_k\longrightarrow \Gamma_L\qquad\text{in }D'(\n O)\text{ and in } L^2(\n O)\text{ as }k\ra\infty.
\]
We can also choose $\{f_k\}$ so that
\[
\Vert f_k\Vert_{L^{\infty}(\n O)}\leq C.
\]
For each $k\in\bb N$, we can write
\[
(\D f_k,\mathbf{h}-\nabla g)=(\nabla f_k,\mathbf{h}-\nabla g)-(i\mathbf{a} f_k,\mathbf{h}-\nabla g)=-(i\mathbf{a}f_k,\mathbf{h}-\nabla g),
\]
and since
\[
|(\mathbf{a} f_k,\mathbf{h}-\nabla g)|\leq C\Vert\mathbf{h}-\nabla g\Vert_{L^2(\n O)},
\]
it follows that for any $\mathbf{h}\in L^2(\n O)^n$,
\[
|(\D f_k,\mathbf{h}-\nabla g_{\mathbf{h}})|\leq C\Vert\mathbf{h}\Vert_{L^2(\n O)},
\]
where $C$ does not depend on $\mathbf{h}$. Consequently, in order to prove that $\D\Gamma_L$ is a bounded linear functional on $L^2(\n O)^n$ (hence, it lies in $L^2(\n O)^n$), it suffices to check that 
\begin{equation}
\Big|(\D\Gamma_L,\nabla g)\Big|\leq C\Vert\nabla g\Vert_{L^2(\n O)},\qquad\text{for each } g\in C_c^{\infty}(\n O).\label{distributionneeded}
\end{equation}
Since
\[
(\D\Gamma_L,\nabla g)=(\D\Gamma_L,\D g)+(\D\Gamma_L,i\mathbf{a}g),
\]
and
\[
\Big|(\D\Gamma_L,\D g)\Big|\leq\Big|(\D\cdot\D\Gamma_L,g)\Big|\leq C\Vert g\Vert_{L^1(\n O)}\leq C\Vert g\Vert_{L^{2^*}(\bb R^n)}\leq C\Vert\nabla g\Vert_{L^2(\n O)},
\]
\begin{multline}\nonumber
\Big|(\D\Gamma_L,i\mathbf{a}g)\Big|=\Big|(\Gamma_L,i\D\cdot(\mathbf{a}g))\Big|\leq C\Vert g\nabla\cdot{\mathbf a}+{\mathbf a}\cdot\nabla g-i|{\mathbf a}|^2g\Vert_{L^1(\n O)}\leq C\Vert\nabla g\Vert_{L^2(\n O)},
\end{multline}
then it's clear that (\ref{distributionneeded}) follows. Hence $\D\Gamma_L\in L^2(\n O)$. This ends the proof of the theorem.\hfill{$\square$}

\begin{remark} The above method of proof instantly generalizes the results of this theorem to the case where $\mathbf{a}\in L^4_{loc}(\bb R^n), \text{div }\mathbf{a}\in L^2_{loc}(\bb R^n)$, and $V\in L^{\infty}_{loc}(\bb R^n)$.
\end{remark}

\subsection{The fundamental solution of the magnetic Schr\"odinger operator, part II}
It is our intent to approximate the magnetic Schr\"odinger operators with rough coefficients by those with smooth coefficients. Take $\mathbf{a}\in L^2_{loc}(\bb R^n)$, and $V\in L^1_{loc}(\bb R^n)$ with $V\geq0$ a.e. on $\bb R^n$. We can construct sequences $\{\mathbf{a}_k\}\subset C_c^{\infty}(\bb R^n),\{V_k\}\subset C^{\infty}(\bb R^n)\cap L^{\infty}(\bb R^n)$ respectively, which converge to $\mathbf{a},V$ in the topology of $L^2_{loc}(\bb R^n), L^1_{loc}(\bb R^n)$, respectively. Moreover, the $V_k$'s can be chosen to be non-negative. We denote by $L_k$ the operator associated to $\mathbf{a}_k$ and $V_k$. For each $k\in\bb N$, Theorem \ref{testfundsol} applies, giving the existence of a fundamental solution $\Gamma_k$ to the operator $L_k$.

We will obtain our desired result in two main steps: First, we show weak convergence of $\dot L_k^{-1}$ to $\dot L^{-1}$ in $\dot{\m V}_{\mathbf{a},V}$. Second, we show convergence of the fundamental solutions $\Gamma_k$ to $\Gamma_L$ in the local topology associated to $\m V_{\mathbf{a},V}(\bb R^n\backslash\{y\})$.

\begin{lemma}\label{testconv} Assume that $\mathbf{a}\in L^2_{loc}(\bb R^n)$, $A\equiv I$, and $V\in L^1_{loc}(\bb R^n)$ with $V\geq0$ a.e. on $\bb R^n$. Fix $f\in L^{\infty}(\bb R^n)$ with compact support. Then
\begin{enumerate}[a)]
\item $\{D_{\mathbf{a}_k}\dot L_k^{-1}f\}$  converges to $\D\dot L^{-1}f$ strongly in $L^2(\bb R^n)$.
\item $\{V_k^{\frac12}\dot L^{-1}_kf\}$ converges to $V^{\frac12}\dot L^{-1}f$ strongly in $L^2(\bb R^n)$.
\item $\{\dot L^{-1}_kf\}$ converges weakly to $\dot L^{-1}f$ in $L^{2^*}(\bb R^n)$, and strongly in the $L^2_{loc}(\bb R^n)$ sense.
\end{enumerate}
\end{lemma}

\noindent\emph{Proof.} By Theorem \ref{firstone}, the sequence $u_k=\dot L_k^{-1}f$ is well-defined, and for each $k\in\bb N$, $u_k\in\dot W^{1,2}(\bb R^n)$. By construction, we have
\begin{equation}\label{uksol}
\int_{\bb R^n}D_{\mathbf{a}_k}u_k\overline{D_{\mathbf{a}_k}\phi}+V_ku_k\overline{\phi}=(f,\phi),
\end{equation}
for each $\phi\in\dot W^{1,2}(\bb R^n)$ and each $k\in\bb N$. Plugging in $\phi=u_k$, we see that
\begin{equation}\nonumber
\Vert u_k\Vert_{\dot{\m V}_{\mathbf{a}_k,V_k}}^2=\int_{\bb R^n}|D_{\mathbf{a}_k}u_k|^2+V_k|u_k|^2=(f,u_k)\leq C\Vert f\Vert_{L^2(\supp f)}\Vert u_k\Vert_{L^2(\supp f)},
\end{equation}
and since by the Sobolev inequality and the diamagnetic inequality we have
\[
\Vert u_k\Vert_{L^2(\supp f)}\leq C\Vert u_k\Vert_{L^{2^*}(\bb R^n)}\leq C\Vert D_{\mathbf{a_k}}u_k\Vert_{L^2(\bb R^n)}\leq C\Vert u_k\Vert_{\dot{\m V}_{\mathbf{a}_k,V_k}},
\]
where $C$ depends on $\supp f$ but not on $k$, then it follows that
\[
\Vert u_k\Vert_{\dot{\m V}_{\mathbf{a}_k,V_k}}\leq C
\]
where $C$ depends on $f$. It follows that the sequences $\{D_{\mathbf{a_k}} u_k\},\{V_k^{\frac12}u_k\}$ are uniformly bounded in $L^2(\bb R^n)$. Applying the diamagnetic inequality for each $k\in\bb N$, it follows that $\{u_k\}$ is a uniformly bounded sequence in $L^{2^*}(\bb R^n)$, so that in particular $\{u_k\}$ is uniformly bounded in the $L^2_{loc}(\bb R^n)$ sense (by this we mean that the sequence is uniformly bounded in $L^2(\Omega)$ for each bounded set $\Omega\subset\bb R^n$). Actually, using the Moser estimate, Theorem \ref{moser}, for each $k\in\bb N$, we can see that $\{u_k\}$ is a uniformly bounded sequence in the $L^{\infty}_{loc}(\bb R^n)$ sense.  Let $g,h,u$ be weak limits in $L^2(\bb R^n),L^2(\bb R^n)$ and $L^{2^*}(\bb R^n)$ of $\{D_{\mathbf{a_k}} u_k\},\{V_k^{\frac12}u_k\}, \{u_k\}$, respectively. Now fix $\psi\in C_c^{\infty}(\bb R^n)$. Diagonalizing, we pass to an indexing set where all three sequences achieve the aforementioned weak limits, and for ease of notation we say that the entire sequences converge. Observe that for each $k\in\bb N$ we have
\begin{gather}\label{dualizing1}
(\D u,\psi)=\int_{\bb R^n}u\overline{\D\psi}=\int_{\bb R^n}u_k\overline{D_{\mathbf{a}_k}\psi}+\int_{\bb R^n}u_k\overline{[\D\psi-D_{\mathbf{a}_k}\psi]}+\int_{\bb R^n}[u-u_k]\overline{\D\psi},
\end{gather}
and since 
\[
\Bigl|\int_{\bb R^n}[u-u_k]\overline{\D\psi}\Bigr|\leq C\Vert u-u_k\Vert_{L^2(\supp\psi)}\longrightarrow0,
\]
\[
\Bigl|\int_{\bb R^n}u_k\overline{[\D\psi-D_{\mathbf{a}_k}\psi]}\Bigr|\leq C\Vert u_k\Vert_{L^2(\supp\psi)}\Vert\mathbf{a}-\mathbf{a}_k\Vert_{L^2(\supp\psi)}\longrightarrow0,
\]
as $k\ra\infty$, then from (\ref{dualizing1}) we have
\[
(\D u,\psi)=\lim\limits_{k\ra\infty}\int_{\bb R^n}u_k\overline{D_{\mathbf{a}_k}\psi}=\lim\limits_{k\ra\infty}(D_{\mathbf{a}_k}u_k,\psi)=(g,\psi),
\]
whence by varying over $\psi\in C_c^{\infty}(\bb R^n)$ we conclude that $\D u\in L^2(\bb R^n)$ and $\D u\equiv g$ in $L^2(\bb R^n)$. Likewise, since $\{V_k^{\frac12}\}$ is a uniformly bounded sequence in the $L^2_{loc}(\bb R^n)$ sense, and
\[
\Bigl|\int_{\bb R^n}V^{\frac12}[u-u_k]\psi\Bigr|\leq\Vert\psi\Vert_{L^{\infty}(\bb R^n)}\Vert V\Vert_{L^1(\supp\psi)}\Vert u-u_k\Vert_{L^2(\supp\psi)}\longrightarrow0,
\]
\begin{align*}
\Bigl|\int_{\bb R^n}[V^{\frac12}-V_k^{\frac12}]u_k\psi\Bigr|&\leq\Vert\psi\Vert_{L^{\infty}(\bb R^n)}\Vert u_k\Vert_{L^2(\supp\psi)}\Bigl(\int_{\supp\psi}|V^{\frac12}-V_k^{\frac12}|^2\Bigr)^{\frac12}\\[4mm]&\leq C\Vert V-V_k\Vert_{L^1(\supp\psi)}^{\frac12}\longrightarrow0,
\end{align*}
as $k\ra\infty$, a similar argument to the one above concludes that $V^{\frac12}u\in L^2(\bb R^n)$ and $V^{\frac12}u\equiv h$ in $L^2(\bb R^n)$. Hence it follows that $u\in\dot{\m V}_{\mathbf{a},V}$. Taking limit in (\ref{uksol}) as $k\ra\infty$ for fixed $\phi\in C_c^{\infty}(\bb R^n)$, we get

\[
\int_{\bb R^n}\D u\overline{\D\phi}+Vu\overline{\phi}=(f,\phi),
\]
which is valid by the aforementioned observations and the fact that $\phi\in C_c^{\infty}(\bb R^n)$. By the uniqueness of such a function $u$, we have proven that $u=\dot L^{-1}f$. Since we defined $u$ to be a weak limit of $\{u_k\}$ in $L^{2^*}(\bb R^n)$, then the first part of c) is proven.

To prove the last part of c), fix $\Omega\subset\bb R^n$ a bounded open set. By the Moser estimate, Theorem \ref{moser}, it follows that $\{u_k\}$ is a uniformly bounded sequence in $L^{\infty}(\Omega)$. Since $\{a_k\}$ is a uniformly bounded sequence in $L^2(\Omega)$ and $\{D_{{\mathbf a}_k}u_k\}$ is uniformly bounded in $L^2(\Omega)$, it therefore follows that $\{\nabla u_k\}$ is a uniformly bounded sequence in $L^2(\Omega)$. By the Rellich-Kondrakov Theorem, a subsequence of $\{u_k\}$ must be strongly convergent in $L^2(\Omega)$, and this limit has no choice but to be $u$, hence the limit is unique and the strong convergence occurs along the whole sequence.

Finally, to prove the strong convergences in a) and b), first observe that
\[
\Bigl|\int_{\bb R^n}f(u_k-u)\Bigr|\leq\Vert f\Vert_{L^2(\supp f)}\Vert u_k-u\Vert_{L^2(\supp f)}\longrightarrow0\quad\text{as }k\ra\infty,
\]
and so by an argument analogous to (\ref{strongconvarg}), we deduce that
\begin{equation}\label{normconv}
\Vert D_{\mathbf{a}_k}u_k\Vert_{L^2(\bb R^n)}^2+\Vert V_k^{\frac12}u_k\Vert_{L^2(\bb R^n)}^2\longrightarrow \Vert\D u\Vert_{L^2(\bb R^n)}^2+\Vert V^{\frac12}u\Vert_{L^2(\bb R^n)}^2
\end{equation}
as $k\ra\infty$. Since the sequence of numbers on the left-hand side of (\ref{normconv}) is non-negative, it follows that both terms converge independently. Suppose that there exists $\ep>0$ such that
\[
\Vert D_{\mathbf{a}_k}u_k\Vert_{L^2(\bb R^n)}^2\longrightarrow\Vert\D u\Vert_{L^2(\bb R^n)}^2+\ep.
\]
Then, owing to (\ref{normconv}), it must be the case that
\[
\Vert V_k^{\frac12}u_k\Vert_{L^2(\bb R^n)}^2\longrightarrow \Vert V^{\frac12}u\Vert_{L^2(\bb R^n)}^2-\ep,
\]
but this is absurd, since by the weak convergence of the sequence $\{V_k^{\frac12}u_k\}$, we have
\[
\Vert V^{\frac12}u\Vert_{L^2(\bb R^n)}^2\leq\liminf\limits_{k\ra\infty}\Vert V_k^{\frac12}u_k\Vert_{L^2(\bb R^n)}^2.
\]
Therefore, there cannot exist such $\ep>0$. Since $\ep<0$ is also absurd by the weak convergence of the sequence $\{D_{\mathbf{a}_k}u_k\}$, it follows that
\[
\Vert D_{\mathbf{a}_k}u_k\Vert_{L^2(\bb R^n)}^2\longrightarrow\Vert\D u\Vert_{L^2(\bb R^n)}^2,\qquad\Vert V_k^{\frac12}u_k\Vert_{L^2(\bb R^n)}^2\longrightarrow \Vert V^{\frac12}u\Vert_{L^2(\bb R^n)}^2
\]
which, together with the weak convergences already shown, imply the respective strong convergences.\hfill{$\square$}


\begin{theorem}\label{fundsolapprox}Assume that $\mathbf{a}\in L^2_{loc}(\bb R^n)$, $A\equiv I$, and $V\in L^1_{loc}(\bb R^n)$ with $V\geq0$ a.e. on $\bb R^n$. Then for each $y\in\bb R^n$, the fundamental solutions $\Gamma_k$ to the operators $L_k$ converge in the weak topology of $\m V_{\mathbf{a},V}(\bb R^n\backslash\{y\})$, and locally in the strong $L^2_{loc}(\bb R^n\times\bb R^n\backslash\{x=y\})$ sense, to $\Gamma_L$, the measurable function of Theorem \ref{firstone}. In particular, $\Gamma_L(\cdot,y)\in\m V_{\mathbf{a},V,loc}(\bb R^n\backslash\{y\})$, and $L\Gamma_L=0$ in the weak sense on $\bb R^n\backslash\{y\}$, so that $\Gamma_L$ is the fundamental solution of the operator $\dot L$. Moreover,
\begin{equation}\label{adjointcompare1}
\Gamma_L(x,y)\equiv\overline{\Gamma_{L^*}(y,x)}
\end{equation}
is true in the $a.e.$ sense on $\bb R^n\times\bb R^n$.
\end{theorem}

\noindent\emph{Proof.} Let $U_1,U_2$ be arbitrary open bounded subsets of $\bb R^n$ such that
\begin{equation}\label{distantsets}
\dist(U_1,U_2)=3r,\quad r>0,
\end{equation}
and let $\phi\in C_c^{\infty}(U_1)$, $f\in C_c^{\infty}(U_2)$. On the set $\Omega:=U_1\times U_2$, $\{\Gamma_k\}$ is a uniformly bounded sequence in $L^p(\Omega)$ for $p\in[1,\infty]$ owing to property b) in Theorem \ref{firstone} (with a norm depending on $U_1,U_2,r$). By Theorem \ref{testfundsol}, for each $y\in U_2$, $\Gamma_k(\cdot,y)$ solves $L\Gamma_k(\cdot,y)=0$ on $\bb R^n\backslash\{y\}$ in the weak sense, and the analogous statement holds for the $y-$variable because of (\ref{adjointcompare}). Cover $U_1$ by a family of balls $\{B_r^m\}$ such that the balls intersect a uniformly finite number of times depending only on dimension (recall that $r$ is given by (\ref{distantsets})). Applying the Caccioppoli inequality (\ref{cacc}) with $R=2r$ on each ball $B_r^m$, it is straightforward that
\[
\Vert D_{\mathbf{a}_k}\Gamma_k(\cdot,y)\Vert_{L^2(B_r^m)}^2\leq C(r)\int_{B_{2r}^m}|\Gamma_k(\cdot,y)|^2.
\]
Consequently,
\begin{align*}
\Vert D_{\mathbf{a}_k}\Gamma_k(\cdot,y)\Vert_{L^2(U_1)}^2&\leq\sum\limits_{m=1}^{\infty}\Vert D_{\mathbf{a}_k}\Gamma_k(\cdot,y)\Vert_{L^2(B_r^m)}^2\leq\sum\limits_{m=1}^{\infty}C(r)\int_{B_{2r}^m}|\Gamma_k(\cdot,y)|^2\\[4mm]&\leq C(r,n)\int_{U_1+r}|\Gamma_k(\cdot,y)|^2\leq C(r,n,U_1),
\end{align*}
whence we see that for each $y\in U_2$, the sequence $\{D_{\mathbf{a}_k}\Gamma_k(\cdot,y)\}$ is uniformly bounded in $L^2(U_1)$. Similarly, thanks to (\ref{adjointcompare}), it is proven that for each $x\in U_1$, $\{D_{\mathbf{a}_k}\Gamma_k(x,\cdot)\}$ is uniformly bounded in $L^2(U_2)$. Combining these results, we obtain that $\{D_{\mathbf{a}_k}\Gamma_k\}$ is uniformly bounded in $L^2(\Omega)$. Therefore, since $\{\Gamma_k\}$ is a uniformly bounded sequence in $L^{\infty}(\Omega)$ (due to (\ref{pointwiseupper})), then from the fact that
\[
D_{\mathbf{a}_k}\Gamma_k=\nabla\Gamma_k-i{\mathbf a}_k\Gamma_k,
\]
it actually follows that $\{\nabla\Gamma_k\}$ is a uniformly bounded sequence in $L^2(\Omega)$. By the Rellich-Kondrachov Theorem, it follows that a subsequence of $\{\Gamma_k\}$ converges strongly in $L^2(\Omega)$. Hence we pass to such an $L^2(\Omega)-$convergent subsequence (which we denote as the whole sequence).

On the other hand, from Lemma \ref{testconv}, it follows that $\{\dot L_k^{-1}f\}$ is a uniformly bounded sequence in $L^{2^*}(\bb R^n)$, so that in particular we can write
\[
\int_{\bb R^n}\Big[\dot L_k^{-1}f-\dot L^{-1}f\Big]\phi(x)\,dx\longrightarrow0
\]
as $k\ra\infty$. Using the kernel representation, we have that
\begin{equation}\label{thisgoesto0}
\int_{\bb R^n}\int_{\bb R^n}\Big[\Gamma_k(x,y)-\Gamma_L(x,y)\Big]f(y)\,dy\phi(x)\,dx\longrightarrow0
\end{equation}
as $k\ra\infty$. Therefore, $\Gamma_L(x,y)$ must be the unique strong limit in $L^2(\Omega)$ of the whole sequence $\{\Gamma_k\}$. After passing to a subsequence, we can assume that $\{\Gamma_k\}$ converges pointwise a.e. on $\Omega$ to $\Gamma_L$. In particular, for a.e. $y\in U_2$, there is a subsequence of $\{\Gamma_k(\cdot,y)\}$, with the indexing set independent of $y$ (but depending on $\Omega=U_1\times U_2$), which converges pointwise a.e. on $U_1$ to $\{\Gamma_L(\cdot,y)\}$.

By the aforementioned discussion and similar argumentation to that in Lemma \ref{testconv}, it must be the case that $\{D_{\mathbf{a}_k}\Gamma_k(\cdot,y)\}$ converges weakly in $L^2(U_1)$ to $\D\Gamma_L(\cdot,y)$, and moreover by the weak$^*$ convergence in $L^{\infty}(U_1)$ of a subsequence of $\{\Gamma_k(\cdot,y)\}$, it therefore follows that $\{V_k\Gamma_k(\cdot,y)\}$ converges to $\{V\Gamma_L(\cdot,y)\}$ on $U_1$ in the sense of distributions. Then, since $\Gamma_L(\cdot,y)\in L^{\infty}(U_1)$, we have shown that $\Gamma_L(\cdot,y)\in\m V_{\mathbf{a},V}(U_1)$. For fixed $\phi\in C_c^{\infty}(U_1)$, taking the limit as $k\ra\infty$ on each identity (\ref{weak}) satisfied by $\Gamma_k(\cdot,y)$, we arrive at the fact that $\Gamma_L$ solves $L\Gamma_L(\cdot,y)=0$ on $U_1$ in the weak sense. By (\ref{adjointcompare}), all results in the $x-$variable are also true in the $y-$variable, and the pointwise a.e. convergence property of $\{\Gamma_k\}$ to $\Gamma_L$, coupled with (\ref{adjointcompare}), gives (\ref{adjointcompare1}) immediately on $U_1\times U_2$. Varying over admissible $U_1,U_2$ finishes the proof of the theorem.\hfill{$\square$}\\

}

\section{Upper bound on the exponential decay of the fundamental solution}\label{upperboundsec}

In the following definition, $\n L$ is either $L+\ep$ for $\ep>0$ or $\dot L$, with $\n H,\n R$ the domain and range respectively of these operators. Specifically, when $\n L=L+\ep$, we write $\n H=D(L)$, $\n R=L^2(\bb R^n)$. When $\n L=\dot L$, we write $\n H=\dot{\m V}$, $\n R=\dot{\m V}'$.

\begin{definition} We say that the operator $\n L$ has the \emph{zero-source local uniform boundedness property} if for each ball $B\subset\bb R^n$ and each function $u$ which solves $\n Lu=0$ in the weak sense on $2B$, it follows that $u\in L^{\infty}(B)$ and
\begin{equation}\label{moser1}
\Vert u\Vert_{L^{\infty}(cB)}\leq C\Bigl(~\fint_{2B}|u|^2\Bigr)^{\frac12},
\end{equation}
where $c, C$ are independent of $B$ and $u$.
\end{definition}

We briefly remark that for the operators under consideration, the diamagnetic inequality guarantees that $u\in L^2_{loc}(\bb R^n)$.

For the fundamental solution $\Gamma$ associated to the operator $\n L$, as given in Definition \ref{fundsoln}, it will be useful to consider the conditions
\begin{equation}\label{pointwiseupper1}
	|\Gamma(x,y)|\leq\frac{C}{|x-y|^{n-2}},\qquad\text{for a.e. }x,y\in\bb R^n,\, x\neq y
\end{equation}
and
\begin{equation}\label{fundamentallower2}
	\frac{C_1}{|x-y|^{n-2}}\leq|\Gamma(x,y)|\leq\frac{C_2}{|x-y|^{n-2}},\qquad\text{for a.e. }x,y\in\bb R^n,\,x\neq y,
\end{equation}
where the constants $C,$ $C_1$, $C_2$ depend on the dimension only. In some situations (strictly speaking, not necessarily including ours), the condition (\ref{pointwiseupper1})  is equivalent to a Moser-type bound \cite{KK}, but we do not explore this direction.

Let us first state a result analogous to Proposition \ref{proposition1} for the operator $\dot L$. The proof is virtually identical to that of Proposition \ref{proposition1}, and is thus omitted.

\begin{proposition}\label{proposition2} Assume that $\mathbf{a}\in L^2_{loc}(\bb R^n)$, $A$ is an elliptic matrix with complex, bounded, measurable coefficients, and $V\in L^1_{loc}(\bb R^n)$ is real-valued with $V\geq0$ a.e. on $\bb R^n$. If $\mathbf{a}\equiv0$, assume $V\in RH_{\frac n2}$, otherwise take assumptions (\ref{magassumptions-intro}). Suppose $ U\subset\bb R^n$ is a bounded open set. Let $u\in\m V_{loc}(\bb R^n\backslash U)$ be a solution to $Lu=0$ in the weak sense on $\bb R^n\backslash U$. Suppose $\phi\in C_c^{\infty}(\bb R^n)$ is such that $\phi\equiv0$ on $2 U$. Let $g\in C^{\infty}(\bb R^n)$ be a non-negative function satisfying $|\nabla g(x)|\leq C_2m(x,V+|\mathbf{B}|)$ for every $x\in\bb R^n$. Then
	\begin{equation}\label{estimate2}
		\int_{\bb R^n}m(x,V+|\mathbf{B}|)^2|u\phi|^2e^{2\ep g}\,dx\leq C\int_{\bb R^n}|u|^2|\nabla\phi|^2e^{2\ep g}\,dx,
	\end{equation}
	for any $\ep\in(0,\ep_0)$, where $\ep_0$ and $C$ depend only on $\lambda,\Lambda, C_2, n, \Vert V+|\mathbf{B}|\Vert_{RH_{\frac n2}},$ and the constants from (\ref{magassumptions-intro}).
\end{proposition}

%
%

The next results follow the lines of the argument in \cite{S1}:

\begin{theorem}\label{Upper Bound} Suppose $\mathbf{a}\in L^2_{loc}(\bb R^n)$, $A$ is an elliptic matrix with complex, bounded coefficients, $V\in L^1_{loc}(\bb R^n)$, and that $L$ is an operator for which there exists a fundamental solution in the sense of Definition \ref{fundsoln} satisfying (\ref{pointwiseupper1}). Moreover, if $\mathbf{a}\equiv0$, assume $V\in RH_{\frac n2}$; otherwise assume (\ref{magassumptions-intro}). Then there exists $\ep>0$ and a constant $C>0$, depending on $L$, such that
\begin{equation}\label{l2upperbound}
\Bigl(\fint_{B(x,\frac1{m(x,V+|\mathbf{B}|)})}|\Gamma(z,y)|^2\,dz\Bigr)^{\frac12}\leq\frac{Ce^{-\ep d(x,y,V+|\mathbf{B}|)}}{|x-y|^{n-2}} \mbox{ for all }x, y \in \rn.
\end{equation} 	
If $L$ satisfies the zero-source local uniform boundedness property, then	
\begin{equation}\label{upperbound}
|\Gamma(x,y)|\leq\frac{Ce^{-\ep d(x,y,V+|\mathbf{B}|)}}{|x-y|^{n-2}} \mbox{ for all }x, y \in \rn.
\end{equation}
\end{theorem}

\bp  Fix $x_0,y_0\in\bb R^n$, $x_0\neq y_0$. Since for each $x,y\in\bb R^n$ we have
\[
d(x,y)\leq C\quad\text{if}\quad|x-y|<\frac{C}{m(x,V+|{\mathbf B}|)},
\]
then in particular we may also assume that $|x_0-y_0|\geq\frac C{m(y_0,V)}$. Indeed, otherwise $d(x_0,y_0)\leq C$, and so from (\ref{pointwiseupper1}), we observe
\[
|\Gamma(x_0,y_0)|\leq\frac{C}{|x-y|^{n-2}}\leq\frac{Ce^{\ep C}e^{-\ep d(x_0,y_0)}}{|x-y|^{n-2}}
\]
which gives (\ref{upperbound}). Furthermore, we can assume that
\begin{equation}\label{nointersect}
B\Bigl(x_0,\frac4{m(x_0,V+|{\mathbf B}|)}\Bigr)\cap B\Bigl(y_0,\frac4{m(y_0,V+|{\mathbf B}|)}\Bigr)=\varnothing.
\end{equation}
Indeed, suppose there exists $z\in B\Bigl(x_0,\frac4{m(x_0,V+|{\mathbf B}|)}\Bigr) \cap B\Bigl(y_0,\frac4{m(y_0,V+|{\mathbf B}|)}\Bigr)$. Then we note
\begin{multline}
|x_0-y_0|\leq|x_0-z|+|z-y_0|\leq\frac4{m(x_0,V+|{\mathbf B}|)}+\frac4{m(y_0,V+|{\mathbf B}|)}\nonumber\\[2mm]\leq8\max\Bigl\{\frac1{m(x_0,V+|{\mathbf B}|)}~,~\frac1{m(y_0,V+|{\mathbf B}|)}\Bigr\}
\end{multline}
which once again implies $d(x_0,y_0)\leq C$.

Per our assumptions, $L\Gamma(\cdot,y_0)=0$ in the weak sense on any ball centered around $x_0$ which does not contain $y_0$. Let $r=\frac1{m(y_0,V+|\mathbf{B}|)}$. Applying Proposition \ref{proposition2} with $u=\Gamma(\cdot,y_0)$ and $U=B(y_0,4r)$, we can carry out the argument of Theorem \ref{l2decay} up to proving (\ref{greatineq}) to establish that
\begin{multline}\label{greatineq1}
\int_{B_M\backslash4 U}m(\cdot,V+|\mathbf{B}|)^2|u|^2e^{2\ep \varphi_j}\\[3mm]\leq C\Bigl\{\int_{4 U\backslash2 U}|u|^2\frac1{|\diam U|^2}e^{2\ep \varphi_j}+\int_{2B_M\backslash B_M}|u|^2\frac1{M^2}e^{2\ep \varphi_j}\Bigr\}.
\end{multline}
where $\varphi_j$ is as in Proposition \ref{dapprox1} with $w=V+|{\mathbf B}|$. Owing to (\ref{pointwiseupper1}) and the fact that $\phi_j$ is uniformly bounded in $2B_M\backslash B_M$, it follows that the second term in the right-hand side of (\ref{greatineq1}) drops to $0$ as $M\ra\infty$. Therefore, we may conclude that
\begin{align}
\int_{\bb R^n\backslash B(y_0,4r)}m(x,V+|\mathbf{B}|)^2|\Gamma(x,y_0)|^2&e^{2\ep\varphi_j(x,y_0)}\,dx\nonumber\\[4mm]&\leq\frac C{r^2}\int_{B(y_0,4r)\backslash B(y_0,2r)}|\Gamma(x,y_0)|^2e^{2\ep\varphi_j(x,y_0)}\,dx.\label{greatfund}
\end{align}
We note that for $x\in B(y_0,4r)\backslash B(y_0,2r)$, by Proposition \ref{dbounded} we have $d(x,y_0,V+|\mathbf{B}|)\leq C$, which by (\ref{closetod}) implies $\varphi(x,y_0)\leq C$, where $\varphi$ is as in Proposition \ref{dapprox}. But of course, by the construction of the $\varphi_{j}$, we see that $\varphi_{j}(x,y_0)\leq\varphi(x,y_0)\leq C$ for all $j$ and $x\in B(y_0,4r)\backslash B(y_0,2r)$. It follows that
\[
\sup\limits_{x\in B(y_0,4r)\backslash B(y_0,2r)}e^{2\ep\varphi_{j}(x,y_0)}\leq C,\qquad\text{for each } j\in\bb N,
\]
so that, using (\ref{pointwiseupper1}),
\[
\int_{\bb R^n\backslash B(y_0,4r)}m(x,V+|\mathbf{B}|)^2|\Gamma(x,y_0)|^2e^{2\ep\varphi_j(x,y_0)}\,dx\leq\frac{C}{r^{n-2}}.
\]
By Fatou's Lemma we have
\begin{equation}\label{star1}
\int_{\bb R^n\backslash B(y_0,4r)}m(x,V+|\mathbf{B}|)^2|\Gamma(x,y_0)|^2 e^{2\ep\varphi(x,y_0)}\,dx\leq\frac{C }{r^{n-2}}.
\end{equation}
Let $R=\frac1{m(x_0,V+|\mathbf{B}|)}$. We claim that $B(x_0,R)\subset \bb R^n\backslash B(y_0,4r)$. Suppose not. Then there is a point $z\in B(x_0,R)\cap\Big(\bb R^n\backslash B(y_0,4r)\Big)^c= B(x_0,R)\cap B(y_0,4r)$, which contradicts (\ref{nointersect}).  Now, using (\ref{closetod}) and (\ref{star1}), we get
\[
\int_{B(x_0,R)}m(x,V+|\mathbf{B}|)^2|\Gamma(x,y_0)|^2 e^{2\ep d(x,y_0)}\,dx\leq\frac{C}{r^{n-2}}.
\]
By the Triangle Inequality we observe
\[
d(x_0,y_0)\leq d(x_0,x)+d(x,y_0),\quad\text{for each } x\in B(x_0,R).
\]
Recall that we have by Proposition \ref{propdbounded} that $d(x_0,x)\leq L$, $L$ a constant depending on $\Vert V+|\mathbf{B}|\Vert_{RH_{\frac n2}}$ but independent of $x_0$ and $x$. It follows that
\[
e^{2\ep d(x,y_0)}\geq e^{-2\ep L}e^{2\ep d(x_0,y_0)}=Ce^{2\ep d(x_0,y_0)},\quad\text{for each } x\in B(x_0,R).
\]
From this fact and the fact that $m(x,V+|\mathbf{B}|)\sim m(x_0,V+|\mathbf{B}|)=R^{-1}$ for every $x\in B(x_0,R)$ owing to (\ref{slow}), we can conclude
\[
\frac1{R^2}\int_{B(x_0,R)}|\Gamma(x,y_0)|^2 e^{2\ep d(x_0,y_0)}\,dx\leq\frac{C}{r^{n-2}}.
\]
Dividing out by $R^{n-2}$ and taking square root we see that
\begin{align}
\Bigl(\frac1{R^n}\int_{B(x_0,R)}|\Gamma(x,y_0)|^2\,dx\Bigr)^{1/2}&\leq\quad\frac{Ce^{-\ep d(x_0,y_0)}}{(Rr)^{(n-2)/2}}\nonumber\\[4mm]&\leq C\Big[m(x_0,V+|\mathbf{B}|)m(y_0,V+|\mathbf{B}|)\Big]^{(n-2)/2}e^{-\ep d(x_0,y_0)}.\label{l2almostthere}
\end{align}

We claim that if $|x-y|m(x,V+|\mathbf{B}|)\geq2$ then
\begin{equation}\label{claim}
d(x,y)\geq c\Big[1+|x-y|m(x,V+|\mathbf{B}|)\Big]^{1/(k_0+1)},
\end{equation}
for some $k_0>0$. To see this, choose $\gamma:[0,1]\ra\bb R^n$ such that $\gamma(0)=x$, $\gamma(1)=y$, and
\[
2d(x,y)\geq\int_0^1m(\gamma(t),V+|\mathbf{B}|)|\gamma'(t)|\,dt
\]
which can be done by the definition of $d$ and the fact that $x\neq y$ necessarily in this case. It follows from (\ref{polygrowth2}) that
\begin{align*}
2d(x,y)\geq c\int_0^1\frac{m(x,V+|\mathbf{B}|)|\gamma'(t)|\,dt}{\big[1+|\gamma(t)-x|m(x,V+|\mathbf{B}|)\big]^{\frac{k_0}{k_0+1}}}.
\end{align*}
The integral on the right-hand side of the above inequality is greater than or equal to the geodesic distance from $x$ to $y$ in the metric
\[ \frac{m(x,V+|\mathbf{B}|)dz}{\big[1+|z-x|m(x,V+|\mathbf{B}|)\big]^{\frac{k_0}{k_0+1}}},\]
which is
\[
\int_0^1\frac{m(x,V+|\mathbf{B}|)|y-x|\,dt}{\big[1+t|y-x|m(x,V+|\mathbf{B}|)\big]^{\frac{k_0}{k_0+1}}}\geq c'\Big[1+|x-y|m(x,V+|\mathbf{B}|)\Big]^{1/(k_0+1)}
\]
and so (\ref{claim}) follows.

Hence, owing to (\ref{nointersect}) and (\ref{claim}), we have for each $\ep'>0$,
\begin{align*}
|x_0-y_0|m(x_0,V+|\mathbf{B}|)\leq1+|x_0-y_0|m(x_0,V+|\mathbf{B}|)&\leq\frac1{c^{k_0+1}}d(x_0,y_0)^{k_0+1}\nonumber\\[4mm]&\leq\frac1{c^{k_0+1}}C_{\ep'/2}e^{(\ep'/2)d(x_0,y_0)}.
\end{align*}
Likewise,
\[
|x_0-y_0|m(y_0,V+|\mathbf{B}|)\leq\frac1{c^{k_0+1}}C_{\ep'/2}e^{(\ep'/2)d(x_0,y_0)}.
\]
Adding the last two inequalities we see
\[
|x_0-y_0|m(y_0,V+|\mathbf{B}|)+|x_0-y_0|m(x_0,V+|\mathbf{B}|)\leq\frac2{c^{k_0+1}}C_{\ep'/2}e^{(\ep'/2)d(x_0,y_0)},
\]
for any $\ep'>0$. Squaring the above inequality we have
\[
2|x_0-y_0|^2m(x_0,V+|\mathbf{B}|)m(y_0,V+|\mathbf{B}|)\leq4c^{-2(k_0+1)}C_{\ep'/2}^2e^{\ep' d(x_0,y_0)}.
\]
Now, given $\ep$ small enough, choose $\ep'=\frac1{n-2}\ep$. Then in view of (\ref{l2almostthere}), we get
\[
\Bigl(\frac1{R^n}\int_{B(x_0,R)}|\Gamma(x,y_0)|^2\,dx\Bigr)^{1/2}\leq\frac{Ce^{-(\ep/2)d(x_0,y_0)}}{|x_0-y_0|^{n-2}},
\]
which is (\ref{l2upperbound}). Since $\Gamma(\cdot,y_0)$ is a weak solution to $Lu=0$ on $B(x_0,R)$, then if the operator $L$ has the zero-source local uniform boundedness property, by (\ref{moser1}) and (\ref{l2upperbound}), we immediately achieve (\ref{upperbound}).\hfill{$\square$}\\

The exponential decay result of the last theorem holds in the pointwise a.e. sense for the fundamental solutions of some of the operators previously considered. In particular, we have
\begin{corollary}\label{upperboundcor} Let $L_1$ be a generalized magnetic Schr\"odinger operator formally given by (\ref{generalized}) where $\mathbf{a}\in L^2_{loc}(\bb R^n)$, $A\equiv I$, $V\in L^1_{loc}(\bb R^n)$ with $V\geq0$ a.e. on $\bb R^n$, and assumptions (\ref{magassumptions-intro}) are satisfied. Then there exists $\ep>0$ and a constant $C>0$, depending on $\Vert V+|{\mathbf B}|\Vert_{RH_{\frac n2}}, n,$ and the constants from (\ref{magassumptions-intro}), such that
\begin{equation}\label{upperbound1}
|\Gamma_{L_1}(x,y)|\leq\frac{Ce^{-\ep d(x,y,V+|\mathbf{B}|)}}{|x-y|^{n-2}} \mbox{ for a.e. }x, y \in \rn.
\end{equation}
Let $L_2$ be a generalized magnetic Schr\"odinger operator formally given by (\ref{generalized}) where $\mathbf{a}\equiv0$, $A$ is an elliptic matrix with real, bounded coefficients, and $V\in RH_{\frac n2}$. Then there exists $\ep>0$ and a constant $C>0$, depending on $\Vert V\Vert_{RH_{\frac n2}}, n,\lambda,$ and $\Lambda$ such that
\begin{equation}\label{upperbound2}
\Gamma_{L_2}(x,y)\leq\frac{Ce^{-\ep d(x,y,V)}}{|x-y|^{n-2}} \mbox{ for all }x, y \in \rn\text{ with }x\neq y.
\end{equation}
\end{corollary}
\noindent\emph{Proof.} In the first setting, the results of Section \ref{fundsolnsec} are true, and so the hypothesis of Theorem \ref{Upper Bound} hold. In the second setting, the theory of fundamental solutions set forth in \cite{DHM} applies, whence the hypothesis of Theorem \ref{Upper Bound} hold. Furthermore, in this case the fundamental solution is actually positive and continuous (see Section \ref{electric}), so (\ref{upperbound}) holds pointwise on $\bb R^n\times\bb R^n\backslash\{x=y\}$.

\section{Lower bound on the exponential decay of the fundamental solution}\label{lowerboundsec}

\subsection{Properties of the generalized Schr\"odinger operator, without the magnetic potential.}\label{electric}
Recall the definition of the operator $L_E$,
\[
L_E=-\text{div }A\nabla~+~V,
\]
which is the operator $L$ with $\mathbf{a}\equiv0$. For this operator with $A$ an elliptic matrix with real, bounded coefficients, and $V\in RH_{\frac n2}$, the theory set forth in \cite{DHM} applies, so that the fundamental solution in the sense of Definition \ref{fundsoln} exists. Actually, its fundamental solution is known to be continuous and positive. Below we present a few lemmas that apply to this operator; the most critical for us in obtaining the lower bound estimate is the scale-invariant Harnack Inequality.


\begin{lemma}\label{Moser}\cite{DHM} (Moser-type Estimate) Assume that $\mathbf{a}\equiv0$, $A$ is an elliptic matrix with real, bounded coefficients, and $V\in L^1_{loc}(\bb R^n)$ with $V\geq0$ a.e. on $\bb R^n$. Let $B_R\subset\bb R^n$ be a ball, and let $u\in\m V_{0,V}(B_R)$ solve $L_Eu=f$ in the weak sense on $B_R$, where $f\in L^q(B_R)$ for some $q>\frac n2$. Then for any $r,\,0<r<R$,
\[
\Vert u\Vert_{L^{\infty}(B_{R/2})}\leq C\Bigl[\Bigl(\fint_{B_R}|u|^r\Bigr)^{1/r}+R^{2-\frac nq}\Vert f\Vert_{L^q(B_R)}\Bigr]
\]
where $C$ depends on $n,p,q,\lambda,\Lambda$.
\end{lemma}

\begin{lemma}\label{Holder cont}\cite{DHM} (H\"older Continuity Estimate) Assume that $\mathbf{a}\equiv0$, $A$ is an elliptic matrix with real, bounded coefficients, and $V\in RH_{\frac n2}$. Let $u$ solve $L_Eu=0$ in the weak sense on a ball $B_{R_0}\subset\bb R^n, R_0>0$.  Then there exists $\eta\in(0,1)$ depending on $R_0$, and $C_{R_0}>0$ such that whenever $0<R\leq R_0$,
\[
\sup\limits_{x,y\in B_{R/2},~x\neq y}~~\frac{|u(x)-u(y)|}{|x-y|^{\eta}}\leq C_{R_0}R^{-\eta}\Bigl[\Bigl(\fint_{B_R}|u|^{2^*}\Bigr)^{1/2^*}+R^{2-\frac np}\Vert V\Vert_{L^p(B_R)}\Bigr].
\]
\end{lemma}

\begin{lemma}\label{Harnack}\cite{CFG} (Scale-Invariant Harnack Inequality)  Assume $\mathbf{a}\equiv0$, $A$ is an elliptic matrix with real, bounded coefficients, and $V\in RH_{\frac n2}$. There exists a small constant $c=c(n,\lambda, \Lambda)$ such that whenever  $B=B(x_0, r)$, $r<\frac c{m(x_0,V)}$, $x_0\in \RR^n$, the following property holds.  For any  $u$ which solves $L_Eu=0$ in the weak sense on $2B$, 
\begin{equation}\label{harnack0}
\sup_{x\in B}u(x)\leq C\inf_{x\in B}u(x),
\end{equation}
with the constant $C>0$ depending on $n, \lambda, \Lambda$ and $ V$ only. 	
\end{lemma}

\begin{remark}\label{rHarnack} We remark that the Harnack inequality, of course, holds for any $0<r<r_0$, $r_0>0$, but typically with the constant growing exponentially in $r_0$ and $\|V\|_{L^p(B(x_0,r_0))}$. The important feature of \eqref{harnack0} is that the constant depends on $n, \lambda, \Lambda$ and $C_V$ only.
\end{remark}

%
%
%

The following lemma is a trivial extension of a particular result of Theorem 1.1 in \cite{GW}:


\begin{lemma}\label{fundamentallower} Suppose that $A$ is an elliptic matrix with real, bounded coefficients. Let $\Gamma_0$ be the unique fundamental solution to the operator $L_0=-\text{\emph{div} }A\nabla$. Then there exist constants $c_n, C_n$ greater than $0$ and depending on $n,\lambda,\Lambda$, such that
\begin{equation}\label{fundamentallower1}
\frac{c_n}{|x-y|^{n-2}}\leq\Gamma_0(x,y)\leq\frac{C_n}{|x-y|^{n-2}}.
\end{equation}
\end{lemma}

It will be useful to consider the conclusion of Lemma \ref{Harnack} in greater generality:

\begin{definition}\label{scinvharnack} We say that the operator $L$ satisfying assumptions (\ref{magassumptions-intro}) has the \emph{$m$-scale invariant Harnack Inequality} if  whenever  $B=B(x_0, r)$, $r<\frac c{m(x_0,V+|\mathbf{B}|)}$, $x_0\in \RR^n$, the following property holds.  For any  $u$ which solves $Lu=0$ in the weak sense on $2B$, 
\begin{equation}\label{harnack1}
\sup_{x\in B}|u(x)|\leq C\inf_{x\in B}|u(x)|,
\end{equation}
with the constant $C>0$ independent of $B$.
\end{definition}

\subsection{Proof of the lower bound estimate}

First, we establish two auxiliary propositions.

\begin{proposition}\label{represent} Let $\mathbf{a}\in L^2_{loc}(\bb R^n)$, and let $A$ be an elliptic matrix with complex, bounded coefficients. Assume that $V\in L^1_{loc}(\bb R^n)$ satisfies (\ref{imVreq}), (\ref{gkcond}) with $c_2\equiv c_4\equiv0$. Let
\[
L_V=-(\nabla-i\mathbf{a})^TA(\nabla-i\mathbf{a})+V,\qquad L_0=-(\nabla-i\mathbf{a})^TA(\nabla-i\mathbf{a}).
\]
Then for each $f\in(\dot{\m V}_{\mathbf{a},0})'$, $\dot L_V^{-1}f$ is well-defined, belongs to $\dot{\m V}_{{\mathbf a},V}$, and the identity
\begin{equation}\label{GG2}
\dot L_0^{-1}f=\dot L_V^{-1}f +\dot L_0^{-1} V\dot L_V^{-1}f
\end{equation}
holds in $\dot{\m V}_{\mathbf{a},0}$. Moreover, if $L_0,$ its adjoint $L_0^*$ and $L_V$ are operators whose inverses have fundamental solutions $\Gamma_0,\Gamma_0^*, \Gamma_V$ respectively which satisfy 
\begin{equation}\label{boundedfund}
\Gamma_V,\Gamma_{0},\Gamma_{0}^*\in L^{\infty}_{loc}(\bb R^n\times\bb R^n\backslash\{x=y\}),
\end{equation}
then the identity
\begin{equation}\label{uniquefund-1}
\Gamma_0(x,y)=\Gamma_V(x,y)+\int_{\bb R^n}\overline{\Gamma_0^*(z,x)}\Gamma_V(z,y)V(z)\,dz
\end{equation}
holds a.e. on $\bb R^n\times\bb R^n$.
\end{proposition}

\bp Let $f\in(\dot{\m V}_{\mathbf{a},0})'$. Since
\[
\Vert u\Vert_{\dot{\m V}_{\mathbf{a},0}}\leq C\Vert\D u\Vert_{L^2(\bb R^n)}\leq C\Vert u\Vert_{\dot{\m V}_{\mathbf{a},V}},
\]
it follows that $\dot{\m V}_{\mathbf{a},0}$ continuously embeds into $\dot{\m V}_{\mathbf{a},V}$, which implies in particular that $f\in(\dot{\m V}_{\mathbf{a},V})'$, and $\dot L_V^{-1}f\in\dot{\m V}_{\mathbf{a},V}\subset\dot{\m V}_{\mathbf{a},0}$. Therefore $\dot L_0\dot L_V^{-1}f\in(\dot{\m V}_{\mathbf{a},0})'$. Hence
\[
(\dot{\m V}_{\mathbf{a},0})'\ni~ f-\dot L_0\dot L_V^{-1}f=(\dot L_V-\dot L_0)\dot L_V^{-1}f=V\dot L_V^{-1}f,
\]
so that $\dot L_0^{-1}V\dot L_V^{-1}f\in\dot{\m V}_{\mathbf{a},0}$. Now, since
\[
\dot L_0\Big(\dot L_V^{-1}f +\dot L_0^{-1} V\dot L_V^{-1}f\Big)=\dot L_0\dot L_V^{-1}f+V\dot L_V^{-1}f=\dot L_V\dot L_V^{-1}f=f,
\]
by the invertibility of $\dot L_0$, we obtain (\ref{GG2}). Now let $\phi\in C_c^{\infty}(\bb R^n)$. Multiplying (\ref{GG2}) by $\overline{\phi}$ and integrating over $\bb R^n$ we observe
\begin{align}
\int_{\bb R^n}\Big[\dot L_0^{-1}f-\dot L_{V}^{-1}f\Big]\overline{\phi}&=\int_{\bb R^n}\Big((\dot L_0)^{-1}V\dot L_{V}^{-1} f\Big)\overline{\phi}\nonumber\\[4mm]&=\int_{\bb R^n}\Big(V\dot L_{V}^{-1} f\Big)\overline{\Big((\dot L_0^*)^{-1}\phi\Big)}.\label{GG4}
\end{align}
Now let $U_1,U_2\subset\bb R^n$ be open sets such that $\dist(U_1,U_2)>0$, and $f\in C_c^{\infty}(U_1),\phi\in C_c^{\infty}(U_2)$. Writing out the integral representations of the operators in (\ref{GG4}) and using Fubini's Theorem (justified since (\ref{boundedfund}), $V\in L^1_{loc}(\bb R^n)$, and $f,\phi$ are smooth, compactly supported, satisfying $\dist(\supp f,\supp g)>0$), we have that
\begin{multline}\label{ints}
\int_{\bb R^n}\int_{\bb R^n}\Big[\Gamma_0(x,y)-\Gamma_V(x,y)\Big]f(y)\overline{\phi(x)}\,dx\,dy\\[3mm]=\int_{\bb R^n}\int_{\bb R^n}\int_{\bb R^n}V(z)\Gamma_V(z,y)\overline{\Gamma_0^*(z,x)}f(y)\overline{\phi(x)}\,dx\,dy\,dz.
\end{multline}
We may rewrite (\ref{ints}) as
\begin{multline}\label{ints2}
\int_{\bb R^n}\int_{\bb R^n}\Big[\Gamma_0(x,y)-\Gamma_V(x,y)-\int_{\bb R^n}V(z)\Gamma_V(z,y)\overline{\Gamma_0^*(z,x)}\,dz\Big]f(y)\overline{\phi(z)}\,dy\,dx=0.
\end{multline}
Since (\ref{ints2}) holds for arbitrary $f,\phi\in C_c^{\infty}(\bb R^n)$ with disjoint supports, then (\ref{uniquefund-1}) holds for a.e. $x,y\in\bb R^n$.\hfill{$\square$}\\

\begin{lemma}\label{estimate1} Suppose that $\mathbf{a}\in L^2_{loc}(\bb R^n)$, $A$ is an elliptic matrix with complex, bounded coefficients, $V\in L^1_{loc}(\bb R^n)$, and that $L$, $L_0:=L-V$, $L_0^*$ are operators for which there exist fundamental solutions $\Gamma\equiv\Gamma_V$, $\Gamma_0$, $\Gamma_0^*$ in the sense of Definition \ref{fundsoln}. Assume that $L$ has the zero-source local uniform boundedness property and that $\Gamma_V,\Gamma_0,\Gamma_0^*$ satisfy (\ref{pointwiseupper1}). Moreover for $p>\frac n2$, if $\mathbf{a}\equiv0$, assume $V\in RH_{p}$; otherwise assume (\ref{magassumptions-intro}) with $\frac n2$ replaced by $p$.  Let $\tilde p=p$ if $n\geq4$ and $\tilde p\in\left(\frac32,\min\{3,p\}\right)$ if $n=3$, and denote $\delta:=2-\frac n{\tilde p}$. Then
\begin{equation}\label{ineq1}
|\Gamma(y,x)-\Gamma_0(y,x)|\leq\frac{C[|x-y|m(y,V+|\mathbf{B}|)]^{\delta}}{|x-y|^{n-2}},
\end{equation}
for a.e. $x,y\in \RR^n$ such that $|x-y|<\frac1{m(y,V+|{\mathbf B}|)}$. Here $C$ depends on $\Vert V+|\mathbf{B}|\Vert_{RH_{p}},p, \lambda,$ $\Lambda, n$ and the constants from (\ref{magassumptions-intro}).
\end{lemma}

\bp Since $\Gamma_V,\Gamma_0,\Gamma_0^*$ satisfy (\ref{pointwiseupper1}), then they also satisfy (\ref{boundedfund}). It follows from Proposition \ref{represent}, the upper bound of $\Gamma(x,y)$ in Theorem \ref{Upper Bound} and (\ref{pointwiseupper1}) that
\begin{equation}\label{ineq2}
|\Gamma(y,x)-\Gamma_0(y,x)|\leq C\int_{\bb R^n}\frac{e^{-\ep d(z,y,V+|\mathbf{B}|)}}{|z-x|^{n-2}|z-y|^{n-2}}V(z)\,dz\leq C(I_1+I_2+I_3),
\end{equation}
where $I_1, I_2,$ and $I_3$ denote the integrals over $B(x,r/2)$, $B(y,r/2)$, and
\[
\Omega:=\Big\{z\in\bb R^n:~|z-x|\geq r/2,~~|z-y|\geq r/2\Big\},
\]
respectively, with $r=|x-y|$. Recall that $|x-y|<\frac1{m(y,V+|\mathbf{B}|)}$. Write $R=\frac1{m(y,V+|\mathbf{B}|)}\sim\frac1{m(x,V+|\mathbf{B}|)}$. For $z\in B(x,r/2)$, we have $|z-y|>r/2$, so
\begin{align}
\nonumber I_1&\leq\frac C{r^{n-2}}\int_{B(x,r/2)}\frac{V(z)\,dz}{|z-x|^{n-2}}\leq\frac C{r^{n-2}}\Bigl[\frac 1{r^{n-2}}\int_{B(x,r/2)}V\Bigr]\\[4mm]&\leq\frac C{r^{n-2}}\Bigl(\frac rR\Bigr)^{2-\frac np}\Bigl[\frac 1{R^{n-2}}\int_{B(x,R)}V+|\mathbf{B}|\Bigr]\leq\frac C{r^{n-2}}\Bigl(\frac rR\Bigr)^{2-\frac np},\label{I1bound}
\end{align}
where the second inequality is due to (\ref{Kato}), the third one is due to Lemma \ref{Vprop}, and the last one is due to Proposition \ref{similar}.
Similarly we achieve
\begin{equation}\label{I2bound}
I_2\leq C\Bigl(\frac rR\Bigr)^{2-\frac np}\frac1{r^{n-2}}.
\end{equation}

To estimate $I_3$, we note that for all $z\in\Omega$,
\[
|z-y|\leq|z-x|+|x-y|=|z-x|+r\leq|z-x|+2|z-x|=3|z-x|,
\]
and so
\begin{align}
I_3&\leq C\int_{|z-y|\geq\frac r2}\frac{e^{-\ep d(z,y,V+|\mathbf{B}|)}V(z)\,dz}{|z-y|^{2n-4}}\nonumber\\[4mm]&\leq C\int_{2R>|z-y|\geq r/2}\frac{V(z)\,dz}{|z-y|^{2n-4}}~+~C\int_{|z-y|\geq2R}\frac{e^{-\ep d(z,y,V+|\mathbf{B}|)}V(z)\,dz}{|z-y|^{2n-4}}\nonumber\\[4mm]&=C(I_{31}+I_{32}).
\end{align}
To estimate $I_{31}$, we proceed similarly to the proof of Proposition \ref{intbound}. Recall $\tilde p$ from the statement of the Lemma. We note that $\tilde p \leq p$ in any dimension and hence
$V+|{\mathbf B}|\in RH_{\tilde p}(\bb R^n)$. Let $q$ be the H\"older conjugate of $\tilde p$. Since $\{z~|~|z-y|\in[r/2,2R]\}\subset B(y,2R)$, by H\"older's Inequality we have
\begin{align}\label{holderI31}
I_{31}\leq\Bigl(\int_{B(y,2R)}V^{\tilde p}\Bigr)^{\frac1{\tilde p}}\Bigl(~\int_{2R>|z-y|\geq r/2}\frac1{|z-y|^{2(n-2)q}}\,dz\Bigr)^{\frac1q}.
\end{align}
Next,
\begin{align}
\Bigl(\int_{B(y,2R)}V^{\tilde p}\Bigr)^{\frac1{\tilde p}}&=|B(y,2R)|^{\frac1{\tilde p}}\Bigl(\frac1{|B(y,2R)|}\int_{B(y,2R)}V^{\tilde p}\Bigr)^{\frac1{\tilde p}}\nonumber\\[4mm]&\leq\Vert V+|\mathbf{B}|\Vert_{RH_{\tilde p}}\big(C(n)\big)^{\frac1{\tilde p}-1}(2R)^{\frac n{\tilde p}-2}\Bigl(\frac1{(2R)^{n-2}}\int_{B(y,2R)}V+|\mathbf{B}|\Bigr).\nonumber
\end{align}

Hence, by the definition of $R$ and Proposition \ref{similar} we have
\begin{equation}\label{Vest}
\Bigl(\int_{B(y,2R)}V^{\tilde p}\Bigr)^{\frac1{\tilde p}}\leq CR^{\frac n{\tilde p}-2}.
\end{equation}

Now, observe that $n+\frac n{\tilde p}>4$ and hence $n-2q(n-2)<0$ for our choice of $\tilde p$. Therefore, 
\begin{equation}\label{estimating1}
\int_{2R>|z-y|\geq r/2}\frac1{|z-y|^{2(n-2)q}}\,dz\leq C(n,q)\,r^{n-2q(n-2)}.
\end{equation}

From (\ref{holderI31}), (\ref{Vest}), and (\ref{estimating1}), we obtain
$$
I_{31}\leq CR^{\frac n{\tilde p}-2}r^{\frac nq-2(n-2)}=C\Bigl(\frac rR\Bigr)^{2-\frac n{\tilde p}}\frac1{r^{n-2}}.
$$

As for $I_{32}$, first observe that $n+\frac{n}{\tilde p}>4$ can be written as
\begin{equation}\label{eqp}
n-2-\Bigl(2-\frac n{\tilde p}\Bigr)>0.
\end{equation}
With this in mind, we split up $\{|z-y|\geq2R\}$ into annuli of radius $2^jR$ for each $j\in\bb N$:
\[
I_{32}=\sum\limits_{j=1}^{\infty}\,\int_{2^jR\leq|z-y|\leq2^{j+1}R}\frac{e^{-\ep d(z,y,V+|\mathbf{B}|)}\,V(z)\,dz}{|z-y|^{2n-4}}\,.
\]
For $z\in\{2^jR\leq|z-y|\leq2^{j+1}R\}$, we observe that
\[
|z-y|m(y,V+|\mathbf{B}|)\geq2^jRm(y,V+|\mathbf{B}|)=2^j.
\]
Hence (\ref{claim}) implies that there exists $\ell>1$ such that
\[
d(z,y,V+|\mathbf{B}|)\geq c\ell^j,\qquad\text{for each } z\in\{2^jR\leq|z-y|\leq2^{j+1}R\}.
\]

Therefore,
\begin{align*}
I_{32}&\leq\sum\limits_{j=1}^{\infty}\frac{e^{-\ep c\ell^j}}{(2^jR)^{2n-4}}\int_{B(y,2^{j+1}R)}V\\[4mm]&\leq\sum\limits_{j=1}^{\infty}\frac{e^{-\ep c\ell^j}}{(2^jR)^{2n-4}}\,C_0^{j+1} \int_{B(y,R)}V+|\mathbf{B}|= \frac{C_0}{R^{n-2}}\sum\limits_{j=1}^{\infty}e^{-\ep c\ell^j}2^{-(2n-4)j}C_0^{j}\leq \frac{C}{R^{n-2}},
\end{align*}
where on the second inequality we invoked Proposition \ref{RHisdoubling} for $V+|\mathbf{B}|$, and $C_0$ is the constant in (\ref{doubling}). Here, $C$ depends on {$\lambda, \Lambda$}, $\ell$, $n$, and $C_0$, hence, on $\lambda,\Lambda, n$ and $\Vert V+|\mathbf{B}|\Vert_{RH_{\tilde p}}$ only. As a result,


\begin{align*}
I_{32}\leq\frac C{R^{n-2}}=C\Bigl(\frac rR\Bigr)^{2-\frac n{\tilde p}} \Bigl(\frac rR\Bigr)^{n-2-\Bigl(2-\frac n{\tilde p}\Bigr)}\frac1{r^{n-2}}\leq C\Bigl(\frac rR\Bigr)^{2-\frac n{\tilde p}}\frac1{r^{n-2}},
\end{align*}
using \eqref{eqp}.

Since $2-\frac n{\tilde p}\leq2-\frac np$ for any $n\geq3$ and $r<R$ by assumption, the above estimations on $I_1, I_2, I_3$ together with (\ref{ineq2}) imply
\[
|\Gamma(y,x)-\Gamma_0(y,x)|\leq C\Bigl(\frac rR\Bigr)^{2-\frac n{\tilde p}}\frac1{r^{n-2}}
\]
which translates to (\ref{ineq1}).\hfill{$\square$}

\begin{theorem}\label{lowerbounding} Suppose that $\mathbf{a}\in L^2_{loc}(\bb R^n)$, $A$ is an elliptic matrix with complex, bounded coefficients, $V\in L^1_{loc}(\bb R^n)$, and that $L$, $L_0:=L-V$, $L_0^*$ are operators for which there exist fundamental solutions $\Gamma\equiv\Gamma_V$, $\Gamma_0$ $\Gamma_0^*$ in the sense of Definition \ref{fundsoln}. Assume that $\Gamma_0$ satisfies (\ref{fundamentallower2}), that $\Gamma_0^*,\Gamma_V$ satisfy (\ref{pointwiseupper1}), that $L$ has the zero-source local uniform boundedness property, and that $L$ satisfies the $m$-scale invariant Harnack Inequality. Moreover, if $\mathbf{a}\equiv0$, assume that $V\in RH_{\frac n2}$; otherwise assume (\ref{magassumptions-intro}). Then there exist constants $c$ and $\ep_2$ depending on $\lambda,\Lambda, \Vert V+|\mathbf{B}|\Vert_{RH_{\frac n2}}, n$ and the constants from (\ref{magassumptions-intro}) such that
\begin{equation}\label{lowerbounded}
|\Gamma(x,y)|\geq\frac{ce^{-\ep_2d(x,y,V+|{\mathbf B}|)}}{|x-y|^{n-2}}.
\end{equation}
\end{theorem}

\bp Fix $x,y\in\bb R^n$. If $|x-y|m(x,V+|\mathbf{B}|)\leq c$ for $c$ small enough, then by Lemma \ref{estimate1} and (\ref{fundamentallower2}),
\begin{equation}\label{ineq3}
|\Gamma(y,x)|\geq|\Gamma_0(y,x)|-|\Gamma(y,x)-\Gamma_0(y,x)|\geq\frac{c'}{|x-y|^{n-2}}.
\end{equation}
Fix $\m C$ as the constant from Proposition \ref{Aprop}.
By the $m-$scale invariant Harnack Inequality (\ref{harnack1}), inequality (\ref{ineq3}) implies that for any $\m C\geq c$,
\begin{equation}\label{ineqwithA}
|\Gamma(y,x)|\geq\frac{\tilde c}{|x-y|^{n-2}},\qquad\text{if }~|x-y|m(x,V+|\mathbf{B}|)\leq\m C
\end{equation}
where $\tilde c$ depends on $\m C$, and hence for such $x,y$ the estimate (\ref{lowerbounded}) trivially  holds, since $\ep_2 d(x,y,V+|\mathbf{B}|)\geq0$. So it suffices to show (\ref{lowerbounded}) in the case that $x,y\in\bb R^n, x\neq y$ satisfy $|x-y|m(x,V+|\mathbf{B}|)>\m C$.
To this end, choose $\gamma:[0,1]\ra\bb R^n$ such that $\gamma(0)=x, \gamma(1)=y$, and
\[
\int_0^1m(\gamma(t),V+|\mathbf{B}|)|\gamma'(t)|\,dt\leq2d(x,y,V+|\mathbf{B}|)
\]
which can be done by the definition of $d$. Let
\[
t_0=\sup\Bigl\{t\in[0,1]~\Big|~|x-\gamma(t)|\leq\frac{\m C}{m(x,V+|\mathbf{B}|)}\Bigr\}
\]
Thus $t_0$ satisfies $|x-\gamma(t_0)|\leq\frac{\m C}{m(x,V+|\mathbf{B}|)}$, and hence by Lemma~\ref{propertiesm} we have $m(\gamma(t_0),V+|\mathbf{B}|)\sim m(x,V+|\mathbf{B}|)$. In the case that
\[
|y-\gamma(t_0)|\leq\frac1{m(\gamma(t_0),V+|\mathbf{B}|)}
\]
we note that
\[
|x-y|\leq|x-\gamma(t_0)|+|y-\gamma(t_0)|\leq\frac{\m C}{m(x,V+|\mathbf{B}|)}+\frac1{m(\gamma(t_0),V+|\mathbf{B}|)}\leq\frac{\m C+C}{m(x,V+|\mathbf{B}|)}
\]
which establishes (\ref{lowerbounded}) due to (\ref{ineqwithA}) (used with $A+C$ in place of $A$).

Now suppose that $|y-\gamma(t_0)|>\frac1{m(\gamma(t_0),V+|\mathbf{B}|)}$. Since $m(\cdot, V+|\mathbf{B}|)$ is locally bounded, we can define a finite sequence $t_0<t_1<\cdots<t_m\leq1$ such that 
\[
t_j=\inf\Bigl\{t\in[t_{j-1},1]~\Big|~|\gamma(t)-\gamma(t_{j-1})|\geq\frac1{m(\gamma(t_{j-1}),V+|\mathbf{B}|)}\Bigr\}.
\]
Respectively, it has the following properties
\begin{equation}\label{proptm1}
\gamma(t)\in B\Bigl(\gamma(t_{j-1}),\frac1{m(\gamma(t_{j-1}),V+|\mathbf{B}|)}\Bigr)\quad\text{for }~t\in[t_{j-1},t_j),~~j=1,2,\ldots,m,
\end{equation}
\begin{equation}\label{proptm2}
|\gamma(t_j)-\gamma(t_{j-1})|=\frac1{m(\gamma(t_{j-1}),V+|\mathbf{B}|)},\qquad j=1,\ldots,m
\end{equation}
and
\begin{equation}\label{proptm3}
\gamma(t)\in B\Bigl(\gamma(t_{m}),\frac1{m(\gamma(t_m),V+|\mathbf{B}|)}\Bigr)\quad\text{for }~t\in[t_m,1].
\end{equation}
It follows that
\begin{align*}
\int_0^1m(\gamma(t),V+|\mathbf{B}|)|\gamma'(t)|\,dt&\geq\sum\limits_{j=1}^{m}\int_{t_{j-1}}^{t_j}m(\gamma(t),V+|\mathbf{B}|)|\gamma'(t)|\,dt\\[3mm]&\geq c\sum\limits_{j=1}^{m-1}m(\gamma(t_{j-1}),V+|\mathbf{B}|)|\gamma(t_j)-\gamma(t_{j-1})|=cm
\end{align*}
{where the second inequality is due to (\ref{slow}) and the last equality is due to (\ref{proptm2}).} Hence
\begin{equation}\label{dm}
m\leq \frac{2}{c} \,d(x,y,V+|\mathbf{B}|).
\end{equation}
By definition of $t_0$, we have that
\[
|x-\gamma(t_j)|\geq\frac{\m C}{m(x,V+|\mathbf{B}|)}
\]
for each $j=0,\ldots,m$. It follows by (\ref{Aexists}) that
\[
x\notin B\Bigl(\gamma(t_j),\frac2{m(\gamma(t_j),V+|\mathbf{B}|)}\Bigr)
\]
for each $j=0,1,\ldots,m$. Let $u(z):=\Gamma(z,x)$. Then we have that $u$ is a weak solution of $Lu=0$ in $B\Bigl(\gamma(t_j),\frac2{m(\gamma(t_j),V+|\mathbf{B}|)}\Bigr)$. By (\ref{proptm2}) and the Harnack Inequality, we can deduce the following chain of inequalities:
\[
|u(\gamma(t_0))|\leq C|u(\gamma(t_1))|\leq\cdots\leq C^m|u(\gamma(t_m))|\leq C^{m+1}|u(y)|,
\]
with the constant $C$ depending on $n, \lambda, \Lambda$ and $\Vert V+|\mathbf{B}|\Vert_{RH_{\frac n2}}$ and the constants from (\ref{magassumptions-intro}) only. This implies
\[
|\Gamma(y,x)|\geq C^{-m-1}|\Gamma(\gamma(t_0),x)|\geq C^{-m-2}\frac1{|\gamma(t_0)-x|^{n-2}}\geq C^{-m-3}\Big(m(x,V+|\mathbf{B}|)\Big)^{n-2},
\]
possibly enlarging the value of $C$ which still depends on the same parameters.
Here,  the second inequality holds due to the fact that $\gamma(t_0)$ satisfies the hypothesis of (\ref{ineqwithA}). From (\ref{dm}) we deduce
\[
|\Gamma(y,x)|\geq C^{-\frac{2}{c}\,d(x,y,V+|\mathbf{B}|)}\Big(m(x,V+|\mathbf{B}|)\Big)^{n-2}.
\]

We can choose $\ep_2>0$ large enough so that
\[
\ep_2-\frac{2}{c}\,\ln C\geq0
\]
in which case we can write
\[
|\Gamma(y,x)|\geq Ce^{-\ep_2d(x,y,V+|\mathbf{B}|)}\Big(m(x,V+|\mathbf{B}|)\Big)^{n-2}.
\]

Finally, from the hypothesis that $m(x,V+|\mathbf{B}|)\geq\frac{\m C}{|x-y|^{n-2}}$ we obtain
\begin{equation}\label{finally}
|\Gamma(y,x)|\geq\frac{ce^{-\ep_2d(x,y,V+|\mathbf{B}|)}}{|x-y|^{n-2}}
\end{equation}
for a.e. $x,y\in\bb R^n$, $x\neq y$.

It is immediate that (\ref{finally}) implies (\ref{lowerbounded}). Indeed, fix $x,y\in\bb R^n$ with $x\neq y$.  Since the right-hand side of (\ref{finally}) is symmetric with respect to $x,y$, it follows that
\[
|\Gamma(x,y)|\geq\frac{ce^{-\ep_2d(y,x,V+|\mathbf{B}|)}}{|y-x|^{n-2}}=\frac{ce^{-\ep_2d(x,y,V+|\mathbf{B}|)}}{|x-y|^{n-2}},
\]
for almost every such $x,y$, as desired.\hfill{$\square$}\\

\begin{corollary}\label{lowerboundcor} Let $L_2$ be a generalized magnetic Schr\"odinger operator formally given by (\ref{generalized}) where $\mathbf{a}\equiv0$, $A$ is a real, bounded, elliptic matrix, and $V\in RH_{\frac n2}$. Then there exist constants $c$ and $\ep_2$ depending on $\lambda,\Lambda, n,$ and $\Vert V\Vert_{RH_{\frac n2}}$ such that its fundamental solution, $\Gamma_{L_2}$, satisfies
\begin{equation}\label{lowerbounded2}
\Gamma_{L_2}(x,y)\geq\frac{ce^{-\ep_2d(x,y,V)}}{|x-y|^{n-2}}.
\end{equation}
\end{corollary}
\noindent\emph{Proof.} Per the results given in Section \ref{electric}, the operator $L_2$ satisfies the hypothesis of Theorem \ref{lowerbounding}. The result follows.\hfill{$\square$}

\end{document}